\documentclass[preprint,sort&compress,12pt]{elsarticle}
\usepackage{bbm}
\usepackage{amsmath}
\usepackage{mathrsfs}
\usepackage{stmaryrd}
\usepackage{bm}
\usepackage{amsfonts}
\usepackage{amssymb}
\usepackage{amsthm}
\usepackage{lineno}
\usepackage[RGB]{xcolor}
\journal{\quad}
\newcommand{\mm}{\mathrm}
\newcommand{\ml}{\mathcal}

\newcommand{\dx}{\mathrm{d}x}

\newcommand{\be}{\begin{equation}}
\newcommand{\bea}{\begin{equation}\begin{aligned}}
\newcommand{\beas}{\begin{equation*}\begin{aligned}}
\newcommand{\eeas}{\end{aligned}\end{equation*}}
\newcommand{\eea}{\end{aligned}\end{equation}}
\newcommand{\ee}{\end{equation}}

\renewcommand{\div}{{\rm div }}

\begin{document}
\begin{frontmatter}
%
% \title{On the Inhibition of Parker Instability by \\
% a Non-slip Boundary Condition}
%%%%%
\title{On the Dynamical Stability and\\
 Instability of Parker Problem
 %\footnote{  We certify that the general content of the manuscript, in whole or in part, is not submitted, accepted, or published elsewhere, including conference proceedings. }
 }

\author[FJ,sd]{Fei Jiang}
\ead{jiangfei0591@163.com}
 \cortext[cor1]{Corresponding author.}
\author[sJ]{Song Jiang}
\ead{jiang@iapcm.ac.cn}
%%%%%%%%%%%%%%%%%%%%%%%%%%%%%%%
\address[FJ]{College of Mathematics and
Computer Science, Fuzhou University, Fuzhou 350108, China.}
\address[sd]{Key Laboratory of Operations Research and Control of Universities in Fujian, Fuzhou, 350108, China.}
\address[sJ]{Institute of Applied Physics and Computational Mathematics, P.O. Box 8009, Beijing 100088, China.}
%\address[ww]{Institute of Applied Mathematics, AMSS, Chinese Academy of Sciences, Beijing 100190, China.}

\begin{abstract}
{We investigate a {Parker problem} for the three-dimensional compressible isentropic viscous magnetohydrodynamic system
with zero resistivity in the presence of a modified gravitational force in a vertical strip domain in which
the velocity of the fluid is non-slip on the boundary, and focus on the stabilizing effect of the (equilibrium) magnetic field
through the non-slip boundary condition. We show that there is a discriminant $\Xi$, depending on the known physical parameters,
 for the stability/instability of the {Parker problem}. More precisely, if $\Xi > 0$, then the {Parker problem} is unstable,
i.e., the Parker instability occurs, while if $\Xi < 0$ and the initial perturbation satisfies some relations,
 then there exists a global (perturbation) solution which decays algebraically to zero in time, i.e., the Parker instability
 does not happen. The stability results in this paper reveal the stabilizing effect of the magnetic field
 through the non-slip boundary condition and the importance of boundary conditions upon the Parker instability,
 and demonstrate that a sufficiently strong magnetic field can prevent the Parker instability from occurring.
 In addition, based on the instability results, we further rigorously verify the Parker instability
 under Schwarzschild's or Tserkovnikov's instability conditions in the sense of Hadamard for a horizontally periodic domain.}
 \end{abstract}
\begin{keyword}
Compressible magnetohydrodynamic flow; Schwarzschild's criterion; Parker instability; magnetic buoyancy instability; stability.
%\MSC[2000] 35Q35\sep  76D03.
%(2000 is the default)
\end{keyword}
\end{frontmatter}

% \linenumbers

\newtheorem{thm}{Theorem}[section]
\newtheorem{lem}{Lemma}[section]
\newtheorem{pro}{Proposition}[section]
\newtheorem{cor}{Corollary}[section]
\newproof{pf}{Proof}
\newdefinition{rem}{Remark}[section]
\newtheorem{definition}{Definition}[section]
% \linenumbers

\section{Introduction}\label{introud}
\numberwithin{equation}{section}

The equilibrium, in which a gas layer in a gravitational field is supported in part by
a vertically decreasing horizontal magnetic field, is unstable \cite{HDGWAGTDX}, and such instability  is called the Parker instability (or the magnetic buoyancy instability in some literatures). The behavior of the Parker instability can be described as follows.
 Suppose that magnetic field lines are disturbed and begin to undulate. The mass in the raised portion of a loop drains down
 along the field lines so that the loop to become lighter than the ambient medium. If the buoyancy at the loop top is larger
 than the restoring magnetics tension, the loop rises further and the instability
 sets in \cite{SKMRFNA}. Parker noted that the falling mass accumulates in the magnetic valleys \cite{ParkerEN1955}, and thus
  explained that some of the ``large-scale" interstellar molecular cloud complexes can be formed in this way, since
the growth rate of the instability is 10 times larger than that of the Jeans's gravitational instability of the gas cloud.

Since Parker's pioneering work \cite{ParkerEN1955}, many physicists  have continued to develop the linear theory and nonlinear numerical
simulation of the Parker instability, see \cite{machida2013dynamo,HDGWAGTDXGKMH,KJRDJTWT,RLFSSGRSABPJFATP,MKADCRHDW}
and the references cited therein. Moreover, the physicists have considered the  Parker instability
to be also one of the possible mechanisms for other various astrophysical phenomenon \cite{RFMGATPA232018}, such as,
the rising and falling motions of gases above the spiral arm \cite{RFMGATPA,RFMGATPA2},　
the molecular loops in the Galactic center \cite{MLGFYS3},
the Barnard loop in Orion, the rise and emergence of magnetic flux tubes in the Sun and other stars
\cite{ParkerEN,RFMGATPA256,barker2012magnetic,RFMGATPA12345} as well as in accretion disks \cite{CGDLJGRADF},
the jets ejected from the centers of active galaxies \cite{PJHA14}, and so on. Recently, Khalzov et.al. first used
the Madison Plasma Couette Experiment to  model the Parker instability \cite{KIVBBPKNFCB}.

It has been also widely investigated how the Parker instability evolves under the effects of other physical factors,
such as rotation \cite{KJRDJTWT},
cosmic rays \cite{kudoh2014magnetohydrodynamic,kuwabara2004nonlinear}, corona \cite{KHHTMRHTSKMS},
self-gravity \cite{CWCMRTTUMSK}, nonuniform gravitational fields \cite{RFMGATPA232018}, random magnetic fields \cite{KJRDT}, and so on.
%It is well-known that the different boundary conditions can %affect
%the stability/instability of  flows, see the thermal instability %\cite{CSHHS}. Motivated by it,
The effect of boundary conditions of the velocity on the evolution of the Parker instability was investigated in our previous article \cite{JFJSJMFM}.
Based on the linearized  motion equations, Jiang et.al. have found a new phenomenon that the non-slip velocity boundary condition,
imposed on the direction of a (equilibrium) magnetic field, can enhance the stabilizing effect of the field,
so that the Parker instability can be prevented under a sufficiently strong magnetic field.
The main aims in this article are to rigorously prove this new phenomenon (i.e., the inhibition effect of a magnetic field on the nonlinear Parker instability through a non-slip boundary condition), and the criterion that gives the nonlinear Parker instability by developing new mathematical techniques
based on the nonlinear motion equations. Before stating our results, we formulate the problem mathematically.

\subsection{Parker problem}
 The verification of the stabilizing effect of a magnetic field on the Parker instability through a non-slip boundary
 condition and the investigation of the criteria leading to the occurrence of the Parker instability can be reduced to the proof
 of stability and instability for a {Parker problem} of the magnetohydrodynamic (MHD) equations, respectively.
 In this article, we consider the following three-dimensional (3D) compressible isentropic viscous MHD equations
 with zero resistivity (i.e., without magnetic diffusivity)
in the presence of a gravitational field in a domain $\Omega\subset \mathbb{R}^3$ read as follows (see, for example,
\cite{cowling1957magnetohydrodynamics,kulikovskiy1965magnetohydrodynamics} on the derivation of the motion equations).
\begin{equation}\label{comequations}\left\{\begin{array}{l}
 \rho_t+\mm{div}(\rho{  v})=0,\\[1mm]
\rho v_t+\rho v\cdot\nabla v+\nabla
   (P+  \lambda  |M|^2/2) =\mu_1\Delta v+\mu_2\nabla\mm{div} v
+\lambda M\cdot \nabla M-\rho g e_3,\\[1mm]
M_t=M\cdot \nabla v-v\cdot\nabla M - M \mm{div }v,\\[1mm]
\mathrm{div} {M}=0.\end{array}\right.\end{equation}
Here the unknowns $\rho:=\rho(x,t)$, ${v}:= {v}(x,t)$ and $M:= {M}(x,t)$  denote the density, velocity
and magnetic field  of the compressible MHD fluid, respectively; $\lambda$ stands for the permeability of vacuum dividing by $4\pi$,
$g>0$ for the gravitational constant,
$e_3=(0,0,1)^{\mm{T}}$ for the vertical unit vector, and $-\rho g {e}_3$ for the gravitational force. $\mu_1>0$ is the coefficient
of shear viscosity and $\mu_2:=\nu+\mu_1/3$ with $\nu$ being the positive bulk viscosity.
The pressure $P$ is usually determined through the equations of state. In this article we focus our study on the case of polytropic gas:
\begin{equation}\label{isetropiceest}
P\equiv P(\rho)=A\rho^\gamma,
\end{equation}
where $\gamma\geq 1$ denotes the adiabatic index and $A>0$ is a constant.
In the system \eqref{comequations} the equation \eqref{comequations}$_1$ is the continuity equation,
\eqref{comequations}$_2$ describes the balance law of momentum, while \eqref{comequations}$_3$ is called the induction equation.
As for the constraint $\mm{div}\, M=0$, it can be seen just as a restriction on the initial value of $M$ since
$(\mm{div}\,M)_t=0$ due to \eqref{comequations}$_3$. We remark that the resistivity is neglected in \eqref{comequations}$_3$,
 and this arises in the physics regime with negligible electrical resistance.

Next, we construct a (magnetohydrostatic) equilibrium state $(\bar{\rho},0,\bar{M})$ to the MHD equations \eqref{comequations}.
Firstly, we choose  a (equilibrium) density profile $\bar{\rho}:=\bar{\rho}(x_3)\in C^{1}(\bar{\Omega})$, which is independent of $(x_1,x_2)$ and satisfies
\begin{eqnarray}\label{0102}
 \inf_{ x\in \Omega}\bar{\rho}>0.
\end{eqnarray}
The condition \eqref{0102} prevents us from treating vacuum.
 Then, for given $\bar{\rho}$ and $g$, we defined a  horizontal magnetic  field profile $\bar{M} :=me_1:= (m ,0,0)^T$ with
\begin{equation}\label{com:01091}
   m\equiv m(x_3):=\pm\sqrt{\frac{2}{\lambda }\left(C
- \bar{P} -F(g\bar{\rho})\right)}, \quad \mbox{a function of }x_3\;\mbox{ only},
\end{equation}
where $\bar{P}:=P(\bar{\rho})$ is the pressure profile, $F(g\bar{\rho})$ denotes a primitive function $g\bar{\rho}$ and $C$ is a positive constant satisfying
$$\inf_{x\in \bar{\Omega}} \{C - \bar{P} -F(g\bar{\rho})\}>0.$$
It is easy to see that \eqref{com:01091} makes sense for a bounded domain $\Omega$, and
\begin{equation}\label{comsteady}
 \bar{P}'=-\lambda  m  m'  -g\bar{\rho},  \end{equation}
 where $':= {\mm{d}}/{\mm{d}x_3}$. Moreover,
\begin{equation}\label{equcomre}
 \nabla  (\bar{P}+\lambda   |\bar{M} |^2/2)=
 -g\bar{\rho} e_3\quad\mbox{and}\quad\mathrm{div} \bar{M} =0, \end{equation}
where  $\lambda \bar{M} \cdot \nabla \bar{M}=0$.
In what follows, we denote such  equilibrium-state $(\bar{\rho},0,\bar{M})$ by $s_e$.

It is well-known that the equilibrium-state  $s_e$
is unstable, if the density profile further satisfies Schwarzschild's (instability) condition \cite{HDGWAGTDX}
 \begin{eqnarray}\label{Newcombinstablity}
 -\bar{\rho}'(x_3^0)<\frac{g\bar{\rho}^2 }{\gamma \bar{P} }\bigg|_{x_3=x^0_3} \quad
 \mbox{ for some }x^0_{3}\in \{x_3~|~(x_1,x_2,x_3)^{\mm{T}}\in \Omega\}.
\end{eqnarray}
%Next we briefly introduce the  relevant  background of Schwarzschild's %condition.
Such condition was not first found in the Parker instability problem, but in the (compressible thermal) convection problem studied by Schwarzschild in 1906 \cite{KSNKGWG}.
We mention that, for the convection problem, the pressure profile in Schwarzschild's condition not only depends on a density profile,
but also on a temperature profile $\bar{T}$ due to the thermal effect, and often takes the form $\bar{P}=R\bar{\rho}\bar{T} $ with $R>0$ being a gas constant.
Later, Tserkovnikov further investigated the convection problem in the presence of a horizontal magnetic field (abbreviated as the magnetic convection (MC) problem),
and obtained an instability condition \cite{YATSHD5} in 1960:
 \begin{equation}\label{Tserkovnikovs}
-g\bar{\rho}'(x_3^0)<\frac{g^2\bar{\rho}^2}{\gamma  \bar{P}+\lambda m^2}\bigg|_{x_3=x^0_3}\quad\mbox{ for some }\; x_3^0\in (-l,l).
\end{equation}
However, in 1961, Newcomb \cite{NWACTOPF} extended  Tserkovnikov's analysis by imposing no constraints on the perturbation wave vectors and found
that Schwarzschild's condition is not only the sufficient and necessary condition for the compressible thermal instability, but also for the instability
in the MC problem, please refer to \cite{YCPCove} for the physical interpretation for Newcomb's and Schwarzschild's conditions.

After Parker's instability work in 1966, physicists further noted that the instability in the MC problem considered by Tserkovnikov and  Newcomb
not only involves the thermal instability, but also the Parker instability \cite{VGMBNHC,TJHNAHC,HDGWAGTDX}.
In particular, for the isentropic case (i.e., one omits the thermal effect and the pressure is of form \eqref{isetropiceest}), Schwarzschild's condition
is equivalent to the magnetic buoyancy condition
\begin{eqnarray}\label{Parker0102}  (m^2)'|_{x_3=x^0_{3}}<0.   \end{eqnarray}
by virtue of the relations \eqref{comsteady} and \eqref{isetropiceest}.
In other words, under Schwarzschild's condition, there exists a region, in which the strength of $\bar{M}$
is vertically decreasing. Thus the corresponding magnetic pressure in \eqref{equcomre} plays a role of the magnetic buoyancy \cite{ParkerEN},
which is able to support more mass against gravity than would be possible in its absence. In view of \eqref{Newcombinstablity} and \eqref{Parker0102},
we can see that the mechanism of the Parker instability refers to the pressure state and the magnetic buoyancy.
To emphasis the mechanism of the magnetic buoyancy, the Parker instability is often called the magnetic buoyancy instability
(or the ballooning instability in fusion plasma physics \cite{SKTTMRSP}).

Finally, we mention a special density profile, i.e., $\bar{\rho}$ satisfies  the Rayleigh-Taylor condition
\begin{eqnarray}\label{01022107}
 \bar{\rho}'|_{x_3=x^0_{3}}>0    \end{eqnarray}
for some $x^0_{3}\in \{x_3~|~(x_1,x_2,x_3)\in \Omega\}$. By \eqref{Newcombinstablity}, the equilibrium-state $s_e$ satisfying \eqref{01022107} is unstable. Since such density
distribution is in close correspondence to the classical case of a heavy fluid supported by a light one, the Parker instability under such case is called
the magnetic Rayleigh-Taylor instability \cite{IHMTSKYTFstru}, or the Kruskal-Schwarzschild instability due to the first investigation of Kruskal and Schwarzschild in 1953 \cite{KMSMSP},
where they further pointed out that curvature of {the magnetic lines} can influence the development of instability.

Now, we introduce the {Parker problem} for the MHD equations around the equilibrium state $s_e$.  Denoting the perturbation to the equilibrium state  $s_e$ by
$$ \varrho=\rho -\bar{\rho},\quad  v= v- {0},\quad N=M-\bar{M} ,$$
 and using the relations in  \eqref{equcomre}, we obtain the perturbation equations:
\begin{equation}\label{0103} \left\{\begin{array}{l}
\varrho_t+\mm{div}((\varrho+\bar{\rho})v)=0, \\[1mm]
(\varrho+\bar{\rho}){v}_t+(\varrho+\bar{\rho}){v}\cdot\nabla v+\nabla (P(\varrho+\bar{\rho})+\lambda |N+\bar{M} |^2/2)\\
\quad =\mu_1\Delta v+\mu_2\nabla\mm{div} v + \lambda ( N+\bar{M} )\cdot \nabla (N+\bar{M} )-(\varrho+\bar{\rho}) g e_3,\\[1mm]
 N_t=( N+\bar{M} )\cdot \nabla v - v\cdot\nabla (N+\bar{M} )-(N+\bar{M} )\mm{div}v,\\[1mm]
 \mathrm{div} N=0.  \end{array}\right.  \end{equation}
We impose the following initial and boundary conditions for \eqref{0103}:
\begin{eqnarray} \label{c0104}   &&
(\varrho,v, {N} )|_{t=0}=(\varrho_0,v_0,N_0)\;\quad\mbox{in } \Omega ,  \\
&& \label{0105}
v|_{\partial\Omega}=0 \;\quad \mbox{ for any }\; t>0.
\end{eqnarray}
%The initial-boundary problem  \eqref{0103}--\eqref{0105} is called %the Parker problem.　
 The linearized equations of  \eqref{0103} around the  equilibrium state $s_e$ read as
\begin{equation}\label{lincom}
\left\{\begin{array}{ll}
 \varrho_t+\mm{div}( \bar{\rho}v)=0, \\[1mm]
  \bar{\rho}v_t +\nabla (P'(\bar{\rho})\varrho+ \lambda
 m N_1)=\mu_1\Delta v+\mu_2\nabla\mm{div}v + \lambda\bar{M}' N_3  + \lambda m \partial_1  N-\varrho  ge_3,\\[1mm]
   N_t=m\partial_1 v -v_3 \bar{M}'-\bar{M}\mm{div}v.\\[1mm]
 \mathrm{div} N=0, \end{array}\right.\end{equation}
 where $P'(\bar{\rho}):=(A x^\gamma_3)'|_{x_3=\bar{\rho}}$, and $N_i$ and $v_i$ denote the $i$-th component of $N$ and $v$, respectively.
The system \eqref{lincom} with initial and boundary conditions \eqref{c0104}--\eqref{0105} constitutes a linearized {Parker problem},
while the initial-boundary problem \eqref{0103}--\eqref{0105} is called the nonlinear {Parker problem}.
The Parker instability, which is shown based on the linearized {Parker problem}, resp. the nonlinear {Parker problem},
is called the linear Parker instability, resp. the nonlinear Parker instability.
At present, only the linear Parker instability under Schwarzschild's condition or Tserkovnikov's condition \eqref{Tserkovnikovs} is mathematically
investigated, also for the case without viscosity.

\subsection{Criterion  for stability/instability}
The linearized equations are convenient to analyze mathematically in order to have an insight into the physical
and mathematical mechanisms of the Parker instability. Moreover, applying the energy principle \cite{BIBFEAKMDKRM} to
the linearized {Parker problem}  \eqref{c0104}--\eqref{lincom}, Jiang et.al. have obtained criteria of stability/instability for the linearized {Parker problem} \cite{JFJSJMFM}.
More precisely, in the case of a bounded domain $\Omega$, one has
 \begin{enumerate}[\quad  \  (1)]
  \item[(1)] if $\Xi<0$, then the linearized {Parker problem} is stable;
  \item[(2)]if $\Xi>0$, then the linearized {Parker problem} is unstable,
\end{enumerate}
where
$$ \Xi:=\sup_{w\in H_0^1(\Omega)}\frac{E (w)}{\int_\Omega\lambda (|\partial_{1}w_{\mm{v}}|^2 +|\mm{div}_{\mm{v}}{w}_{\mm{v}}|^2)\mm{d}x},$$
and
\begin{eqnarray} && \label{20160306cond}
E(w):=\int_\Omega g\left(\bar{\rho}'+\frac{g \bar{\rho}^2}{\gamma\bar{P}}\right)w_3^2\mm{d} x -\int_\Omega \frac{1}{\gamma\bar{P}}
\left( g\bar{\rho}w_3 - {\gamma \bar{P}}\mm{div}w\right)^2\dx   \nonumber \\
&& \qquad \qquad -\int_\Omega  \lambda m^2 (|\partial_{1}w_{\mm{v}}|^2 +|\mm{div}_{\mm{v}}w_{\mm{v}}|^2)\dx ,
\end{eqnarray}
here and in what follows $w_{\mm{v}}:=(w_2,w_3)$ and $\mm{div}_{\mm{v}}w_{\mm{v}}:=\partial_2 w_2+\partial_3 w_3$ for $w=(w_1,w_2,w_3)$.

We mention that by virtue of \eqref{comsteady}, $E(w)$ can be rewritten as
\begin{equation} \label{05031703new}
\begin{aligned}
E(w)=&-\int_\Omega \frac{\lambda g\bar{\rho} mm'w_3^2}{ \gamma \bar{P}}\mm{d} x-\int_\Omega \frac{1}{\gamma\bar{P}}\left( g\bar{\rho}w_3 - {\gamma \bar{P}}\mm{div}w\right)^2\dx \\
& -\int_\Omega  \lambda m^2 (|\partial_{1}w_{\mm{v}}|^2 +|\mm{div}_{\mm{v}}w_{\mm{v}}|^2  ) \dx.
\end{aligned}
\end{equation}

 {In view of the energy functionals \eqref{20160306cond} and \eqref{05031703new}, we can see that Schwarzschild's condition or the magnetic buoyancy condition
may make $E(w)$ to be positive for some $w$, thus contributing to the occurrence of the Parker instability.
In particular, if $\Omega$ is an infinite layer domain, Newcomb in 1961 found that Schwarzschild's condition leads to the linear Parker instability in some sense
by using the Fourier analysis method \cite{NWACTOPF}. However, in the case of a bounded domain, Jiang et.al. \cite{JFJSJMFM} found the stabilizing effect of a strong
(equilibrium) magnetic field $\bar{M}$ through a non-slip boundary condition upon the Parker instability, even if
the density profile $\bar{\rho}$ satisfies the Schwarzschild's condition. In fact, if we define
$$\bar{n}:=\sup_{w\in H_0^1(\Omega)} \sqrt{\frac{\int_\Omega g(\bar{\rho}'+ g\bar{\rho}^2/\gamma \bar{P})w_3^2\mm{d} x}{\int_\Omega\lambda
(|\partial_{1}w_{\mm{v}} |^2 + |\mm{div}_{\mm{v}}w_{\mm{v}}|^2) \mm{d}x}},$$
then $\bar{n}\in [0,+\infty)$, since $\Omega$ is bounded in the $x_1$-direction.
In view of \eqref{com:01091}, for given $\bar P$, $\bar{\rho}$, $g$ and $\lambda $, we can choose $m$ satisfying
$\inf_{x\in\Omega}|m | >\bar{n}$. Then, it is easy to verify that such $m$ satisfies $\Xi<0$. This infers that, under the non-slip boundary condition,
a sufficiently strong $\bar{M}$ has a remarkable stabilizing effect to prevent the Parker instability from occurring.
 It is well-known that a magnetic field has the stabilizing effect upon the Parker instability, but can not prevent the Parker instability from occurring
in the case of an infinite layer domain. However, the non-slip velocity boundary condition in a bounded domain,
imposed in the direction of the magnetic field, can enhance the stabilizing effect of the magnetic field,
so that Schwarzschild's criterion \eqref{Newcombinstablity} of the Parker instability fails.

In this article, we further extend the above linear results to the nonlinear case by developing new mathematical techniques.
In other words, we will rigorously verify that $\Xi$ is also the discriminant for instability/stability of the nonlinear
{Parker problem} \eqref{0103}--\eqref{0105} under some additional conditions. Moreover,
for a horizontally periodic domain, we shall provide a rigorous mathematical proof of the nonlinear Parker instability under Schwarzschild's
or Tserkovnikov's condition in the sense of Hadamard. The detailed results of nonlinear stability/instability will be presented in
Section \ref{0608241825}.

We end this section by deriving a upper-bound for $\bar{n}$, which may be useful in experimental researches and numerical simulations.
From the definition of $\bar{n}$ we see that $\bar{n} \leq \varpi \chi$, where
 $$\varpi:=\sqrt{{\|g\bar{\rho}'+g^2\bar{\rho}^2/\gamma \bar{P} \|_{L^\infty(\Omega)}}/\lambda}\;\;\mbox{ and }\;\;
 \chi:=\sqrt{\frac{\int_\Omega w_3^2\mm{d}x}{\int_\Omega |\partial_{1}w_{\mm{v}} |^2\mm{d}x}}. $$
  Let $a:=\inf\{x_1~|~(x_1,x_{\mm{v}})\in \Omega\}$, $b:=\sup\{x_1~|~(x_1,x_{\mm{v}})\in \Omega\}$  and
  $$ \Omega':=\{(x_1,x_\mm{v})\in \mathbb{R}^3~|~x_{\mm{v}}\in \mathcal{T},\ a<x_1<b\},  $$
  where  $\mathcal{T} :=(2\pi L_1\mathbb{T})\times(2\pi L_2\mathbb{T}) $, $\mathbb{T}=\mathbb{R}/\mathbb{Z}$, and $2\pi L_1$, $2\pi L_2>0$
  are the periodicity lengths. Then, one can choose sufficiently large $L_1$ and $L_2$, such that
  $H_0^1(\Omega)$ can be regarded as a subspace of $H(\Omega')$ by horizontally periodic translation. Thus one has
  %%%%%%%%%%%%%%%%%%%%%%%%%%%%%%%%%%%%%%%%%%%%%%%%%
\begin{equation}\label{05022101}
\bar{n}\leq \varpi\sup_{w\in H_0^1(\Omega')}\sqrt{\frac{\int_\Omega w_3^2\mm{d}x}{\int_\Omega |\partial_{1}w_{\mm{v}} |^2\mm{d}x}}.
\end{equation}
  On the other hand, similarly to the derivation of \cite[Proposition 5.1]{JFJSJMFM}, we use the Fourier analysis method to infer that
\begin{equation}\label{05022103}
\sup_{w\in H_0^1(\Omega')}\sqrt{\frac{\int_\Omega w_3^2\mm{d}x}{\int_\Omega  |\partial_{1}w_{\mm{v}} |^2\mm{d}x}}
  =\sup_{\psi\in H_0^1(a,b)}\sqrt{\frac{\int_a^b \psi^2\mm{d}x_1}{\int_a^b |\partial_{1}\psi |^2\mm{d}x_1}},     \end{equation}
where it is easy to see that
\begin{equation}  \label{05022102}
\sup_{\psi\in H_0^1(a,b)}\sqrt{\frac{\int_a^b \psi^2\mm{d}x_1}{\int_a^b |\partial_{1}\psi |^2\mm{d}x_1}}=\frac{b-a}{\pi}.  \end{equation}
Consequently, we can deduce from \eqref{05022101}--\eqref{05022102} that
$\bar{n}\leq (b-a)\varpi/\pi$. This means that ${\Xi}<0$ for $\inf_{x\in\Omega}|m | > {(b-a)\varpi}/{\pi}$.

%%%%%%%%%%%%%%%%%%%%%%%%%%%%%%%%%%%%%
%%%%%%%%%%%%%%%%%%%%%%%%%%%%%%%%%%%%%
\section{Main results}\label{0608241825}
In this section we state the main results of this paper on instability/stability of the nonlinear
{Parker problem} \eqref{0103}--\eqref{0105}.

\subsection{Reformulation of the nonlinear stability}
In general, it is difficult to directly show the existence of a unique global-in-time solution to the {Parker problem} \eqref{0103}--\eqref{0105}
defined on a general bounded domain when $\Xi<0$, since the magnetic field is difficult to control. To circumvent such difficulty,
similarly to \cite{JFJSWWWOA,TZWYJGw}, we switch our analysis to that in Lagrangian coordinates.
We mention that such transformation method have been also used in the proof of global well-posedness of incompressible
MHD equations, please refer to \cite{XLPZ10221525,FHLXLPZ10221526,FHLPZ10221525,HABPZOT10221525}.

 To show the stability of the nonlinear {Parker problem} in Lagrangian coordinates, we assume that the domain is a vertical strip, i.e.,
 \begin{equation}  \label{02011959101604091846}
\Omega:=\{(x_1,x_2,x_3)\in \mathbb{R}^3~|~x_1\in (0,l)\}\;\mbox{ with }l >0.  \end{equation}
To make the expression \eqref{com:01091} sense, in this article we modify the gravitational constant $g$ to be
\begin{equation}   \label{modififorg}
0\leq g:=g(x_3)\in C^6_0(\mathbb{R}) \end{equation}
and suppose that
\begin{equation}
\label{rho20160326}
\bar{\rho}\in B^7(\mathbb{R}):=\left\{f\in C^7(\mathbb{R})~\left|~\sup_{x_3\in \mathbb{R}}
\left|\frac{\mm{d}^if}{\mm{d}x^i}\right|<\infty\;\mbox{ for }\; 0\leq i\leq 7\right.\right\}.
\end{equation}
Under the conditions  \eqref{02011959101604091846}--\eqref{rho20160326} and \eqref{0102},
 there always exists a $m$ constructed by  \eqref{com:01091}. Moreover, $m$ satisfies
\begin{equation}
\label{rho20160326mB6}m\in B^7(\mathbb{R}).
\end{equation}
%%%
We further assume that there is an invertible mapping $\zeta_0:=\zeta_0(y):\Omega\to \Omega$, such that
\begin{equation}  \label{zeta0inta}
(\zeta_0-y)|_{\partial\Omega}=0\mbox{ and }\det(\nabla \zeta_0(y))\neq 0\mbox{ for any }y\in \bar{\Omega}.
\end{equation}
Then, defining the flow map $\zeta$ as the solution to
\begin{equation*}\label{0114}
\left\{       \begin{array}{l}
 \zeta_t (y,t)=v(\zeta(y,t),t) \\
\zeta (y,0)= {\zeta_0},   \end{array}    \right.
\end{equation*}
we denote the Eulerian coordinates by $(x,t)$ with $x=\zeta(y,t)$, whereas $(y,t)\in \Omega\times\mathbb{R}^+$ stand for the
Lagrangian coordinates. In order to switch back and forth from Lagrangian to Eulerian coordinates, we assume that
$\zeta(\cdot,t)$ is invertible and $\Omega=\zeta(\Omega, t)$.

We define now the unknowns in Lagrangian coordinates by
\begin{equation*}\label{0115}
(\sigma, u, B)(y,t)=(\rho,v, M)(\zeta(y,t),t),\qquad (y,t)\in \Omega\times\mathbb{R}^+.
\end{equation*}
Thus, the evolution equations for $\sigma$, $u$ and $B$ in Lagrangian coordinates read as
\begin{equation}\label{0116}
\left\{       \begin{array}{ll}
\zeta_t=u, \\[1mm]
\sigma_t +\sigma\mm{div}_{\mathcal{A}}u=0,\\
\sigma u_t -\mu_1 \Delta_{\ml{A}}u-\mu_2\nabla_{\mathcal{A}}\mm{div}_{\mathcal{A}}u
+\nabla_{\ml{A}}( {P}(\sigma)+\lambda |B|^2/2) = \lambda  B\cdot \nabla_{\ml{A}} B-\sigma \tilde{g}e_3, \\[1mm]
B_t-B\cdot\nabla_{\ml{A}}u+B\mm{div}_{\mathcal{A}}u=0, \\[1mm]
\div_\ml{A}B=0
\end{array}     \right.
\end{equation}
with initial and boundary conditions
$$(u,\zeta-y)|_{\partial\Omega}=0\;\mbox{ and }\; (\zeta,\sigma, u,B)|_{t=0}=(\zeta_0,\sigma_0,u_0,B_0). $$
Here $\tilde{g}:=g(\zeta_3)$, the matrix $\mathcal{A}:=(\ml{A}_{ij})_{3\times 3}$ via
$\ml{A}^T=(\nabla \zeta)^{-1}:=(\partial_j \zeta_i)^{-1}_{3\times 3}$, and
the differential operators $\nabla_{\ml{A}}$, $\mm{div}\ml{A}$ and $\Delta_{\ml{A}}$ are defined by $\nabla_{\ml{A}}f:=(\ml{A}_{1k}\partial_kf,
\ml{A}_{2k}\partial_kf,\ml{A}_{3k}\partial_kf)^T$, $\mm{div}_{\ml{A}}(X_1,X_2,X_3)^T:=\ml{A}_{lk}\partial_k X_l$ and $\Delta_{\mathcal{A}}f:=\mm{div}_{\ml{A}}\nabla_{\ml{A}}f$.
%%%%  for appropriate $f$ and $X$.
It should be noted that  we have used the Einstein summation convention over repeated indices, and $\partial_k=\partial_{y_k}$.
Additionally, in view of the definition of $\mathcal{A}$, one can deduce the following two important properties:
\begin{eqnarray} &&
\label{AklJ=0} \partial_l (J\mathcal{A}_{kl})=0, \\
&& \label{AklJdeta}
  \partial_i\zeta_k\ml{A}_{kj}= \ml{A}_{ik}\partial_k\zeta_j=\delta_{ij},
 \end{eqnarray}
 where $J=\det (\nabla \zeta)\neq 0$, $\delta_{ij}=0$ for $i\neq j$ and $\delta_{ij}=1$ for $i=j$. The relation
  \eqref{AklJ=0} is often called the geometric identity.  In addition, it is easy to check that, by the boundary conditions of $\zeta$ and $\zeta_0$,
   \begin{equation}\label{201810301047}
  J\mathcal{A}e_1|_{\partial\Omega}=J_0\mathcal{A}_0e_1|_{\partial\Omega}=e_1.
  \end{equation}
%we can see that $\mathcal{A}=(A^*_{ij})_{3\times 3}$, where
%$A^{*}_{ij}$ is the algebraic complement minor of $(i,j)$-th entry $\partial_j %\zeta_i$. Since $ \partial_kA^{*}_{ik}=0$,
% Moreover, we can use \eqref{} to wr \begin{equation}\label{diverelation} \mm{div}_{\mathcal{A}}u=\partial_l (\mathcal{A}_{kl}u_k)=0,
%\end{equation} which will be used in the derivation of temporal derivative %estimates. We can further deduce from \eqref{} the geometric identity

Our next goal is to eliminate $(\sigma,B)$ in \eqref{0116} by expressing them in terms of $\zeta$, and this can be
achieved in the same manner as in \cite{JFJSWWWOA,TZWYJGw,WYTIVNMI}. We mention that this idea was also used in the proof of the global well-posedness for the Cauchy problem
of incompressible or compressible MHD fluids without magnetic diffusivity \cite{HXPGETDCMF,HXPGETDCMFLin,HABPZOT10221525}; please refer to \cite{RHPYZYZHU,XXRZYXIAngZFZhang,RXXWJHXZYZZF,FHLPZ10221525} for other relevant results of global well-posedness.
For the reader's convenience, we give the derivation here. It follows from \eqref{0116}$_1$ that
\begin{equation}\label{Jtdimau}
J_t=J\mm{div}_{\mathcal{A}}u,             \end{equation}
which, together with \eqref{0116}$_2$, yields that
\begin{equation}\label{0122sd}\partial_t(\sigma J)=0.   \end{equation}
Applying $J\ml{A}_{jl}$ to \eqref{0116}$_4$, we can use \eqref{AklJdeta} and \eqref{Jtdimau} to infer that
\begin{equation*}     \begin{aligned}
J\ml{A}_{jl}\partial_t B_j=&J B_i\ml{A}_{ik}(\partial_t \partial_k \zeta_j ) \ml{A}_{jl}- \mathcal{A}_{jl} B_j J\mm{div}_{\ml{A}}u  \\
 =&-JB_i\mathcal{A}_{ik}\partial_k \zeta_j\partial_t \ml{A}_{jl}-J_t \ml{A}_{jl}B_j =-JB_j\partial_t\ml{A}_{jl}-J_t \ml{A}_{jl}B_j,
\end{aligned}           \end{equation*}
which implies $\partial_t(J\ml{A}_{jl}B_j)=0$, i.e.,
\begin{equation}\label{JAN=0}
\partial_t(J\mathcal{A}^T B) =0.
\end{equation}
In addition, applying $\mm{div}$-operator to the above identity and using the geometric identity, we obtain
\begin{equation}\label{tJdivnes}\partial_t(J\mm{div}_\mathcal{A} B)=\partial_t\mm{div}(J\mathcal{A}^T B)=0.\end{equation}

To obtain the time-asymptotical stability of the  equilibrium state, we naturally expect
$$ (\zeta,\sigma,u,B)(t)\to (y,\bar{\rho},0, \bar{M}{})\mbox{ as }t\to \infty. $$
Hence, from \eqref{AklJ=0} and \eqref{0122sd}--\eqref{tJdivnes} it follows that
\begin{equation*}   \label{imporeation}
\begin{aligned}  & J\mm{div}_\mathcal{A} B=J_0\mm{div}_{\mathcal{A}_0} B_0=\mm{div} \bar{M}{}=0,\\
&\sigma J= \sigma_0J_0=\bar{\rho}\mbox{ and } J\mathcal{A}^T B=J_0\mathcal{A}^T_0 B_0= \bar{M} ,
\end{aligned}    \end{equation*}
which implies that
\begin{align}   &  \label{abjlj0}
\mm{div}_\mathcal{A} B=0,\;\;\sigma=\bar{\rho}J^{-1}\;\mbox{ and }\; B= m{} J^{-1}\partial_1\zeta ,
\end{align}
provided the initial data $(\sigma_0,B_0,\zeta_0)$ satisfy
\begin{align}
&\label{abjlj0i}\sigma_0=\bar{\rho}J^{-1}_0\mbox{ and } B_0= m J^{-1}_0\partial_1\zeta_0.
\end{align}
Here $\mathcal{A}_0$ and $J_0$ denote the initial data of $\mathcal{A}$ and $J$, respectively.
We mention that, by \eqref{201810301047} and the fact  $J\mathcal{A}^T B=J_0\mathcal{A}^T_0 B_0$, we have $ B_1 =B_1^0$ on $\partial\Omega$,
 where $B^0_1$ is the first component of $B_0$.

Let
\begin{equation}  \label{relation201604091949}
\eta:=\zeta-y,\mbox{ i.e., }\zeta=\eta + y. \end{equation}
Next we use $\eta$ to represent the (generalized) Lorentz force, the pressure term and the modified gravity term. By a straightforward computation,
we can split the Lorentz force into
\begin{equation}\label{0131}
\begin{aligned}&B\cdot \nabla_{\ml{A}} B -\nabla_{\mathcal{A}} |B|^2/2
%=B\cdot \nabla_{\mathcal{A}}B-\nabla_{\mathcal{A}}B_kB_k
=Lo_1+\mathcal{N}_1  - \nabla_{\mathcal{A}} |\tilde{M}{}|^2/2,
\end{aligned}\end{equation}
where $\tilde{M}:=\bar{M}(\zeta_3(y,t))$,
\begin{equation*}\label{0131def}  \begin{aligned}
Lo_1:= & \bar{M} \cdot \nabla (B-\bar{M}{}) + (B- \bar{M}{})\cdot \nabla \bar{M}{} - m  m 'J^{-1}\partial_1\eta_3    e_1\\
& -(\bar{M} ^T \nabla (B -\bar{M} ) +   (B -\bar{M} )^T \nabla \bar{M})^{T} +\nabla (m m '\eta_3)
\end{aligned}\end{equation*}
and
\begin{equation*}\label{0131defd}
\begin{aligned}  \mathcal{N}_1:=&
 \bar{M} \cdot (\nabla_{\mathcal{A}}-\nabla )(B-\bar{M} )+(B-\bar{M} )\cdot \nabla_{\mathcal{A}}(B-\bar{M} )\\
&  + (B- \bar{M} )\cdot (\nabla_{\mathcal{A}}-\nabla) \bar{M}  -(\bar{M} ^T(\nabla_{\mathcal{A}}-\nabla )(B -\bar{M} )\\
& + (B -\bar{M} )^T\nabla_{\mathcal{A}} (B-\bar{M} ) + (B -\bar{M} )^T(\nabla_{\mathcal{A}}-\nabla)\bar{M})^T  \\
& + m m' J^{-1}( \partial_1\eta_2\partial_2\eta_3- \partial_1\eta_3 \partial_2\eta_2  )e_1 \\
&+(\nabla_{\mathcal{A}}-\nabla)( |\tilde{M} |^2- |\bar{M} |^2)/2+ \nabla  ( |\tilde{M} |^2- |\bar{M} |^2-2m m '\eta_3 )/2.
\end{aligned}\end{equation*}
On the other hand, by the expression of $B$ in \eqref{abjlj0} and $\bar{M}$, one has
\begin{equation}  \label{20160309MH}
\begin{aligned}
 Lo_1 &=  m \partial_1 (B-m e_1)+  B_3  m 'e_1 - J^{-1} m  m '\partial_1\eta_3 e_1 +\nabla ( (m m '\eta_3 -m  (B_1-m )) \\
&= m^2\partial_1[J^{-1}\partial_1\eta+(J^{-1}-1)e_1] +\nabla( m m '\eta_3 -m^2 (J^{-1}\partial_1\eta_1+J^{-1}-1)) \\
& =\lambda^{-1}\mathcal{L}_M+\mathcal{N}_2,
\end{aligned}   \end{equation}
where
\begin{eqnarray*}
\mathcal{L}_M &=&\lambda(m^2\partial_1^2\eta -m^2\partial_1\mm{div}\eta e_1 +\nabla (m m '\eta_3+m^2   \mm{div}_\mm{v}\eta_\mm{v})), \\
\mathcal{N}_2 &=& m^2\partial_1((J^{-1}-1)\partial_1\eta + (J^{-1}-1+\mm{div}\eta)e_1) \\
&& \quad -\nabla( m^2( (J^{-1}-1)\partial_1\eta_1+J^{-1}-1+\mm{div}\eta)).
\end{eqnarray*}
Thus, inserting \eqref{20160309MH} into \eqref{0131}, we get
\begin{equation}  \label{leolasexipre}
B\cdot \nabla_{\ml{A}} B -\nabla_{\mathcal{A}}(|B|^2/2)= \lambda^{-1}\mathcal{L}_M - \nabla_{\mathcal{A}} |\tilde{M} |^2/2+\mathcal{N}_1+\mathcal{N}_2.
\end{equation}

Now we turn to dealing with the pressure term. It is easy to see that
$$P(\sigma)=P(\bar{\rho}J^{-1})=\bar{P}
+P'(\bar{\rho})\bar{\rho}(J^{-1}-1)+\int_{\bar{\rho}}^{\bar{\rho}J^{-1}}(\bar{\rho}
J^{-1}-z)P''(z)\mm{d}z.$$
Applying $\nabla_{\mathcal{A}}$-operator to the above identity, we find that
\begin{equation}  \label{press20160309MH}
\nabla_{\mathcal{A}}P(\sigma) = \nabla_{\mathcal{A}}P(\tilde{\rho})- \nabla (P'(\bar{\rho})\mm{div}(\bar{\rho} \eta )) -\mathcal{N}_3,
\end{equation}
where $\tilde{\rho}:=\bar{\rho}(\zeta_3)$, and
$$ \begin{aligned}
\mathcal{N}_3:=& (\nabla-\nabla_{\mathcal{A}}  )( P'(\bar{\rho})\bar{\rho}(J^{-1}-1)) +\nabla  ( P(\tilde{\rho})-\bar{P} -\bar{P}'\eta_3)\\
& + (\nabla_{\mathcal{A}}-\nabla)(P(\tilde{\rho})-\bar{P}) -\nabla (P'(\bar{\rho})\bar{\rho}(J^{-1}-1+\mm{div} \eta))
-\nabla_{\mathcal{A}}\int_{\bar{\rho}}^{\bar{\rho}J^{-1}}(\bar{\rho} J^{-1}-z)P''(z)\mm{d}z.
\end{aligned}$$

Finally, we represent the modified gravity term as follows.
\begin{equation}
\label{gravit20160309}
-\sigma \tilde{g}e_3:=  g \mm{div}(\bar{\rho}\eta)e_3- \tilde{\rho}
\tilde{g}e_3  +\mathcal{N}_4,   \end{equation}
where $\mathcal{N}_4:=(\bar{\rho}(1-J^{-1}) +\tilde{\rho}-\bar{\rho})(\tilde{g}-g)e_3+(\tilde{\rho} -\bar{\rho}'\eta_3 -\bar{\rho} J^{-1}-\bar{\rho}\mm{div}\eta)ge_3$.

Summing up the above calculations, we see that, if the initial data $\sigma_0$, $B_0$, $\zeta_0$ satisfy
\eqref{abjlj0i} and \eqref{zeta0inta}, then we can use the relations \eqref{leolasexipre}--\eqref{gravit20160309} and \eqref{relation201604091949} to
transform \eqref{0116} into the following evolution equations for $(\eta,u)$:
\begin{equation}   \label{peturbationequation}
  \left\{       \begin{array}{l}
\eta_t=u, \\[1mm]
 \bar{\rho} J^{-1}u_t -\mu_1 \Delta_{\ml{A}}u- \mu_2\nabla_{\mathcal{A}}\mm{div}_{\mathcal{A}}u
- \nabla (P'(\bar{\rho})\mm{div}(\bar{\rho} \eta ))
 =  g\mm{div}( \bar{\rho}\eta) e_3+\mathcal{L}_M+\mathcal{N},
\end{array}     \right.
\end{equation}
and $(\sigma,B)$ is given by \eqref{abjlj0}, where we have utilized the equilibrium state $s_e$ in Lagrangian coordinates
$ \nabla_{\mathcal{A}}(P(\tilde{\rho})+\lambda   |\tilde{M} |^2/2)=  -\tilde{\rho}\tilde{g} e_3$ for \eqref{peturbationequation},
and denoted $\mathcal{A}^T=(I+\nabla \eta)^{-1}$, $I=(\delta_{ij})_{3\times 3}$ and $\mathcal{N}:=\lambda (\mathcal{N}_1+ \mathcal{N}_2)+\mathcal{N}_3+\mathcal{N}_4$.
The associated initial and boundary conditions read as
\begin{equation}\label{defineditcon}
(\eta,u)|_{t=0}=(\eta_0,u_0), \qquad (\eta,u)|_{\partial\Omega}=0.
\end{equation}

In this paper, we call the initial-boundary problem \eqref{peturbationequation}--\eqref{defineditcon} the transformed {Parker problem}.
Compared with the original {Parker problem} \eqref{0103}--\eqref{0105}, the transformed {Parker problem} enjoys {a fine energy structure},
so that one can establish the stabilizing effect of a horizontal magnetic field with the help of non-slip boundary condition by the energy method.

\subsection{Nonlinear stability}
      Before stating our main result on the transformed {Parker problem}, we introduce some notations used throughout this paper. We denote
\begin{equation*}  \begin{aligned}& \mathbb{R}^+_0:=[0,\infty),\quad \int:=\int_\Omega,\quad
L^p:=L^p (\Omega):=W^{0,p}(\Omega)\;\mbox{ for }\; 1< p\leq \infty, \\
& {H}^1_0:=W^{1,2}_0(\Omega),\quad {H}^k:=W^{k,2}(\Omega ),\quad \|\cdot\|_k:=\|\cdot\|_{H^k(\Omega)}\;\mbox{ for }\;k\geq 0,\\
&\partial_\mm{v}^i\mbox{ deontes }\partial_{2}^{\alpha_1}\partial_{3}^{\alpha_2}\;\mbox{ for any }\alpha_1+\alpha_2=i,\quad
\|\cdot\|_{i,k}^2:=\sum_{\alpha_1+\alpha_2=i} \|\partial_{2}^{\alpha_1}\partial_{3}^{\alpha_2}\cdot\|_{k}^2, \\
&\|{f}\diamond g
\|_{i,k}^2:=\sum_{\alpha_1+\alpha_2=i} \sum_{\beta_1+\beta_2+\beta_3\leq k}\| f \partial_{3}^{\beta_3}\partial_{1}^{ \alpha_1+\beta_1}\partial_{2}^{ \alpha_2+\beta_2} g
\|_{0}^2,\quad \|{f}\diamond g\|^2_{\underline{i},k}:=\sum_{0\leq j\leq i}\|{f}\diamond g\|_{j,k}^2, \\
&\|\cdot\|^2_{\underline{i},k}:=\sum_{0\leq j\leq i}\|\cdot\|_{j,k}^2, \qquad
a\lesssim b\;\mbox{ means that }\; a\leq cb\;\mbox{ for some constant }c>0,\\
&  C_j^{j-l}\mbox{ denotes the  number of }(j-l)\mbox{-combinations from a given set }S\mbox{ of }j \mbox{ elements}. \end{aligned}\end{equation*}
The letter $c$ will denote a generic constant which may depend on the domain $\Omega$ and the physical parameters,
such as $\lambda $, $m$, $g$, $\mu$ and $\bar{\rho}$ in the original perturbation equations \eqref{0103}.
It should be noted that a product space $(X)^n$ of vector functions is still denoted by $X$,
for example, a vector function $ {u}\in (H^3)^3$ is denoted by ${u}\in H^3$ with norm
$\| {u}\|_{H^3}:=(\sum_{k=1}^3\|u_k\|_{H^3}^2)^{1/2}$. Finally, we define some functionals:
%%%%%%%%%%%%%%%%%%%%%%%%%%%%%%%%%%%%%
\begin{equation*}     \begin{aligned}
&\mathcal{E}^{L}(t):=\| \eta(t)\|_{4}^2+\|u(t)\|_{3}^2+\|u_t(t)\|_{1}^2,\quad
\mathcal{E}^{H}(t):=\|\eta(t)\|^2_{7}  +\sum_{k=0}^3 \|\partial_t^k u(t)\|_{6-2k}^2 , \\
&\mathcal{D}^{L}(t):=\| (\partial_1\eta,\mm{div}\eta)(t)\|^2_{3} + \sum_{k=0}^2\|\partial_t^ku(t)\|_{4-2k}^2 , \quad
\mathcal{D}^{H}(t):=   \| (\partial_1\eta,\mm{div}\eta)(t)\|^2_{6} +\sum_{k=0}^3\|\partial_t^k u(t)\|_{7-2k}^2 , \\
 &\mathcal{G}_1(t)=\sup_{0\leq \tau\leq t}\mathcal{E}^H(\tau)+\int_0^t
\mathcal{D}^H(\tau)\mm{d}\tau,\quad \mathcal{G}_2(t)=\sup_{0\leq \tau\leq t}(1+\tau)^{3}\mathcal{E}^L(\tau).
\end{aligned}          \end{equation*}

Now, our stability result of the transformed {Parker problem} reads as follows.
%%%%%%%%%%%%%%%%%%%%%%%%%%%%%%%%%%%%%%
\begin{thm}\label{thm3}
Let $\Omega$ be a vertical strip domain,  $(\bar{\rho},g)$ satisfy \eqref{rho20160326}, \eqref{modififorg} and \eqref{0102}, and
$m$ be given by \eqref{com:01091}. If $\Xi<0$,
then there is a sufficiently small $\delta >0$, such that for any
$(\eta_0,u_0)\in H^7\times H^6$ satisfying that
\begin{enumerate}
  \item[(1)] $ \| \eta_0 \|_7^2+\|u_0\|_{6}^2 \leq\delta $;
  \item[(2)] $\zeta_0:=y+\eta_0$ satisfies \eqref{zeta0inta};
  \item[(3)]  $(\eta_0, u_0)$ satisfies the compatibility conditions on boundary
   (i.e., $ \partial_t^j u(x,0)|_{\partial\Omega}=0$ for $j=0,1,2$),
\end{enumerate}
there exists a unique global solution $(\eta,u)\in C(\mathbb{R}^+_0,H^{7}\times H^6)$
to the transformed {Parker problem} \eqref{peturbationequation}--\eqref{defineditcon}. Moreover, $(\eta,u)$ enjoys
the following stability estimate:
\begin{equation}\label{1.19}
 \mathcal{G}(\infty):= \mathcal{G}_1(\infty) + \mathcal{G}_2(\infty) \lesssim  \| \eta_0 \|_7^2+\|u_0\|_{6}^2.
\end{equation}
Here the constant $\delta$ depends on $\Omega$ and the physical parameters in the original perturbation equations \eqref{0103}.
\end{thm}
\begin{rem}
  Since $\eta=0$ on boundary, then, if $\delta$ in Theorem 2.1 is sufficiently small, we can further have (referring to Lemma 4.2 in \cite{JFJSOMITN} for example)
 \begin{align}&
   \zeta:=\eta+y  : \overline{\Omega}\to \overline{\Omega} \mbox{ is a homeomorphism mapping},
\\&    \zeta  : \Omega \to \Omega \mbox{ are }C^5\mbox{-diffeomorphic mapping}. \label{2018031adsadfa21601xx}
\end{align}
Thus one can recover a stability result in  Eulerian coordinates  from Theorem \ref{thm3} by an inverse transformation of Lagrangian coordinates, referring to Theorem 1.2 in \cite{JFJSJMFMOSERT}.
\end{rem}

Next we briefly describe the basic idea in the proof of Theorem \ref{thm3}.
By the energy method, there exist two functionals $\tilde{\mathcal{E}}^{L}$ and $\mathcal{Q}$ of $(\eta,u)$
satisfying the lower-order energy inequality for the {Parker problem} (see Proposition \ref{pro:0301n})
\be   \label{latdervin}
  \frac{\mm{d}}{\mm{d}t}\tilde{\mathcal{E}}^{L} + \mathcal{D}^L\leq \mathcal{Q}\mathcal{D}^L, \ee
in which we can use the stability condition $\Xi<0$ to show that the functional $\tilde{\mathcal{E}}^{L}$
 is equivalent to $\mathcal{E}^L$. Unfortunately, we can
not close the energy estimates only based on \eqref{latdervin}, since $\mathcal{Q}$ can not be controlled by $\tilde{\mathcal{E}}^{L}$.
However, we observe that the structure of (\ref{latdervin}) is very similar to the one of the surface wave problem studied
in \cite{GYTIDAP,GYTIAED}, where Guo and Tice developed a two-tier energy method to overcome this difficulty.
In the spirit of the two-tier energy method, we shall look after a higher-order energy inequality to
match the lower-order energy inequality \eqref{latdervin}. Since $\tilde{\mathcal{E}}^L$ contains $\|\eta\|_4$,
 we can see that the higher-order energy includes at lest $\|\eta\|_7$. Thus, similarly to \eqref{latdervin},
 we are able to establish the higher-order energy inequality (see Proposition \ref{pro:0301}):
\be   \label{latdervins}
  \frac{\mm{d}}{\mm{d}t}\tilde{\mathcal{E}}^{H}+ {\mathcal{D}}^H\leq \sqrt{\mathcal{E}^{L}}\|\eta\|_7^2,  \ee
where the functional $\tilde{\mathcal{E}}^{H}$ is equivalent to $\mathcal{E}^H$ by the stability condition.
In the derivation of the \emph{a priori} estimates, we have $\mathcal{Q}\lesssim \mathcal{E}^H$, and thus
 \eqref{latdervin} implies (see Proposition \ref{pro:0301n})
 \be
\label{starlatdervin}
  \frac{\mm{d}}{\mm{d}t}\tilde{\mathcal{E}}^{L}+\mathcal{D}^L\leq 0.\ee
Consequently, with the help of the two-tier energy method, \eqref{latdervins} and \eqref{starlatdervin},
we can deduce the global-in-time stability estimate \eqref{1.19}.

We should remark that one of the main novelties in this article lies in that we further develop the stability condition to show that the energy
functionals $\tilde{\mathcal{E}}^{L}$  and $\tilde{\mathcal{E}}^{H}$ constructed in \eqref{latdervin} and \eqref{latdervins} are equivalent to $\mathcal{E}^L$
and $\mathcal{E}^H$, respectively. In fact, the energy of the {Parker problem} includes the term $-E(w)$ for $w=\partial_{\mm{v}}^j u$ and $\partial_t^j \eta$.
Therefore, we naturally expect that $-E(w)$ should be positive under the stability condition. Otherwise, the two-tier energy method will fail.
Exploiting the stability condition and the positivity of the term $\int P'(\bar{\rho})\bar{\rho}|\mathrm{div}w|^2\mm{d}y$, we can derive that $-E(w)$ is equivalent to
$\|(w,\partial_1w,\mm{div}_\mm{v}w_\mm{v},\mm{div}w)\|_0^2$ (see Lemma \ref{lem:0601}), which plays an important role in the proof of Theorem \ref{thm3}.

We mention that the basic idea of the two-tier energy method was also used to show the existence of the global stability solutions
to the MHD problem \cite{TZWYJGw}, and the stabilizing effect of the magnetic field in the magnetic Rayleigh-Taylor problem \cite{JFJSJMFMOSERT}.

Finally, in view of Theorem \ref{thm3}, we immediately obtain the existence of a unique global-in-time solution to the original {Parker problem},
which represents the strong stabilizing effect of the magnetic field through the non-slip boundary condition in the original {Parker problem}.
%%%%%%%%%%%%%%%%%%%%%%%%%%%%%%
\begin{thm}\label{thm3orig} Let $\Omega$ be a vertical strip domain, $(\bar{\rho},g)$ satisfy \eqref{rho20160326}, \eqref{modififorg} and \eqref{0102}, and
$m$ be given by \eqref{com:01091}. If $\Xi<0$,
then there is a sufficiently small $\delta >0$, such that for any $(\varrho_0, v_0,N_0)\in H^{6}$ satisfying that
\begin{enumerate}
    \item[(1)]  there exists an invertible mapping $\zeta_0:=\zeta_0(x):\Omega\to \Omega$ such that
\eqref{zeta0inta};
 \item[(2)]  $ \| \zeta_0 -x\|_7^2+\|v_0\|_{6}^2 \leq\delta$;
\item[(3)]$\rho_0(\zeta_0)=\bar{\rho}(\det(\nabla \zeta_0))^{-1}$, $(N_0+\bar{M})(\zeta_0)=
 m (\det(\nabla \zeta_0))^{-1}\partial_1\zeta_0 $ and   $\div_{\ml{A}_0}(
 (N_0+\bar{M})(\zeta_0))=0$, where
$\ml{A}^T_0=(\nabla \zeta_0)^{-1}$;
  \item[(4)]the initial data $(\varrho_0, v_0, N_0)$ satisfies necessary compatibility conditions
  (i.e., $ \partial_t^jv(x,0)|_{\partial\Omega}=0$ for $j=0,1,2$),
\end{enumerate}
there exists a unique global solution $(\varrho,v,N)\in C(\mathbb{R}^+,H^6)$ to the original {Parker problem} \eqref{0103}--\eqref{0105}.
Moreover, $(\varrho, v,N)$ enjoys the following stability estimate:
\begin{equation}\label{1.19st}
\begin{aligned}
 &\sup_{0\leq t<\infty}\left(\|(\varrho,N)\|_6^2+\sum_{k=0}^3 \|\partial_t^k {v}(t)\|_{6-2k}^2\right)\\
&\qquad +\sup_{0\leq t< \infty}(1+t)^{3}(\|(\varrho,v,N)\|_{3}^2+\|v_t\|_{1}^2)\leq   \|\zeta_0-x \|_7^2+\|v_0\|_{6}^2.
\end{aligned}
\end{equation}
Here the constant $\delta$ depends on $\Omega$ and the physical parameters in the original perturbation equations  \eqref{0103}.
\end{thm}

\subsection{Nonlinear instability}
      As aforementioned, $\Xi$ is the discriminant for the instability/stability of the linearized {Parker problem}.
 Next we give an existence result of a local unstable solution to the original {Parker problem} for $\Xi>0$, which, together with Theorem \ref{thm3orig},
 implies that $\Xi$ is also the discriminant for the instability/stability of the original {Parker problem}.
 %%%%%%%%%%%%%%%%%%%%%%%%%%%%%%%%%%%%%%%%%%%%
\begin{thm}\label{thm:0101jfjsww101315}
Let $\Omega$ be a vertical strip domain,  $(\bar{\rho},g)$ satisfy \eqref{rho20160326}, \eqref{modififorg} and \eqref{0102}, and
$m$ be given by \eqref{com:01091}. If $\Xi> 0$, then the {Parker problem} \eqref{0103}--\eqref{0105}
 is unstable, that is, there are positive constants $\varepsilon$ and $\iota$, and a quaternion $(\tilde{\varrho}_0,\tilde{ v}_0,\tilde{ {N}}_0, { v}_{\mm{r}})\in H^3$,
such that for any $\delta\in (0,\iota)$ and the initial data
$(\varrho_0, v_0, {N}_0):=\delta(\tilde{\varrho}_0,\tilde{v}_0, \tilde{ {N}}_0)+(0,\delta^2 v_{\mm{r}},0)$,
there is a unique classical solution $(\varrho ,v,N)$ of the {Parker problem} \eqref{0103}--\eqref{0105} on $[0,T^{\max})$, but
\begin{equation} \label{201602082017MH}
\|(\varrho ,u,N)(T )\|_{3}\geq {\varepsilon}\;\;\mbox{ for some time }\; T\in (0,T^{\max}),  \end{equation}
where $({\varrho}_0,{ v}_0,{ {N}}_0)$ satisfies the compatibility conditions $\partial_t^i v(x,0)|_{\partial\Omega}=0$ ($i=0,1$) and $\mathrm{div} N_0=0$,
$T^{\max}$ denotes the maximal time of existence of the solution $(\varrho, v, {N})$, and $\varepsilon$, $\iota$,
and $(\tilde{\varrho}_0,\tilde{ v}_0,\tilde{ {N}}_0,v_{\mm{r}})$
depend on $\Omega$ and the physical parameters in the original perturbation equations \eqref{0103}.
\end{thm}　
{Theorem \ref{thm:0101jfjsww101315} reveals that the Parker instability will occur for $\Xi>0$. Roughly speaking, Theorem \ref{thm:0101jfjsww101315}
  is proved based on a new version of the bootstrap method.} The method of bootstrap instability probably started
from Guo and Strauss' works \cite{GS1995,GS1995neG}. Then, various versions of the bootstrap method have been developed in the study of dynamical instability
of various physical models, please refer to \cite{GYHCSDDUI,HHJGY,JFJSO2014,FSNPVVNC,FSSWVMNA,FrrVishikM} for more details.
 Unfortunately, the known versions of the bootstrap method can not be directly applied to the {Parker problem} here due
 to the presence of the magnetic field and non-slip boundary condition. In this paper we develop a new version of the bootstrap method
 to overcome such difficulties and establish the Parker instability. The basic idea of our new bootstrap method can be found
 after Lemma \ref{lem:201602091151MH}, and the key technique
mainly lies in the derivation of the error estimate \eqref{0508} in Lemma \ref{lem:0401}, which is obtained by exploiting
the structure of the original perturbation equations. Moreover, we develop an approach based on the elliptic theory to construct
the compatible initial data of the nonlinear unstable solution by using the initial data of the linear unstable solution,
see Lemma \ref{lem:modfied}. We mention that Theorem \ref{thm:0101jfjsww101315} also holds for a bounded $C^4$-domain
by applying our new bootstrap method, and moreover, the smoothness requirements \eqref{modififorg}--\eqref{rho20160326mB6} can be relaxed.

\subsection{Horizontally periodic domains}\label{201605171439}
     Now we further consider the case of a horizontally periodic domain, i.e.,
\begin{equation}    \label{201605161646}
\Omega:=\{ x=(x_h,x_3)\in \mathbb{R}^3~|~x_h\in \mathcal{T},\ -l<x_3<l\}\;\; \mbox{ with }\; l>0,
\end{equation}
  where $x_h:=(x_1,x_2)$, $\mathcal{T} :=(2\pi L_1\mathbb{T})\times(2\pi L_2\mathbb{T}) $, $\mathbb{T}=\mathbb{R}/\mathbb{Z}$, and $2\pi L_1$, $2\pi L_2>0$
  are the periodicity lengths. Then we can verify that $\Xi>0$ {under Schwarzschild's condition with sufficiently large $L_1$ or Tserkovnikov's condition.}
The detailed verification will be presented in Section \ref{201605171506}. Hence, we can follow the proof of Theorem \ref{thm:0101jfjsww101315} to establish
the nonlinear Parker instability under Schwarzschild's or Tserkovnikov's condition.
%%%%%%%%%%%%%%%%%%%%%%%%%%%%%%%%%%%%%%%%%%%%%%%%%

Before stating our nonlinear Parker instability results in a horizontally periodic domain, we introduce some notations:
$$\begin{aligned}
&\chi(\psi):= \int_{-l}^l
\left({\lambda m^2}\left(|\psi'|^2+\frac{\psi^2 }{L^2_2}\right)-\left( \frac{g^2\bar{\rho}^2}{\gamma  \bar{P}}+g\bar{\rho}'\right)
\psi^2\right)\mm{d}x_3,\\
 &\xi_{3D}:=\sup_{\psi \in H^1_0(-l,l) }\sqrt{
\frac{\sqrt{ \chi^2(\psi)+\frac{4}{L_2^2}\int_{-l}^l  \lambda m^2\psi^2\mm{d}x_3\int_{-l}^l
\left(\frac{g^2\bar{\rho}^2}{\gamma  \bar{P}}+g\bar{\rho}'
\right)\psi^2 \mm{d}x_3}-
 \chi(\psi)}{2 \int_{-l}^l  \lambda m^2\psi ^2\mm{d}x_3}},\\
&\kappa(l):=\sup_{\psi\in H^1_0(-l,l)}\sqrt{\frac{\int_{-l}^l
 \left(\frac{g^2\bar{\rho}^2}{\gamma \bar{P}}+g\bar{\rho}'\right)\psi^2\mm{d}x_2}
{\int_{-l}^l {\lambda  m^2}|\psi'|^2\mm{d}x_2}},\\
& \xi_{2D}:=\sup_{\psi\in H^1_0({-l},l)}\sqrt{\frac{\int_{-l}^l
  \left(\left(\frac{g^2\bar{\rho}^2}{\gamma \bar{P}}+g\bar{\rho}'\right)\psi^2-{\lambda  m^2}|\psi'|^2\right)\mm{d}x_2}{\int_{-l}^l \lambda m^2 \psi^2\mm{d}x_2}},
\end{aligned} $$
where we have omitted to mention that $\psi\in H_0^1(-l,l)$ should make the denominators and square roots sense in the definitions of $\xi_{3D}$ and $\xi_{2D}$.
Next we state the main results.
%%%%%%%%%%%%%%%%%%%%%%%%%%%%%%%%%%%%%%%%%%%%%%%%%%%
\begin{thm}\label{thm:201605151841}
Let $l\in (0,\infty)$, $\Omega$ be defined by
\eqref{201605161646}, $g$ be a constant, $\bar{\rho}:=\bar{\rho}(x_3)\in C^4(\bar{\Omega})$ satisfy \eqref{0102}, and
$m$ be given by \eqref{com:01091}. If  one of the following two assumptions holds
\begin{enumerate}
  \item[(1)]$\bar{\rho}$ satisfies Schwarzschild's  condition \eqref{Newcombinstablity}, and $L_1>\xi_{3D}^{-1}$;
  \item[(2)] $\bar{\rho}$ satisfies Tserkovnikov's condition \eqref{Tserkovnikovs};
\end{enumerate}
then, the {Parker problem} \eqref{0103}--\eqref{0105} is unstable as in Theorem \ref{thm:0101jfjsww101315}.\end{thm}
\begin{rem}
By virtue of Schwarzschild's condition, there exists an open interval $I\subset (-l,l)$, such that
$$-\bar{\rho}'<\frac{g\bar{\rho}^2}{\gamma \bar{P}}\quad \mbox{for any  }x_3\in I.$$
Then we choose a function $\psi_0\in H^1_0(-l,l)$ such that
$$\psi_0(x_3)>0\quad\mbox{in }\; I\mbox{ and }\psi_0(x_3)=0\mbox{ in }({-l},l) \setminus I,$$
whence,
$$\int_{-l}^l \left(\frac{g^2\bar{\rho}^2}{\gamma  \bar{P}}+g\bar{\rho}' \right)\psi^2_0 \mm{d}x_3>0,$$
which implies that $\xi_{3D}$ must be a positive constant. Moreover, for given $\bar{\rho}$, we have
$$\xi_{3D}\to 0\quad\mbox{ as }|m|\to \infty.$$
\end{rem}
\begin{rem}\label{201608031854}
It should be noted that Schwarzschild's condition is equivalent to the magnetic buoyancy condition. Hence, Theorem \ref{thm:201605151841} tells us that if there is a point,
at which the magnetic buoyancy points oppositely to the direction of the gravity field in the  equilibrium state $s_e$,
then the {Parker problem} is unstable. If the direction of the magnetic buoyancy is in line everywhere with the direction of the gravity field
in the  equilibrium state, i.e., Schwarzschild's   condition fails, then we have $\Xi<0$  for a vertical strip domain with non-slip boundary condition, which immediately implies that
the {Parker problem} is stable by virtue of Theorem \ref{thm3orig}. We can observe that the proof of Theorem \ref{thm3orig} strongly depends on the non-slip condition
in the horizontal direction. A question arises whether this stability conclusion in Theorem \ref{thm3orig} can be generalized
to a horizontally periodic domain. We shall further investigate this question in a separate article.
\end{rem}
%%%%%%%%%%%%%%%%%%%%%%%%%%%%%%%%%
\begin{rem} If $l=\infty$ and $g\in C_0^3(\mathbb{R})$, we can establish a Gronwall-type inequality as
(3.39) in \cite[Proposition 3.2]{JFJSWWWN} for the classical solution $(\varrho,u,N)\in C^0([0,T),H^3)$ of the {Parker problem}.
Thus we can deduce a more precise result on the nonlinear Parker instability for both 2D and 3D cases by a standard bootstrap method,
in which the third component of the velocity is unstable in $L^2$ as  (1.10) in \cite[Theorem 1.1]{JFJSWWWN}.
\end{rem}

Similarly, we can also establish the following nonlinear Parker instability result in the 2D case.
%%%%%%%%%%%%%%%%%%%%%%%%%%%%%%%%%%%%%%%%%%%%%
\begin{thm}\label{201605161704n} Let $l\in (0,\infty)$, $\Omega:=2\pi L_1\mathbb{T}\times ({-l},l )\subset \mathbb{R}^2$, $g$ be a constant, $\bar{\rho}:=\bar{\rho}(x_2)\in C^4(\Omega)$ satisfy \eqref{0102}, and $m$ be given by \eqref{com:01091}. If $\kappa>1$ and $L_1>\xi_{2D}^{-1}$, then the {Parker problem}
\eqref{c0104}--\eqref{lincom} in the 2D case is unstable as in Theorem \ref{thm:201605151841}.
\end{thm}
\begin{rem}
It is interesting to notice that the definition of $\kappa$ is very similar to the critical number $M_\mathrm{c}(l)$ of the
3D incompressible magnetic RT problem around a vertical equilibrium magnetic field, where
  \begin{equation*}\begin{aligned}
M_\mathrm{c}(l):=\sqrt{\sup_{\psi\in H^1_0({-l},l )}
\frac{g\int_{-l}^l   \bar{\rho}' |\psi |^2\mm{d}s}
{\lambda\int_{-l}^l |\psi' |^2\mm{d}s}}>0\mbox{ for }\bar{\rho}'(s)>0\mbox{ for some }s\in (-l,l) \end{aligned}\end{equation*}
and $M_\mathrm{c}(l)<\infty$ for $l<\infty$, please refer to \cite[Theorem 1.1]{JFJSWWWN}.
Moreover, $M_\mathrm{c}(l)$ is also the critical number $M_\mathrm{c}(l)$ of the
2D incompressible magnetic RT problem around a horizontal equilibrium magnetic field, please refer to \cite{WYC}.
\end{rem}
\begin{rem}
Under Schwarzschild's  condition, we can have that
 $$0<\kappa(l)<\infty \mbox{ for }l<\infty.  $$
Moreover, $\kappa\to 0$ as $m\to \infty$ for a given $\bar{\rho}$. On the other hand, if $\kappa\in (0,1]$, we can use the Fourier analysis
method to obtain $E(w)\leq 0$ for any $w\in H_0^1$, which implies $\Xi\leq 0$. Thus we see that
a sufficiently large horizontal magnetic field has the stabilizing effect in the 2D case.
This conclusion agrees with the result for the 2D incompressible magnetic RT problem around a horizontal equilibrium magnetic field.
In addition, if $l=\infty$ and $0\leq g\in C_0^0(\mathbb{R})$, under Schwarzschild's  condition, we can derive that
$$ \kappa(\infty)=\left\{
      \begin{array}{ll}
        \mbox{ infinity}, & \hbox{ if }\int_{\mathbb{R}}\left(\frac{g^2\bar{\rho}^2}{\gamma \bar{P}}+g\bar{\rho}'
\right)\mm{d}x_2>0; \\[3mm]
        \mbox{ a positive real number}, & \hbox{ if }\int_{\mathbb{R}}\left(\frac{g^2\bar{\rho}^2}{\gamma \bar{P}}+g\bar{\rho}'
\right)\mm{d}x_2<0,
      \end{array}
    \right.  $$　
please refer to \cite[Proposition 2.1]{JFJSWWWN} for the proof.
\end{rem}

The rest sections are  mainly devoted to the proof of Theorems \ref{thm3}--\ref{201605161704n}. In
Section \ref{lowerorder} we first derive the lower-order energy inequality  \eqref{starlatdervin} of the transformed {Parker problem}. Then
 in Section \ref{sec:03}, we derive the higher-order energy inequality \eqref{latdervins}. Finally, based on the previous two energy
 inequalities, we show Theorem \ref{thm3} by applying the two-tier energy method,
  and further deduce Theorem \ref{thm3orig} from Theorem \ref{thm3} in Section \ref{sec:02}.
In Section \ref{sec:06},  we develop a new version of the bootstrap method to prove Theorem \ref{thm:0101jfjsww101315}.
Finally, we verify that $\Xi>0$ under the assumptions of Theorem \ref{thm:201605151841}--\ref{201605161704n}
in Section \ref{201605171506}.

\section{Lower-order energy inequality}\label{lowerorder}
In this section we derive the lower-order energy inequality for the transformed {Parker problem}. To this end, let $(\eta,u)$ be a solution
of the transformed {Parker problem}, such that
\begin{equation}\label{aprpioses}
\sqrt{\sup_{0\leq \tau\leq T}\mathcal{E}^H(\tau)}\leq \delta\in (0,1)\;\;\mbox{ for some  }T>0,
\end{equation}
where $\delta$ is sufficiently small. It should be noted that the smallness depends on the domain $\Omega$ and the physical parameters in the perturbation equations \eqref{0103}.
%% and  will be repeatedly used in what follows.
Moreover, we assume that the solution $(\eta,{u})$ possesses proper regularity, so that the procedure of formal calculations makes sense.
In the calculations that follow, we shall repeatedly use Cauchy-Schwarz's inequality, H\"older's inequality,
and the embedding inequalities (see \cite[4.12 Theorem]{ARAJJFF})
\begin{align}
&\label{esmmdforinftfdsdy}  \|f\|_{L^p}\lesssim \| f\|_{1}\quad\mbox{for }2\leq p\leq 6,\\
&\label{esmmdforinfty}   \|f\|_{L^\infty}\lesssim\| f\|_{2}.\end{align}
and the interpolation inequality in $H^j$ (see \cite[5.2 Theorem]{ARAJJFF})
\begin{equation}
\label{inteposlai}
\|f\|_j\lesssim \|f\|_0^{1-\frac{j}{i}}\|f\|_i^{\frac{j}{i}}\leq C_\epsilon\|f\|_{0} +\epsilon\|f\|_{{j+1}}
\end{equation}
for any $0\leq j< i$ and any constant $\epsilon>0$, where the constant $C_\epsilon$ depends on $\Omega$ and $\epsilon$.
In addition, we shall also repeatedly use the following two estimates:
\begin{equation}
\label{fgestims}
 \|fg\|_j\lesssim \left\{     \begin{array}{ll}
                      \|f\|_1\|g\|_{1} & \hbox{ for }j=0;  \\   \|f\|_j\|g\|_{2} & \hbox{ for }0\leq j\leq 2;  \\
                    \|f\|_2\|g\|_j+\|f\|_j\|g\|_2& \hbox{ for }3\leq j\leq 5;\\
 \|f\|_2\|g\|_j+\|f\|_{5}\|g\|_{j-3}+\|f\|_{j}\|g\|_{2}& \hbox{ for }6\leq j\leq 7
                    \end{array}         \right.
\end{equation}
and
\begin{equation}
\label{insoslai}\|f\|_0\lesssim \|\partial_1 f\|_0\mbox{ for }f\in H^1_0,\end{equation}
where \eqref{fgestims} can be easily verified by H\"older's inequality and the embedding inequalities \eqref{esmmdforinftfdsdy}--\eqref{esmmdforinfty}.

Before deriving the lower-order energy inequality defined on $(0,T]$, we first give some preliminary estimates, temporal derivative estimates,
$y_\mm{v}$-derivative estimates (i.e., the estimates of partial derivatives with respect to $y_2$ and $y_3$) and $y$-derivative estimates
 (i.e., the estimates of partial derivatives with respect to $y_1$, $y_2$ and $y_3$) in sequence.

\subsection{Preliminary estimates}
In this subsection we introduce some preliminary estimates on $J$, $J^{-1}$, $\mathcal{A}$ and $\mathcal{N}$, which will be repeatedly
used in estimating $(\eta, u)$ later.
\begin{lem}\label{estimatelem0101} The following estimates hold.
\begin{align}
\label{Jdetemrinat}&1\lesssim \| J\|_{L^\infty} \lesssim1,
\\& \label{Jdetemrinatnesw}  \|J-1\|_{l}\lesssim  \| \eta\|_{l+1},\\& \label{Jtestismts}\|\partial_t^i J\|_{j}\lesssim \sum_{  k=0}^{i-1}\|\partial_t^{k} u\|_{j+1} ,\\
\label{Jdetemrinat22}&1\lesssim \| J^{-1}\|_{L^\infty} \lesssim 1,\\
 &\label{Jdetemrinat2nes}\| J^{-1}-1\|_{l}\lesssim \| \eta\|_{l+1},\\
 &\label{Jtestismtsinver}\|\partial_t^i J^{-1} \|_j\lesssim \sum_{  k=0}^{i-1}\|\partial_t^k u\|_{j+1},
  %& \|J\|_j\lesssim\|\nabla \eta\|,
\end{align}where $0\leq l\leq 6$, $1\leq  i\leq 3$ and $0\leq  j\leq  7-2i$.
 \end{lem}
 \begin{pf} Recalling the definition of $J=\det(\nabla \eta+I)$ and using the expansion theorem of determinants, we find that
\begin{equation}
\label{Jdetmentionsa}
J=1+\mm{div}\eta+P_2(\nabla \eta)+P_3(\nabla \eta),
\end{equation}
 where $P_i(\nabla \eta)$ ($i=2,3$) denotes the homogeneous polynomial of degree $i$ with respect to $\partial_j\eta_k$ for $1\leq j$, $k\leq 3$.
 Using  \eqref{fgestims}, \eqref{esmmdforinfty}, and the smallness condition \eqref{aprpioses}, we immediately
 get \eqref{Jdetemrinatnesw} and \eqref{Jdetemrinat}. Similarly, we easily obtain \eqref{Jtestismts} from \eqref{Jdetmentionsa} and
 \eqref{peturbationequation}$_1$.

  By \eqref{Jdetmentionsa} and \eqref{Jdetemrinat}, we have \eqref{Jdetemrinat22} and
%%  \begin{equation*}   \label{JinverfsofJdsffd}
$J^{-1}=(1+\mm{div}\eta+P_2(\nabla \eta)+P_3(\nabla \eta))^{-1}$,
% \end{equation*}
which implies
$$
J^{-1}-1=-(J^{-1}-1)\left(\mm{div}\eta+P_2(\nabla \eta)+P_3(\nabla \eta)\right)-\left(\mm{div}\eta+P_2(\nabla \eta)+P_3(\nabla \eta)\right).
 $$
Thus, we obtain \eqref{Jdetemrinat2nes} by using \eqref{fgestims} and the smallness condition.

Finally, noting that
$$ \begin{aligned}   &J^{-1}_t=-J^{-2}J_t,\quad J^{-1}_{tt}=2J^{-3}J_t^2-J^{-2}J_{tt},\quad
J^{-1}_{ttt}=-6J^{-4}J_t^3+6J^{-3}J_{t}J_{tt}-J^{-2}\partial_t^3J,\\
&\|J^{-1}\psi\|_l\lesssim \|(J^{-1}-1)\psi\|_l+\|\psi\|_l\lesssim \|\psi\|_l\quad\mbox{ for }\; 0\leq l\leq 6,
\end{aligned}$$
we make use of \eqref{Jtestismts} and \eqref{fgestims} to get \eqref{Jtestismtsinver} immediately.
 This completes the proof.   \hfill$\Box$
\end{pf}
%%%%%%%%%%%%%%%%%%%%%%%%%%
\begin{lem}\label{estimatelem0102} It holds that
\begin{align} & \label{aimdse}
\|\mathcal{A}\|_{L^\infty} \lesssim 1 ,  \\
&\label{prtislsafdsfs}\| \partial_{t}^i\ml{A}\|_j \lesssim \sum_{k=0}^{i-1} \|\partial_t^k  u\|_{j+1},\\
&\label{prtislsafdsfsfds}\|\tilde{\mathcal{A}}\|_{l} \lesssim \|  \eta\|_{l+1},
\end{align}
where $\tilde{\mathcal{A}}:=(\mathcal{A}-I)$, $0\leq l\leq 6$, $1\leq i\leq 3$ and $0\leq j\leq 7-2i$.
\end{lem}
\begin{pf} Recalling the definition of $\mathcal{A}$, we see that $\mathcal{A}=(A^*_{kl})_{3\times3}J^{-1}$,
where $A^{*}_{kl}$ is the algebraic complement minor of $(k,l)$-th entry in the matrix $\nabla\eta+I$ and
the polynomial of degree $1$ or $2$ with respect to $\partial_m\eta_n$ for $1\leq m$, $n\leq 3$.
Thus, employing \eqref{Jdetemrinat22}, \eqref{esmmdforinfty} and the smallness condition, we obtain \eqref{aimdse}.

Using \eqref{Jdetemrinat2nes}, \eqref{fgestims}, \eqref{peturbationequation}$_1$ and the smallness condition, we have
$$  \begin{aligned}
&   \|\partial_t^{i} A^*_{kl}\|_j \lesssim  \sum_{m=0}^{i-1} \|\partial_t^m  u\|_{j+1} ,\quad
\| A^*_{kl}\partial_t^iJ^{-1}\|_{j}\lesssim (1+ \|\eta\|_{6})\|\partial_t^i J^{-1}\|_j\lesssim  \|\partial_t^i J^{-1}\|_j,\\
&\|\partial_t^i A^*_{kl}J^{-1}\|_{j}\lesssim  \|\partial_t^i A^*_{kl} (J^{-1}-1)\|_j
+\|\partial_t^i A^*_{kl}\|_j\lesssim \|\partial_t^i A^*_{kl}\|_j,\quad 1\leq i\leq 3,\; 0\leq  j\leq  7-2i.
\end{aligned}$$
Thus we can further make use of the above three estimates, \eqref{Jtestismtsinver} and \eqref{fgestims} to deduce that
$$\|\partial_t^i\mathcal{A}\|_{j}=\left\|\sum_{m=0}^i C_i^m(\partial_t^{i-m} A^*_{kl})_{3\times3}\partial_t^{m} J^{-1}\right\|_{j}
\lesssim \sum_{k=0}^{i-1} \|\partial_t^k  u\|_{j+1},$$
which yields \eqref{prtislsafdsfs}.

  We proceed to evaluating $\tilde{\mathcal{A}}$. Since $\delta$ is assumed to be so small that the following power series holds.
$$\mathcal{A}^T=I-\nabla \eta+(\nabla\eta)^2\sum_{i=0}^\infty (-\nabla\eta)^i=I-\nabla \eta+(\nabla\eta)^2\mathcal{A}^T,$$
whence,
\begin{equation}
\label{201603301551}
\tilde{\mathcal{A}}^T =(\nabla\eta)^2\mathcal{A}^T-\nabla \eta=(\nabla\eta)^2\tilde{\mathcal{A}}^T+(\nabla\eta)^2-\nabla \eta.
\end{equation}
Hence, we can use \eqref{fgestims} to deduce \eqref{prtislsafdsfsfds} from \eqref{201603301551} for sufficiently small $\delta$.
The proof is complete. \hfill$\Box$
\end{pf}

\begin{lem} It holds that
 \begin{align}&
\label{N21556201630266} \|\mathcal{N}\|_j\lesssim \| \eta\|_{j-1}\|\eta\|_7\quad \mbox{ for }j=3\mbox{ and }5, \\
 &\label{N215562016302661}\|\mathcal{N}_t\|_2\lesssim \|\eta\|_4\|u\|_4,\\
% &\label{N215562016302662}\|\partial_t^i\mathcal{N}\|_0\lesssim %(\|\eta\|_4+\|u\|_3)\|u_t\|_2\mbox{ for },\\
 &\label{N215562016302661Ntt}\|\mathcal{N}_{tt}\|_0\lesssim \|\eta\|_4(\|u\|_2+\|u_t\|_2)+\|u\|_3^2,
 \\&\label{Nu1743MHtdfsNsg} \|\mathcal{N}_t\|_{3}+ \| \mathcal{N}_{tt}\|_{1}+
\|\partial_t^3 \mathcal{N}\|_{0}  \lesssim \sqrt{\mathcal{E}^L\mathcal{D}^H}.
\end{align}    \end{lem}
 \begin{pf}   Since $\mathcal{N}$ is a linear combination of $\mathcal{N}_i$ for $i=1$, $\ldots$, $4$, it suffices to verify the above estimates
 with $\mathcal{N}_i$ in places of $\mathcal{N}$. By a straightforward calculation, we find that
  $$\begin{aligned}
&\nabla_{\mathcal{A}}-\nabla=\nabla_{\tilde{\mathcal{A}}}, \\
&B-\bar{M}=m((J^{-1}-1 )\partial_1 \eta+\partial_1 \eta+(J^{-1}-1)e_1),  \\
&|\tilde{M}|^2-\bar{M}^2 = 2\int_{0}^{\eta_3} (mm')(s+y_3)\mm{d}s, \\
&  |\tilde{M}|^2-\bar{M}^2-2mm'\eta_3 = 2  \int_{0}^{\eta_3} \int_{0}^s(mm')'(\tau+y_3)\mm{d}\tau\mm{d}s,\\
&\label{JinverfsofJ1307}
J^{-1}-1+\mm{div}\eta=J^{-1}\left((J-1)\mm{div}\eta-P_2(\nabla \eta)-P_3(\nabla \eta)\right),\\
 &P(\tilde{\rho})-\bar{P} =  \int_{0}^{\eta_3} \frac{\mm{d}}{\mm{d}y_3}P(\bar{\rho}(s+y_3))\mm{d}s,\\
 &P(\tilde{\rho})-\bar{P} -\bar{P}'\eta_3=
 \int_{0}^{\eta_3}\int_{0}^s \frac{\mm{d}^2}{\mm{d}y_3^2}P(\bar{\rho}(\tau+y_3))\mm{d}\tau\mm{d}s,\\
& \int_{\bar{\rho}}^{\bar{\rho}J^{-1}}(\bar{\rho}
J^{-1}-z)P''(z)\mm{d}z=\int_{0}^{\bar{\rho}(J^{-1}-1)}(\bar{\rho}
(J^{-1}-1)-z)P''(z+\bar{\rho})\mm{d}z,\\
& \tilde{\rho}-   \bar{\rho}'\eta_3-\bar{\rho} J^{-1}-\bar{\rho}\mm{div}\eta\\
 & \quad = \int_0^{\eta_3}\int_{0}^{s} \bar{\rho}''(\tau+y_3)\mm{d}\tau\mm{d}s
 +\bar{\rho}(J^{-1}-1)\mm{div}\eta+\bar{\rho}J^{-1}\left(P_2(\nabla \eta)+P_3(\nabla \eta)\right).
 \end{aligned}$$
Putting the above relations into the expressions of $\mathcal{N}_1$--$\mathcal{N}_4$, we use Lemmas \ref{estimatelem0101}--\ref{estimatelem0102}, \eqref{inteposlai}--\eqref{fgestims} and \eqref{modififorg}--\eqref{rho20160326mB6}
to see that \eqref{N21556201630266}--\eqref{Nu1743MHtdfsNsg} hold with $\mathcal{N}_i$
in place of $\mathcal{N}$ for $i=1$, $\ldots$, $4$.  This completes the proof. \hfill $\Box$\end{pf}

Finally, we derive an important estimate from the stability condition.
\begin{lem}\label{lem:0601} If $\Xi<0$, then
%, there exists a constant $c$, depending  on the domain $\Omega$ and other physical %parameters in the perturbation equations, such that
\begin{equation*}
\| (w,\partial_1w,\mm{div}_\mm{v}w_\mm{v},\mm{div}w)\|_0^2\lesssim -  E (w), \end{equation*}
where $E(w)$ is defined by \eqref{20160306cond}.
\end{lem}
\begin{pf}
Since $\Xi<0$, one has $-\Xi\lambda\| (\partial_{1}w_\mm{v},\mm{div}_\mm{v}w_\mm{v})\|^2_0\leq -E(w)$, which, together with \eqref{insoslai}, yields
 \begin{equation}   \label{lem:0601fd}
 \| (w_\mm{v} ,\partial_{1}w_\mm{v}, \mm{div}_\mm{v}w_\mm{v})\|^2_0\lesssim -  {E (w)}.
\end{equation}
 On the other hand, by virtue of \eqref{0102},
 $$ \begin{aligned}
\|\mathrm{div}w\|^2_0 \lesssim  & \int \left( P'(\bar{\rho})\bar{\rho}|\mathrm{div}w|^2
 + \lambda m^2 (|\partial_{1}w_\mm{v}|^2 +|\mm{div}_\mm{v}w_\mm{v}|^2  ) \right) \dx\\
=&\int g\bar{\rho}'w_3^2\mm{d} x +\int 2 g\bar{\rho}\mm{div}ww_3\dx - E(w),
\end{aligned}$$
which, combined with \eqref{modififorg}, \eqref{rho20160326} and the relation
\begin{equation}   \label{divfefdMH}
\mm{div}w=\partial_1w_1+\mm{div}_{\mm{v}}w_{\mm{v}}, \end{equation}
results in
 \begin{equation}
 \label{lem:0601fdds}c_1\|(\partial_1 w_1,\mathrm{div}w)\|^2_0 \leq c_2\|(w_3,\mm{div}_{\mm{v}}w_{\mm{v}})\|_0^2 - E(w),
 \end{equation}
where the positive constants $c_{1}$ and $c_2$ depend on $\Omega$ and the physical parameters. Therefore,
the desired conclusion follows from \eqref{lem:0601fdds}, \eqref{lem:0601fd} and \eqref{insoslai} immediately. \hfill $\Box$
\end{pf}
\subsection{Temporal derivative estimates}
In this subsection we establish  the estimates of temporal derivatives. To this end, we apply $\partial_t^j$
to \eqref{peturbationequation} to get
 \begin{equation}\label{01s06p}
\begin{cases}  \partial_t^{j+1}\eta=\partial_t^{j}u, \\[1mm]
 \bar{\rho} J^{-1}\partial_t^{j+1}u -\mu_1 \Delta_{\ml{A}}\partial_t^{j}u-\mu_2
 \nabla_{\mathcal{A}}\mm{div}_{\mathcal{A}}\partial_t^{j}u
-\nabla (P'(\bar{\rho}) \mm{div}(\bar{\rho}\partial_t^{j}\eta))  \\
 \quad = g\mm{div} (\bar{\rho}\partial_t^{j }\eta) e_3 +  \partial_t^{j }\mathcal{L}_M+\partial_t^{j }\mathcal{N} +  N^{t,j}_u  ,  \\[1mm]
 \end{cases}
\end{equation}
where
$$\begin{aligned}
&N^{t,j}_u:=\sum_{0\leq k<j,\ 0\leq l\leq j} ( \mu_1 C_{j}^{k+l}C_{k+l}^{k}\partial_t^{j-k-l }\mathcal{A}_{il}\partial_{l}(\partial_t^{l}\mathcal{A}_{ik}\partial_t^{k}
\partial_ku)\\
&\qquad\qquad +\mu_2  C_{j}^{k+l}C_{k+l}^{k}(\partial_t^{j-k-l }\mathcal{A}_{il}\partial_{l}(\partial_t^{l}\mathcal{A}_{sk}\partial_t^{k}
\partial_ku_s))_{3\times1})-\sum_{ k=0}^{j-1} \bar{\rho} \partial_t^{j-k}J^{-1}\partial_t^{k+1}u.
\end{aligned}$$
 Moreover, making use of \eqref{prtislsafdsfsfds}, \eqref{prtislsafdsfs}, \eqref{Jtestismtsinver}, \eqref{fgestims}, and interpolation inequality,
 we easily infer the following estimates on $N^{t,j}_u$:
 \begin{align}
& \label{badiseqin31ass}
\|N^{t,1}_u \|_2 \lesssim \|u\|_3(\|u\|_4+\|u_t\|_2),  \\
& \label{badiseqin31new}\|N^{t,2}_u\|_0 \lesssim \|u\|^2_3+(\|u\|_3+\|u_{t}\|_1)\|u_{t}\|_2+\|u\|_3\|u_{tt}\|_0,\\
&\label{badiseqin3nes1}\|N^{t,3}_u\|_0
 \lesssim   \sqrt{\mathcal{E}^L} \mathcal{D}^H +\|u_t\|_3^2\lesssim   \sqrt{\mathcal{E}^L} \mathcal{D}^H.
 \end{align}
%\sqrt{\mathcal{E}^L} \mathcal{D}^H +\|u_t\|_2\|u_t\|_3.
 Then we can further deduce the following estimate.
 \begin{lem}\label{badiseqin}
 It holds that
 \begin{equation*}
 \frac{\mm{d}}{\mm{d}t} \|(\sqrt{\mu_1J} \nabla_{\ml{A}}   u_t,\sqrt{ \mu_2J} \mm{div}_{\mathcal{A}}u_t)\|^2_{0}+ c\|
 u_{tt}\|^2_{0}\lesssim  \|   u\|_2^2+\sqrt{\mathcal{E}^H}  \mathcal{D}^L. \end{equation*}
\end{lem}
\begin{pf} Multiplying \eqref{01s06p}$_2$ with $j=1$ by $J u_{tt}$, integrating (by parts) the resulting equality over $\Omega$
and using \eqref{AklJ=0}, we obtain
\begin{align}
  &\frac{1}{2}\frac{\mm{d}}{\mm{d}t}
  \int{J}(\mu_1 |\nabla_{\ml{A}}   u_t|^2
 +  \mu_2  |\mm{div}_{\mathcal{A}}u_t|^2)\mm{d}y+\int {\bar{\rho}}  |u_{tt}|^2\mm{d}y\nonumber\\
 & =
\int J(\nabla( P'(\bar{\rho})\mm{div}(\bar{\rho} u ) ) + g\mm{div} (\bar{\rho}u)e_3  +  \partial_t  \mathcal{L}_M )\cdot  u_{tt} \mm{d}y+  \int J  \mathcal{N}_t\cdot  u_{tt}\mm{d}y\nonumber\\
 & \quad+  \int  J N^{t,1}_u\cdot  u_{tt}\mm{d}y +{\mu_1}\int J\nabla_{\ml{A}}u_{t}: \nabla_{\ml{A}_t} u_{t}\mm{d}y +\frac{\mu_1}{2}\int J_t|\nabla_{\ml{A}}u_{t}|^2\mm{d}y \nonumber\\
  &\label{indeforut}\quad +{\mu_2}\int J \mm{div}_{\ml{A}}u_{t}: \mm{div}_{\ml{A}_t} u_{t}\mm{d}y+\frac{\mu_2}{2}\int J_t |\mm{div}_{\ml{A}}u_{t}|^2\mm{d}y=:\sum_{k=1}^7I^{L}_k.
 \end{align}
On the other hand, making use of \eqref{badiseqin31ass}, \eqref{N215562016302661}, \eqref{prtislsafdsfs}, \eqref{aimdse},
\eqref{Jtestismts} and \eqref{Jdetemrinat}, the seven integral terms $I^{L}_1$--$I^{L}_7$ can be bounded as follows.
 $$\begin{aligned}
%%  \label{ewqofJF1}
 & I^{L}_1\lesssim  \|u \|_{2}\| u_{tt}  \|_0,\quad I^{L}_2\lesssim \|  \mathcal{N}_t\|_0\|  u_{tt}\|_0
 \lesssim  \sqrt{ {\mathcal{E}^H}}\mathcal{D}^L, \quad I^{L}_3\lesssim \|N^{t,1}_u\|_0\| u_{tt}\|_0\lesssim \sqrt{\mathcal{E}^H} \mathcal{D}^L,\\
 %% \label{ewqofJF22110}
 & I^{L}_4\lesssim\|\nabla_{\ml{A}}u_{t}\|_0 \|\nabla_{\ml{A}_t} u_{t}\|_0 \lesssim\sqrt{\mathcal{E}^H} \mathcal{D}^L, \quad
 %% \label{ewqofJF2MH}
 I^{L}_5\lesssim \| J_t\|_2\|\nabla_{\ml{A}}u_{t} \|_0^2\lesssim  \sqrt{\mathcal{E}^H} \mathcal{D}^L, \\
 %% \label{ewqofJF2MHn}
 & I^{L}_6\lesssim \|\mm{div}_{\ml{A}}u_{t}\|_0 \| \mm{div}_{\ml{A}_t} u_{t}\|_0\lesssim \sqrt{\mathcal{E}^H} \mathcal{D}^L,\quad
%%  \label{ewqofJF2MHns}
 I^{L}_7\lesssim \|J_t\|_2 \|\mm{div}_{\ml{A}}u_{t}\|_0^2\lesssim\sqrt{\mathcal{E}^H} \mathcal{D}^L.
\end{aligned}$$
Thus, plugging the above seven estimates into \eqref{indeforut} and applying Cauchy-Schwarz's inequality, one obtains
Lemma \ref{badiseqin} immediately. \hfill$\Box$
\end{pf}

\subsection{$y_\mm{v}$-derivative estimates}
In this subsection we establish  the $y_\mm{v}$-derivative estimates. To this end, we rewrite \eqref{peturbationequation}
as the following non-homogeneous linear form:
\begin{equation}\label{s0106pnnnn}
\begin{cases}
  \eta_t=  u, \\
 \bar{\rho}J^{-1}u_t -\mu_1 \Delta u- \mu_2\nabla \mm{div} u -\nabla (P'(\bar{\rho})\mm{div}(\bar{\rho} \eta ))
  =  g \mm{div}(\bar{\rho}\eta) e_3+\mathcal{L}_M+\mathcal{N} +　 N_u,
\end{cases}  \end{equation}
where
$$N_u:=\mu_1 \left(\mm{div}_{\tilde{\mathcal{A}}}\nabla_{\tilde{\mathcal{A}}}u+\mm{div} \nabla_{\tilde{\mathcal{A}}}u
+\mm{div}_{\tilde{\mathcal{A}}}\nabla u\right) +\mu_2\left(\nabla_{\tilde{\mathcal{A}}}\mm{div}_{\tilde{\mathcal{A}}}u
+\nabla\mm{div}_{\tilde{\mathcal{A}}}u+\nabla_{\tilde{\mathcal{A}}}\mm{div} u\right). $$
Moreover, we employ \eqref{prtislsafdsfsfds}, \eqref{prtislsafdsfs}, \eqref{fgestims}, and the interpolation inequality to
control the term $N_u$ as follows.
\begin{align}  \label{lemm3601524Nu1612}&  \| N_u\|_3   \lesssim  \|\eta\|_3\|u\|_5 +\|\eta\|_5\|u\|_3 ,\\
  \label{Nu1743MHnew} &\| N_u\|_4   \lesssim   \|\eta\|_3 \|u\|_6  +\|\eta\|_6 \|u\|_3,  \\
  \label{Nu1743MH} &  \| N_u\|_5  \lesssim   \sqrt{\mathcal{E}^L}\|(\eta,u)\|_7+\|\eta\|_6\|u\|_4\lesssim \sqrt{\mathcal{E}^L}\|(\eta,u)\|_7,\\
 \label{Nu1743MHt}& \|\partial_t  N_u\|_{0}  \lesssim \sqrt{\mathcal{E}^H\mathcal{D}^L}, \quad  \|\partial_t^j N_u\|_{4-2j} \lesssim \sqrt{\mathcal{E}^L\mathcal{E}^H} ,\\
  \label{Nu1743MHtdfssd}& \|\partial_t^j N_u\|_{5-2j}  \lesssim \sqrt{\mathcal{E}^L\mathcal{D}^H}+
\|u\|_4^2+\|\eta\|_5\|u_t\|_3\lesssim \sqrt{\mathcal{E}^L\mathcal{D}^H},\quad j=1,2. \end{align}
Thus, we have the following bounds on $\eta$.
%%%%%%%%%%%%%%%%%%%%%%%%%%%%%%%%%%%%%%%%%
\begin{lem}  \label{ssebadiseqin}
It holds that
\begin{equation*}  \begin{aligned}
& \frac{\mm{d}}{\mm{d}t}\left( \int \bar{\rho}J^{-1} \partial_\mm{v}^j \eta \cdot  \partial_\mm{v}^j u\mm{d}y
+ \frac{1}{2} \|\partial_\mm{v}^j (\sqrt{\mu_1}\nabla\eta,\sqrt{\mu_2}\mm{div}\eta)\|_{0}^2\right) + c\|\partial_\mm{v}^j (\partial_1\eta ,\mm{div}\eta)\|_0^2\\
& \quad\lesssim \mm{sign}(j)(\|(\partial_1\eta ,\mm{div}\eta, u_t)\|_{\underline{j-1},0}^2)+\|  \partial_\mm{v}^j  u\|^2_{0}
 + \sqrt{\mathcal{E}^H} \mathcal{D}^L,\qquad j=0,\cdots,3, \end{aligned}
 \end{equation*}
where $\mm{sign}(j)=1$ when $j\neq 0$ and $=0$ when $j=0$.
\end{lem}
\begin{pf} Here we only prove the case $j=3$ and the rest three cases can be shown in the same manner.
Applying $\partial_\mm{v}^3$ to \eqref{s0106pnnnn}$_2$, multiplying the resulting equality by $\partial^3_\mm{v}\eta$,
and then using \eqref{s0106pnnnn}$_1$, we find that
$$\begin{aligned}
&\bar{\rho}J^{-1} \partial_t(\partial_\mm{v}^3\eta\cdot\partial_\mm{v}^3 u )-(\mu_1\Delta \partial_\mm{v}^3 \eta_t
+\mu_2\nabla \mm{div}\partial_\mm{v}^3 \eta_t+\partial_\mm{v}^3 \nabla (P'(\bar{\rho})\mm{div}(\bar{\rho}\eta ) )\cdot \partial_\mm{v}^3\eta\\
& =\partial_\mm{v}^3(  g\mm{div}(\bar{\rho}\eta) e_3 + \mathcal{L}_M+ \mathcal{N} + N_u)\cdot\partial_\mm{v}^3\eta +\bar{\rho}J^{-1}|\partial_\mm{v}^3  u|^2-\sum_{ k=0}^2 (\partial_\mm{v}^{3-k}(\bar{\rho}J^{-1})\partial_\mm{v}^k u_t)\cdot \partial_\mm{v}^3 \eta.\end{aligned}$$
Integrating (by parts) the above  identity over $\Omega$, we have
\begin{equation}
\label{estimforhoedsds1}
\begin{aligned}
&\frac{\mm{d}}{\mm{d}t}\left(\int \bar{\rho} J^{-1}\partial_\mm{v}^3 \eta \cdot \partial_\mm{v}^3 u \mm{d}y+\frac{\mu_1}{2} \int|\nabla \partial_\mm{v}^3 \eta|^2\mm{d}y+ \frac{\mu_2}{2}\int_\Omega |\mm{div} \partial_\mm{v}^3 \eta|^2\mm{d}y\right)  \\
&=\int \partial_\mm{v}^3( g\mm{div}(\bar{\rho}\eta))\partial_\mm{v}^3\eta_3 \mm{d}y
-\int \partial_\mm{v}^3(P'(\bar{\rho})\mm{div}(\bar{\rho} \eta ) )\mm{div} \partial_\mm{v}^3\eta \mm{d}y+
\int  \partial_\mm{v}^3 \mathcal{L}_M\cdot  \partial_\mm{v}^3\eta \mm{d}y  \\
  &\quad+\int  \partial_\mm{v}^3\mathcal{N} \cdot  \partial_\mm{v}^3\eta\mm{d}y+\int\partial_\mm{v}^3 N_u\cdot\partial_\mm{v}^3\eta \mm{d}y
  -\sum_{ k=0}^2\int  (\partial_\mm{v}^{3}(\bar{\rho}J^{-1}  u_t)- \bar{\rho}J^{-1} \partial_\mm{v}^3 u_t)\cdot \partial_\mm{v}^3 \eta \mm{d}y  \\
 &\quad+\int \bar{\rho} J^{-1}_t\partial_\mm{v}^3 \eta \cdot \partial_\mm{v}^3 u \mm{d}y +\int {\bar{\rho}J^{-1}}|\partial_\mm{v}^3 u|^2\mm{d}y\leq \sum_{k=1}^7 J^{L}_{k}+c\|  \partial_\mm{v}^3  u\|^2_{0},
\end{aligned}
\end{equation}
 where the first seven integrals on the right hand of \eqref{estimforhoedsds1} are denoted by $J^L_1$--$J^L_7$, respectively.

In view of \eqref{modififorg} and \eqref{rho20160326}, we see that $J^L_1$ can be bounded as follows.
\begin{equation}  \label{estrhoe}     \begin{aligned}
 J^{L}_{1}= & \int g\bar{\rho}'|\partial_\mm{v}^3\eta_3|^2\mm{d}y + \int g\bar{\rho}\partial_\mm{v}^3\eta_3\mm{div}\partial_\mm{v}^3\eta\mm{d}y \\
& +\sum_{k=0}^2\int\partial_\mm{v}^3\eta_3\left((\partial_\mm{v}^{3}( g\bar{\rho}'\eta_3)-g\bar{\rho}'\partial_\mm{v}^3\eta_3)
+ (\partial_\mm{v}^{3}(g\bar{\rho}\mm{div}  \eta) -g\bar{\rho}\mm{div} \partial_\mm{v}^3 \eta )\right)\mm{d}y   \\
 \leq & \int g\bar{\rho}'|\partial_\mm{v}^3\eta_3|^2 \mm{d}y + \int g\bar{\rho}\partial_\mm{v}^3\eta_3\mm{div}\partial_\mm{v}^3 \eta\mm{d}y
 + c \|\partial_\mm{v}^3\eta_3\|_0 \| (\eta ,\mm{div}\eta)\|_{\underline{2},0},
 \end{aligned}\end{equation}
%%  Similarly to \eqref{estrhoe}, one has
\begin{equation}   \label{estrhoeMH}
 \begin{aligned}
J^{L}_{2}= &  - \int P'(\bar{\rho})\bar{\rho} | \mm{div}\partial_\mm{v}^3\eta |^2\mm{d}y - \int \bar{P}'
 \partial_\mm{v}^3\eta_3 \mm{div}\partial_\mm{v}^3 \eta\mm{d}y    \\
   &  -\sum_{ k=0}^2 \int \left(\partial_\mm{v}^{3}(P'(\bar{\rho})\bar{\rho} \mm{div}\eta)-P'(\bar{\rho})\bar{\rho} \mm{div}\partial_\mm{v}^3\eta
+\partial_\mm{v}^{3}( \bar{P}'\eta_3)- \bar{P}' \partial_\mm{v}^3\eta_3 \right)\mm{div}\partial_\mm{v}^3\eta  \mm{d}y　 \\
   \leq& - \int {P'(\bar{\rho})\bar{\rho}} |\mm{div}\partial_\mm{v}^3\eta |^2\mm{d}y -\int \bar{P}'
 \partial_\mm{v}^3\eta_3  \mm{div}\partial_\mm{v}^3 \eta\mm{d}y
 + c \| \mm{div}\partial_\mm{v}^3 \eta\|_0 \|(\eta,\mm{div}\eta )\|_{\underline{2},0}
 \end{aligned}\end{equation}
 and
 \begin{equation}  \label{201601292111}
\begin{split}    \begin{aligned}
 J^{L}_{3}= &\lambda \int (m^2\partial_1^2\partial_\mm{v}^3 \eta -m^2\partial_1\mm{div}\partial_\mm{v}^3 \eta e_1
  +\nabla (m m '\partial_\mm{v}^3 \eta_3+m^2  \mm{div}_\mm{v}\partial_\mm{v}^3 \eta_\mm{v}  )\cdot \partial_\mm{v}^3 \eta\mm{d}y\\
 &-\lambda \sum_{k=0}^2\int (\partial_\mm{v}^{3}( m^2(\partial_1  \eta
 - \mm{div} \eta e_1))-  m^2\partial_\mm{v}^3(\partial_1 \eta
 - \mm{div} \eta e_1 ))\cdot \partial_1\partial_\mm{v}^3 \eta    \\
  &\qquad \qquad + (\partial_\mm{v}^{3} (m m ' \eta_3 )- m m '\partial_\mm{v}^3 \eta_3
  +  \partial_\mm{v}^{3}(m^2  \mm{div}_\mm{v} \eta_\mm{v})  -m^2  \mm{div}_\mm{v}\partial_\mm{v}^3 \eta_\mm{v}   )\mm{div} \partial_\mm{v}^3 \eta)\mm{d}y    \\
 \leq  &-\lambda\int m^2(|\partial_1 \partial_\mm{v}^3 \eta_\mm{v}|^2+
 |\mm{div}_\mm{v}\partial_\mm{v}^3 \eta_\mm{v}|^2) \mm{d}y-\lambda\int   m m '\partial_\mm{v}^3\eta_3\mm{div}\partial_\mm{v}^3\eta\mm{d}y\\
&+ c \| \partial_\mm{v}^3(\partial_1\eta,\mm{div} \eta)\|_0 \|(\eta,\partial_1\eta ,\mm{div}\eta)\|_{\underline{2},0}.
 \end{aligned}   \end{split}\end{equation}
Inserting the above three estimates into \eqref{estimforhoedsds1}, we deduce by \eqref{insoslai} and \eqref{comsteady} that
\begin{equation}  \label{estimforhoedsds1st}
\begin{aligned}
&\frac{\mm{d}}{\mm{d}t}\left(\int \bar{\rho} J^{-1}\partial_\mm{v}^3 \eta \cdot \partial_\mm{v}^3 u \mm{d}y+\frac{1}{2}\|\partial_\mm{v}^3
(\sqrt{\mu_1}\nabla   \eta,\sqrt{\mu_2}\mm{div} \eta)\|^2_0\right)  -  E (\partial_\mm{v}^3 \eta)\\
& \quad \leq c\left(\|  \partial_\mm{v}^3 u\|^2_{0}+  \| \partial_\mm{v}^3( \partial_1  \eta,\mm{div}\eta)\|_{0} 　
\|(\partial_1\eta ,\mm{div}\eta)\|_{\underline{2},0}\right)+\sum_{k=4}^7 J^{L}_{k}.
\end{aligned}   \end{equation}

On the other hand, by \eqref{lemm3601524Nu1612}, \eqref{N21556201630266},  \eqref{Jtestismtsinver}, \eqref{Jdetemrinat2nes} and \eqref{esmmdforinfty},
we have
 \begin{align}
& \label{lemm3601524}
 J^L_4  \lesssim \|\partial_\mm{v}^3\mathcal{N}\|_0\|\partial_\mm{v}^3\eta 　\|_0 \lesssim \sqrt{\mathcal{E}^H} \mathcal{D}^L,\\
 &J^L_5 \lesssim\|\partial_\mm{v}^3 N_u\|_0\| \partial_\mm{v}^3\eta \|_0 \lesssim  \sqrt{\mathcal{E}^H} \mathcal{D}^L,   \label{lemm3601524+1} \\
 &J^L_6 \lesssim  (\|J^{-1}-1\|_{5}+1)\| u_t\|_{\underline{2},0}\| \partial_\mm{v}^3 \eta\|_0
 \lesssim \| u_t\|_{\underline{2},0}\| \partial_\mm{v}^3 \eta\|_0,    \label{lemm3601524+2} \\
 &\label{lemm36015241}J^L_7 \lesssim \|J^{-1}_t\|_2\|\partial_\mm{v}^3 \eta\|_0\| \partial_\mm{v}^3 u\|_0 \lesssim  \sqrt{\mathcal{E}^H} \mathcal{D}^L.
\end{align}
Consequently, putting the above four estimates  into \eqref{estimforhoedsds1st}, and using Lemma \ref{lem:0601}, \eqref{insoslai} and Cauchy-Schwarz's inequality,
we deduce the desired conclusion for the case $j=3$.   \hfill$\Box$
\end{pf}

Similarly, we can also establish the $y_\mm{v}$-derivative estimates on $u$.
\begin{lem}   \label{ssebadsdiseqinsd}
It holds that
$$ \frac{\mm{d}}{\mm{d}t}\left(\|\sqrt{ \bar{\rho}J^{-1} }\partial_\mm{v}^j u\|^2_0-E (\partial_{\mm{v}}^j\eta)\right)
+ c\| \nabla\partial_\mm{v}^j u\|_{0}^2 \lesssim \mm{sign}(j)\| ( \partial_1\eta ,\mm{div}\eta,u_t)\|_{\underline{j-1},0}^2
+\sqrt{\mathcal{E}^H}\mathcal{D}^L $$
for $j=0,1,2,3$.
\end{lem}
\begin{pf} We only show the case $j=3$ and the rest cases can be proved in the same manner.
We apply $\partial_\mm{v}^3$ to \eqref{s0106pnnnn}$_2$ and multiply the resulting equation by $\partial^3_\mm{v} u$ to get
$$ \begin{aligned}
& {\bar{\rho}J^{-1}} \partial_\mm{v}^3u_t\cdot\partial_\mm{v}^3u -(\mu_1\Delta \partial_\mm{v}^3 u +\mu_2\nabla\mm{div}\partial_\mm{v}^3 u
+\partial_\mm{v}^3\nabla (P'(\bar{\rho})\mm{div}(\bar{\rho} \eta ))\cdot \partial_\mm{v}^3 u  \\
& =\partial_\mm{v}^3 (g\mm{div}(\bar{\rho}\eta) e_3 +\mathcal{L}_M +\mathcal{N} +N_u)\cdot\partial_\mm{v}^3 u -\sum_{k=0}^2
  \partial_\mm{v}^{3-k}(\bar{\rho}J^{-1})\partial_\mm{v}^k u_t\cdot\partial_\mm{v}^3 u. \end{aligned}$$
Integrating (by parts) the above  identity over $\Omega$, one has
\begin{align}
&\frac{1}{2}\frac{\mm{d}}{\mm{d}t} \int
 \bar{\rho} J^{-1}|\partial_\mm{v}^3 u|^2\mm{d}y +{\mu_1}\int|\nabla \partial_\mm{v}^3 u|^2\mm{d}y+{\mu_2}\int|\mm{div} \partial_\mm{v}^3 u|^2\mm{d}y\nonumber \\
&= \int \partial_\mm{v}^3 (g\mm{div}(\bar{\rho}\eta))\partial_\mm{v}^3 u_3\mm{d}y
-\int \partial_\mm{v}^3  (P'(\bar{\rho})\mm{div}(\bar{\rho} \eta ))
\mm{div}\partial_{\mm{v}}^3u\mm{d}y+\int \partial_\mm{v}^3\mathcal{L}_M\cdot \partial_\mm{v}^3 u\mm{d}y\nonumber \\
&\qquad+\int \partial_\mm{v}^3 \mathcal{N} \cdot \partial_\mm{v}^3 u\mm{d}y+　 \int \partial_\mm{v}^3 N_u\cdot \partial_\mm{v}^3 u　\mm{d}y
+\frac{1}{2}\int \bar{\rho} J^{-1}_t|\partial_\mm{v}^3 u|^2 \mm{d}y　\nonumber \\
&\label{estimforhoe1nwe}\qquad-\sum_{ k=0}^2\int
 \partial_\mm{v}^{3-k}(\bar{\rho}J^{-1})\partial_\mm{v}^k u_t\cdot \partial_\mm{v}^3 u\mm{d}y =:\sum_{k=1}^7 K^{L}_{k}.
\end{align}

Similarly to \eqref{estrhoe}--\eqref{201601292111}, the three integrals $K^L_1$--$K^L_3$ can be controlled as follows.
 \begin{align}
&\label{estrhoenn201604111600n}
 \begin{aligned}
 K^{L}_{1}\leq &\frac{1}{2} \frac{\mm{d}}{\mm{d}t}\left(\int g\bar{\rho}'|\partial_\mm{v}^3\eta_3|^2 \mm{d}y+    2\int g\bar{\rho}  \partial_\mm{v}^3  \eta_3\mm{div}\partial_\mm{v}^3 \eta\mm{d}y\right) - \int g\bar{\rho}  \partial_\mm{v}^3  \eta_3\mm{div}\partial_\mm{v}^3u\mm{d}y \\
 &+ c \| \partial_\mm{v}^3 u_3\|_{0} \|  ( \eta ,\mm{div}\eta)\|_{\underline{2},0}, \end{aligned} \\[1mm]
& \begin{aligned}    K^{L}_{2}\leq & - \frac{1}{2} \frac{\mm{d}}{\mm{d}t} \int P'(\bar{\rho})\bar{\rho}
| \mm{div}\partial_\mm{v}^3\eta |^2\mm{d}y   - \int   \bar{P}'
 \partial_\mm{v}^3\eta_3  \mm{div}\partial_\mm{v}^3 u\mm{d}y   \\
 &+ c \| \mm{div} \partial_\mm{v}^3 u\|_{0} \| (\eta,\mm{div}\eta)\|_{\underline{2},0}
\end{aligned}\end{align}
and
\begin{equation}   \label{estrhoenn201604111601n}
 \begin{aligned}
 K^{L}_{3}\leq &-\lambda\int\Big( (m^2\partial_1 \partial_{\mm{v}}^3\eta
-m^2\mm{div}\partial_{\mm{v}}^3\eta  e_1)\cdot \partial_1\partial_{\mm{v}}^3 u \\
& \qquad \quad  \ +(m^2 \mm{div}\partial_{\mm{v}}^3\eta  + m m '\partial_{\mm{v}}^3 \eta_3 -m^2 \partial_1\partial_{\mm{v}}^3 \eta_1 )\mm{div}\partial_{\mm{v}}^3u\Big)\mm{d}y\\
 & +c\| \partial_\mm{v}^3(\partial_1  u,\mm{div}  u)\|_{0} \| ( \eta,\partial_1\eta ,\mm{div}\eta)\|_{\underline{2},0}\\
 \leq &-\frac{\lambda}{2}\frac{\mm{d}}{\mm{d}t}\int m^2 (|\partial_1\partial_{\mm{v}}^3 \eta_\mm{v}|^2
+ |\mm{div}_{\mm{v}}\partial_{\mm{v}}^3\eta_{\mm{v}} |^2
  )\mm{d}y -\lambda\int m m '\partial_{\mm{v}}^3\eta_3 \mm{div}\partial_{\mm{v}}^3u\mm{d}y\\
 &  +c\| \partial_\mm{v}^3(\partial_1  u,\mm{div}  u)\|_{0} \|  ( \eta,\partial_1\eta ,\mm{div}\eta)\|_{\underline{2},0}.
\end{aligned}
\end{equation}
Substituting the above three estimates into \eqref{estimforhoe1nwe}, using \eqref{insoslai} and \eqref{comsteady}, we conclude that
\begin{equation}
\label{estimforho4e1nwe}
 \frac{1}{2}\frac{\mm{d}}{\mm{d}t}\left(\|\sqrt{ \bar{\rho}J^{-1} }\partial_\mm{v}^3 u\|^2_0-E (\partial_{\mm{v}}^3\eta)\right) + \|\partial_\mm{v}^3(\sqrt{\mu_1}\nabla u, \sqrt{\mu_2}\mm{div}  u)\|_{0}^2
 \leq c  \|  ( \partial_1\eta ,\mm{div}\eta)\|_{\underline{2},0}^2+\sum_{k=4}^7 K^{L}_{k}.  \end{equation}

On the other hand, similarly to \eqref{lemm3601524}--\eqref{lemm36015241}, one obtains
\begin{align*}\label{201604031210}
& \sum_{k=4}^7 K^{L}_{k}\lesssim\| u_t\|_{\underline{2},0}\| \partial_\mm{v}^3 u\|_0+\sqrt{\mathcal{E}^H} \mathcal{D}^L.
\end{align*}
Consequently, inserting the above inequality into \eqref{estimforho4e1nwe}, and using
\eqref{insoslai} and Cauchy-Schwarz's inequality, we obtain Lemma \ref{ssebadsdiseqinsd}.     \hfill$\Box$
\end{pf}

\subsection{$y$-derivative estimates}
    In this subsection we derive the $y$-derivative estimates of $\eta$ and $u$. Similarly to the incompressible magnetic RT problem
 see \cite[Lemma 2.4]{JFJSJMFMOSERT}), we focus on the term $\partial_1^2\eta$ and extract the following equations from
the second and third components of \eqref{s0106pnnnn}$_2$:
\begin{equation}   \label{Stokesequson}
\begin{aligned}
& \!\!\! -\mu_1 \partial_1^2 \partial_t \eta_\mm{v}-\lambda m^2\partial_1^2\eta_\mm{v}
=\mu_1 \Delta_\mm{v} u_\mm{v} + \mu_2\nabla_\mm{v} \mm{div} u  \\
& +\nabla_\mm{v} (P'(\bar{\rho})\mm{div}(\bar{\rho} \eta )+\lambda m m '\eta_3+\lambda m^2   \mm{div}_\mm{v}\eta_\mm{v})
 +g \mm{div}(\bar{\rho}\eta)(0,1)^{\mm{T}}-\bar{\rho}J^{-1}\partial_t u_\mm{v}+{N}_\mm{v} +　 N_{u,\mm{v}}.
\end{aligned}\end{equation}
Meanwhile, the first component of \eqref{s0106pnnnn}$_2$ can be written as
\begin{equation}  \label{Stokesequson1}
\begin{aligned} &- \mu \partial_1 \mm{div}\eta_t - P'(\bar{\rho}) \bar{\rho}\partial_1 \mm{div}\eta 　\\
& \quad = \mu_1 \Delta_\mm{v} u_1 +  \partial_1 (\bar{P}' \eta_3  + \lambda m m' \eta_3
-\mu_1 \mm{div}_\mm{v} u_\mm{v}) -\bar{\rho}J^{-1}\partial_t u_1 + {N}_1 +　 N_{u,1},
\end{aligned}\end{equation}
where $\mu:=\mu_1+ \mu_2$, $\mathcal{N}:=( {N}_1,{N}_\mm{v})$ and $N_{u}:=(N_{u,1},N_{u,\mm{v}})$.

Noting that the order of $\partial_1$ in the linear part on the right hand side of \eqref{Stokesequson} is lower than that on the left hand side,
this feature provides a possibility that the $y_1$-derivative estimates of $\eta_\mm{v}$ can be converted to the $y_{\mm{v}}$-derivative
estimates of $\eta$. Similarly, \eqref{Stokesequson1} also provides a possibility that the $y_1$-derivative estimates of $\eta_1$ can be converted
to the the $y_{\mm{v}}$-derivative estimates of $\eta$ by the relation $\mm{div}\eta=\partial_1\mm{\eta}_1+\mm{div}_{\mm{v}}\eta_{\mm{v}}$.
Based on these basic observations, we can establish the following $y$-derivative estimates on $\eta$.
%%%%%%%%%%%%%%%%%%%%%%%%%%%%%%%%%%%%%%%%%%%%%%%%%%%%%%%%%%%%%%%
\begin{lem}   \label{lem:dfifessim}  We have
 \begin{equation*}
\frac{\mm{d}}{\mm{d}t}\mathcal{H}_2(\eta) + \|  (\partial_1\eta,\mm{div}\eta)\|_{3}^2+\|u\|_{4}^2
 \lesssim \|(\partial_1 \eta,\mm{div}\eta,\nabla u)\|_{\underline{3},0}^2+ \| u_t\|_{2}^2+ {\mathcal{E}^H}\mathcal{D}^L,   \end{equation*}
 where the energy functional $\mathcal{H}_2(\eta)$ satisfies
\begin{equation}\label{Hetaprofo}
 \|\partial_1^2\eta \|_{2}^2- h_2\| \partial_1 \eta_{\mm{v}}\|_{\underline{3},0}^2\lesssim {\mathcal{H}}_{2}(\eta) \end{equation}
for some positive constant $h_2$, depending on $\Omega$ and other physical parameters.
\end{lem}
\begin{pf}  Applying $\|\cdot \|_{j,i-j}^2$ to \eqref{Stokesequson1}, using \eqref{Jdetemrinatnesw}, \eqref{fgestims} and Cauchy--Schwarz's inequality, we deduce that for $0\leq i\leq 5$,
\begin{equation*}    \begin{aligned}
 &   \frac{\mm{d}}{\mm{d}t}\|  \sqrt{\mu P'(\bar{\rho}) \bar{\rho}}\diamond\partial_1 \mm{div}\eta
\|_{j,i-j}^2+(\| ({P'(\bar{\rho})\bar{\rho}})\diamond\partial_1\mm{div}\eta \|_{j,i-j}^2 +\| {\mu}\partial_1\mm{div}u\|_{j,i-j}^2 )　\\
&\lesssim \left\|  \mu_1 \Delta_\mm{v} u_1 +  \partial_1 (\bar{P}' \eta_3  + \lambda m m ' \eta_3-\mu_1 \mm{div}_\mm{v} u_\mm{v})-\bar{\rho}J^{-1}\partial_t u_1+{N}_1 +　 N_{u,1}\right\|_{j,i-j}^2+\|\partial_1\eta\|_{j,i-j}^2
\\
&\lesssim \|  \partial_1 \eta  \|_{\underline{j},i-j}^2  +  \|  u \|_{j+1,i-j+1}^2 +\| ( u_t,\mathcal{N},N_{u }) \|_{i}^2,
\end{aligned}\end{equation*}
which, together with \eqref{divfefdMH} and \eqref{0102}, gives
\begin{equation}   \label{201602011}
\begin{aligned}
& \frac{\mm{d}}{\mm{d}t}\| \sqrt{\mu P'(\bar{\rho}) \bar{\rho}}\diamond\partial_1 \mm{div}\eta\|_{j,i-j}^2
+c(\|\partial_1^2 (\eta_1,u_1)\|_{j,i-j}^2 +\|\partial_1 (\mm{div}\eta,\mm{div}u)\|_{j,i-j}^2 )  　\\
&\lesssim \|\partial_1\eta\|_{\underline{j+1},i-j}^2 +\|u\|_{j+1,i-j+1}^2 +\|(u_t,\mathcal{N},N_u)\|_{i}^2.
\end{aligned}\end{equation}

Similarly, one easily deduces from \eqref{Stokesequson} that
\begin{equation}    \label{20160201}
\begin{aligned}
& \frac{\mm{d}}{\mm{d}t}\|\sqrt{\mu_1\lambda}m\partial_1^2 \eta_\mm{v}\|_{j,i-j}^2+c\|\partial_1^2(\eta_\mm{v},u_\mm{v})\|_{j,i-j}^2 \\
& \lesssim \|\mu_1 \Delta_\mm{v} u_\mm{v} + \mu_2\nabla_\mm{v} \mm{div} u
+\nabla_\mm{v} (P'(\bar{\rho})\mm{div}(\bar{\rho} \eta )+\lambda m m '\eta_3+\lambda m^2   \mm{div}_\mm{v}\eta_\mm{v}) \\
& \qquad +g \mm{div}(\bar{\rho}\eta)(0,1)^{\mm{T}}-\bar{\rho}J^{-1}\partial_t u_\mm{v}+{N}_\mm{v} +　 N_{u,\mm{v}}\|_{j,i-j}^2\\
& \lesssim \|(\eta,\mm{div}_{\mm{v}}\eta_{\mm{v}},\mm{div}\eta)\|^2_{\underline{j+1},i-j} +\|u\|_{j+1,i-j+1}^2 +\|(u_t,\mathcal{N},N_{u })\|_{i}^2.
\end{aligned}\end{equation}
Thus, using \eqref{divfefdMH} and \eqref{insoslai}, we infer from \eqref{201602011} and \eqref{20160201} that
\begin{equation} \label{omdm12dfs2MJ}\begin{aligned}
 &  \frac{\mm{d}}{\mm{d}t} (\|\sqrt{\mu_1 \lambda} m\partial_1^2 \eta_{\mm{v}}\|_{j,i-j}^2
 +\| \sqrt{\mu P'(\bar{\rho}) \bar{\rho}}\diamond\partial_1 \mm{div}\eta)
\|_{j,i-j}^2 )\\
 &\quad + c(\|(\partial_1\eta,\mm{div}\eta)\|_{j ,i-j+1}^2+\|u\|_{j ,i-j+2}^2 )\lesssim   \|(\partial_1 \eta,\mm{div}\eta,\nabla u)\|_{\underline{j+1},i-j}^2+ \| ( u_t,\mathcal{N},N_{u }) \|_{i}^2.
 \end{aligned}
 \end{equation}
Adding up \eqref{omdm12dfs2MJ} from  $0$ to $j$, one gets
\begin{equation}
\label{20101281537} \begin{aligned}
 &  \frac{\mm{d}}{\mm{d}t} (\|\sqrt{\mu_1 \lambda} m\partial_1^2 \eta_{\mm{v}}\|_{\underline{j},i-j}^2+\| \sqrt{\mu P'(\bar{\rho}) \bar{\rho}}\diamond\partial_1 \mm{div}\eta
\|_{\underline{j},i-j}^2 )\\
 &\quad + c (\|(\partial_1\eta, \mm{div}\eta )\|_{\underline{j},i-j+1}^2+\|u\|_{\underline{j} ,i-j+2}^2 )\lesssim   \|(\partial_1\eta , \mm{div}\eta, \nabla u)\|_{\underline{j+1},i-j}^2+\| ( u_t,\mathcal{N},N_{u }) \|_{i}^2.
 \end{aligned}
 \end{equation}

In addition, we can use the relation \eqref{divfefdMH}
and the estimate
$$ \|\partial_1\eta_{\mm{v}}\|_{{j+1},i-j}^2\lesssim \|\partial_1^2\eta_{\mm{v}}\|_{{j+1},i-j-1}^2
 +\| \partial_1 \eta_{\mm{v}}\|_{\underline{i+1},0}^2\quad\mbox{for }0\leq j<i, $$
 to find that
\begin{equation}
\label{20101281538}
\begin{aligned}
&\|\sqrt{\mu_1 \lambda} m\partial_1^2 \eta_{\mm{v}}\|_{\underline{i},0}^2
+\|\sqrt{\mu P'(\bar{\rho}) \bar{\rho}}\diamond\partial_1 \mm{div}\eta\|_{\underline{i},0}^2
 \\
 &\gtrsim \| (\partial_1^2 \eta_{\mm{v}},\partial_1 \mm{div}\eta )\|_{i,0}^2\gtrsim\|\partial_1^2\eta \|_{i,0}^2
-\tilde{h}_{i,i}\| \partial_1 \eta_{\mm{v}}\|_{{i+1},0}^2
\end{aligned}
\end{equation}
and
\begin{equation}\label{20101281539}\begin{aligned} &\| \sqrt{\mu_1 \lambda} m\partial_1^2
\eta_{\mm{v}}\|_{\underline{j},i-j}^2
+\|\sqrt{\mu P'(\bar{\rho}) \bar{\rho}}\diamond\partial_1 \mm{div}\eta )\|_{\underline{j},i-j}^2
\gtrsim \|(\partial_1^2\eta_{\mm{v}},\partial_1\mm{div}\eta \|_{j,i-j}^2 \\
&\gtrsim\| \partial_1^2\eta \|_{j,i-j}^2-\tilde{h}_{i,j}\left( \| \partial_1^2 \eta\|_{{j+1},i-j-1}^2
+\| \partial_1 \eta_{\mm{v}}\|_{\underline{i+1},0}^2\right),\quad 0\leq j<i \end{aligned}\end{equation}
for some positive constants $\tilde{h}_{i,j}$.
Thus, we can further deduce from \eqref{20101281537}--\eqref{20101281539} that
\begin{equation}
\label{omdm12dfs2}
   \frac{\mm{d}}{\mm{d}t}{\mathcal{H}}_{i}(\eta)+ \|(\partial_1\eta, \mm{div}\eta)\|_{i+1}^2+\|u\|_{i+2}^2
 \lesssim \|(\partial_1\eta ,\mm{div}\eta ,\nabla u)\|^2_{\underline{i+1},0} + \| ( u_t,\mathcal{N},N_u)\|_{i}^2,
  \end{equation}
 where  ${\mathcal{H}}_i(\eta):=\sum_{j=0}^i {h}_{i,j}( \|\sqrt{\mu_1 \lambda} m\partial_1^2 \eta_{\mm{v}}\|_{\underline{j},i-j}^2+\|\sqrt{\mu P'(\bar{\rho})
 \bar{\rho}}\diamond\partial_1 \mm{div}\eta \|_{\underline{j},i-j}^2)$
satisfies
$$  \| \partial_1^2\eta\|_{i}^2- h_i\| \partial_1 \eta_{\mm{v}}\|_{\underline{i+1},0}^2\lesssim {\mathcal{H}}_{i}(\eta) $$
for some positive constants ${h}_{i,j}$ and $h_i$.

If we take $i=2$ in \eqref{omdm12dfs2}, we get
 \begin{equation}\label{dfifessim16431}
\frac{\mm{d}}{\mm{d}t}\mathcal{H}_2(\eta) + \|(\partial_1\eta,\mm{div}\eta)\|_{3}^2+ \|u\|_{4}^2
\lesssim \|(\partial_1 \eta,\mm{div}\eta,\nabla u)\|_{\underline{3},0}^2 +
\|u_t\|_{2}^2 +\|\mathcal{N}\|_{2}^2 +\|N_{u}\|_{2}^2.
 \end{equation}
On the other hand, by virtue of \eqref{lemm3601524Nu1612} and \eqref{N21556201630266},
$$ \|\mathcal{N} \|_{2}^2  +  \|N_{u }\|_{2}^2 \lesssim {\mathcal{E}^H  {\mathcal{D}}^L}, $$
which combined with \eqref{dfifessim16431} gives the lemma.   \hfill$\Box$
\end{pf}

Now, we turn to the derivation of $y$-derivative estimates of $u$. If we apply $\partial_t^j$ to \eqref{s0106pnnnn}$_2$, we obtain
$$  \begin{cases}-\mu_1 \Delta \partial_{t}^j u- \mu_2\nabla \mm{div} \partial_{t}^j u = f_j ,\\[1mm]
 \partial_{t}^ju|_{\partial\Omega}=0,  \end{cases}$$
where
\begin{align}
f_j:=&\partial_{t}^j(\nabla (P'(\bar{\rho})\mm{div}(\bar{\rho}\eta ))
 +   g \mm{div} (\bar{\rho} \eta) e_3+\mathcal{L}_M + \mathcal{N} + N_u)
\nonumber  \\
 &-\bar{\rho} J^{-1} \partial_{t}^{j+1}u-\bar{\rho} \sum_{ k=0}^{j-1}C_j^k\partial_{t}^{j-k} J^{-1}\partial_t^{k+1}u.\nonumber
 \end{align}
Hence, one can apply the classical elliptic regularity theory to the above Lam\'e system \cite{KYLH3295,USTICYCHJHJJ} to get
\begin{equation}\label{02181}\begin{aligned}
\|\partial_t^j u\|_{i-2j+2}^2  \lesssim &\|f_j\|_{i-2j}  \\
\lesssim&\|\partial_{t}^j (\nabla (P'(\bar{\rho})\mm{div}(\bar{\rho} \eta ))
+g \mm{div} (\bar{\rho} \eta) e_3+\mathcal{L}_M) -\bar{\rho}J^{-1}\partial_{t}^{j+1}u\|_{i-2j}^2 +S_{j,i},
 \end{aligned}   \end{equation}
where $S_{0,i}:=\|(\mathcal{N},N_u)\|_{i}^2$ and
$$ S_{j,i}:=\|\partial_t^j (\mathcal{N},N_u)\|_{i-2j}^2 +\Big\| \bar{\rho}
\sum_{k=0}^{j-1}\partial_{t}^{j-k} J^{-1}\partial_t^{k+1}u \Big\|_{i-2j}^2\mbox{ for }j\geq 1. $$

% We mention that the regularity theory for the Lam\'e system in different domains a bounded domain %in $\mathbb{R}^3$ with smooth boundary, or
% in the whole space $\mathbb{R}^3$, or in the half space $\mathbb{R}^3_+$, or in an exterior domain %in $\mathbb{R}^3$, can be founded in \cite[Section 5]{USTICYCHJHJJ}.
%{If one applies the regularity theory for the Lam\'e system in the domain $\mathbb{R}^3$ or %$\mathbb{R}^3_+$ \cite[Section 5]{USTICYCHJHJJ}, one can obtain the
%same regularity for the Lam\'e system defined in a strip domain by following the derivation of (3.7) %in \cite{KYLH3295}.}
%Based on \eqref{02181}, we can establish the following estimates.
 %%%%%%%%%%%%%%%%%%%%%%%%%%%%%%%%%%%%%%%%%%%%%%%%%%%%%%%%%%%%%%%%%%%%%%%
\begin{lem}\label{lem:dfifessimell} We have
\begin{align}\label{dfifessim} &   \|u\|_{3}^2 \lesssim\|\eta\|_{3}^2 +\|u_t\|_{1}^2 +{\ml{E}^H}\mathcal{E}^L, \\
 &\label{dfifeddssimlas}\|u_t\|_{2}^2  \lesssim  \|u\|_2^2+\|u_{tt}\|_{0}^2  + {\ml{E}^H}\ml{D}^L.
 \end{align}  \end{lem}
\begin{pf}  In view of \eqref{02181} with $(i,j)=(1,0)$, one has
$$  \|u\|_{3}^2  \lesssim   \|  \eta\|_{3}^2 +\|u_t\|_{1}^2  +S_{0,1}.   $$
On the other hand, it follows from \eqref{lemm3601524Nu1612} and \eqref{N21556201630266} that
$$ S_{0,1}= \| (\mathcal{N},N_u)\|_{1}^2  \lesssim  {\mathcal{E}^H { \mathcal{E}}^L }.  $$
Therefore, we obtain \eqref{dfifessim} from the above two estimates.

To show \eqref{dfifeddssimlas}, we take $(i,j)=(2,1)$ in \eqref{02181} to find that
$$ \| u_t\|_{2}^2 \lesssim \|u\|_2^2+\|u_{tt}\|_{0}^2 +S_{1,2}. $$
Utilizing \eqref{Nu1743MHt}, \eqref{N215562016302661}, \eqref{Jtestismtsinver} and  \eqref{fgestims}, we have
$$ S_{1,2}= \|\partial_t ( \mathcal{N}, N_u)\|_{0}^2+\left\|\bar{\rho}   J^{-1}_tu_t\right\|_{0}^2 \lesssim {\mathcal{E}^H{\mathcal{D}}^L }. $$
The above two estimates give then \eqref{dfifeddssimlas}. The proof is complete. \hfill $\Box$
\end{pf}

\subsection{Lower-order energy inequality}
Now we are ready to build the lower-order energy inequality. In what follows the letters $c_{i}^L$, $i=1,\ldots,9$, will denote generic positive constants
which may depend on $\Omega$ and the physical parameters in the perturbation equations.
%%%%%%%%%%%%%%%%%%%%%%%%%%%%%%
\begin{pro}  \label{pro:0301n}
Under the assumption \eqref{aprpioses}, if $\delta$ is sufficiently small, then
there is an energy functional $\tilde{\ml{E}}^L$ that is equivalent to $\mathcal{E}^L$, such that
$$ \frac{\mm{d}}{\mm{d}t} \tilde{\mathcal{E}}^{L}+\mathcal{D}^{L}\leq  0\mbox{ on }(0,T]. $$
\end{pro}
\begin{pf}
Noting that \eqref{inteposlai} still holds in the two-dimensional case, we have
$$ \|  (\partial_1\eta ,\mm{div}\eta)\|_{  \underline{2},0}^2\leq \epsilon
 \|(\partial_1\eta ,\mm{div}\eta)\|_{3,0}^2+C_\epsilon\|  (\partial_1\eta ,\mm{div}\eta)\|_0^2.$$
Making use of the above interpolation inequality, Lemmas \ref{ssebadiseqin}--\ref{ssebadsdiseqinsd} and \eqref{insoslai}, we see that
there is a constant $\kappa_{0,1}^L\geq 1$, such that
\begin{equation}  \label{energy1}
\begin{aligned}
 \frac{\mm{d}}{\mm{d}t}\tilde{\mathcal{E}}_1^L+  c_{1}^L \|(\partial_1  \eta, \mm{div} \eta ,\nabla   u)\|_{\underline{3},0}^2  \leq c_{2}^L( \|  (\partial_1\eta ,\mm{div}\eta)\|_{0}^2 + \|u_t\|_2^2 + \sqrt{\mathcal{E}^H} \mathcal{D}^L)\end{aligned}
 \end{equation}
 holds for any $\kappa^L_1\geq \kappa_{0,1}^L$, where
$$\begin{aligned}
&\tilde{\mathcal{E}}_1^L:= \sum_{\alpha_1+\alpha_2\leq 3}\bigg(\int \bar{\rho}J^{-1}\partial_{2}^{\alpha_1}\partial_{3}^{\alpha_2}
\eta\cdot\partial_{2}^{\alpha_1}\partial_{3}^{\alpha_2} u\mm{d}y -\kappa^L_1 E(\partial_{2}^{\alpha_1}\partial_{3}^{\alpha_2}\eta)
+\kappa^L_1\|\sqrt{\bar{\rho}J^{-1}}\partial_{2}^{\alpha_1}\partial_{3}^{\alpha_2}u\|_0^2\bigg) \\
&\qquad\quad  +\frac{1}{2}\|(\sqrt{\mu_1} \nabla \eta,\sqrt{\mu_2}\mm{div}  \eta)\|_{\underline{3},0}^2,
\end{aligned}$$
and $\kappa_{0,1}^L$ depends on $\Omega$ and the physical parameters.

Putting the estimate \eqref{energy1} and Lemma \ref{lem:dfifessim} together, we arrive at
$$ \frac{\mm{d}}{\mm{d}t}(\mathcal{H}_2(\eta)) +c_3^L\tilde{\mathcal{E}}_1^L) +\|(\partial_1\eta, \mm{div}\eta)\|_{3}^2+\|u\|_{4}^2
\leq c_4^L (\|(\partial_1\eta ,\mm{div}\eta )\|_{0}^2 +\|u_t\|_2^2+ \sqrt{\mathcal{E}^H} \mathcal{D}^L), $$
which, together with \eqref{dfifeddssimlas}, Lemma \ref{badiseqin} and the interpolation inequality, yields
 \begin{equation}   \label{Hetdsa1502} \begin{aligned}
& \frac{\mm{d}}{\mm{d}t}( \mathcal{H}_2(\eta)+c_3^L\tilde{\mathcal{E}}_1^L+c_5^L \|(\sqrt{\mu_1} \nabla_{\ml{A}}
 u_t,\sqrt{ \mu_2} \mm{div}_{\mathcal{A}}u_t)\|^2_{0}) +c_6^L(\|(\partial_1\eta, \mm{div}\eta)\|_{3}^2+\|u\|_{4}^2+\|u_{tt}\|_0^2) \\
 &\quad\leq c_7^L( \|(\partial_1\eta ,\mm{div}\eta,u)\|_{0}^2 +\sqrt{\mathcal{E}^H}\mathcal{D}^L).
 \end{aligned}  \end{equation}
So, for a given $\kappa^L_1$, we can further deduce from \eqref{Hetdsa1502} and Lemmas \ref{ssebadiseqin}--\ref{ssebadsdiseqinsd}
with $j=0$ that there is a constant $\kappa_{0,2}^L>1$, such that
 \begin{equation}  \label{emdsldsnew}
\frac{\mm{d}}{\mm{d}t} \tilde{\mathcal{E}}^{L}_2+c_7^L\tilde{\mathcal{D}}^{L}\leq  c_8^L\sqrt{\mathcal{E}^H}\mathcal{D}^L
\quad\mbox{ for any }\kappa^L_2>\kappa^L_{0,2},   \end{equation}
where
$$ \begin{aligned}  &\begin{aligned}
\tilde{\mathcal{E}}^{L}_2:= & \mathcal{H}_2(\eta)+c_3^L\tilde{\mathcal{E}}_1^L+c^{L}_5\|(\sqrt{\mu_1J} \nabla_{\ml{A}}   u_t,
   \sqrt{\mu_2J}\mm{div}_{\mathcal{A}}u_t)\|^2_{0}.
 + (\kappa_2^L)^2( \|\sqrt{\bar{\rho}J^{-1} }  u \|_{0}^2-E(\eta) )\\
 &+\kappa_2^L\Big( \int\bar{\rho}J^{-1}\eta\cdot u\mm{d}y+\frac{1}{2}\|(\sqrt{\mu_1}\nabla \eta,\sqrt{\mu_2}\mm{div}\eta)\|_0^2\Big),
 \end{aligned} \\
 & \tilde{\mathcal{D}}^{L}:=  \|(\partial_1\eta,\mm{div}\eta)\|_{3}^2+\|u\|_{4}^2 + \| u_{tt}\|^2_{0}, \end{aligned}$$
and $\kappa_{0,2}^L$ depends on $\kappa^L_1$, $\Omega$ and the physical parameters.

Obviously, to prove the proposition, it suffices to show
\begin{equation}  \label{emdsldsnewestima}
\tilde{\ml{E}}^L_2\mbox{ is equivalent to } \mathcal{E}^L\mbox{ for some }\kappa^L_1\mbox{ and }\kappa^L_2
\end{equation}
and
 \begin{equation}  \label{emdsldsne1529w}
  {\ml{D}}^L\leq c_{9}^L \tilde{\ml{D}}^L.     \end{equation}
We next verify these two facts.

Thanks to \eqref{prtislsafdsfsfds}, \eqref{Jdetemrinat2nes}, \eqref{Jdetemrinat} and \eqref{insoslai}, it is easy to see that
for sufficiently small $\delta$,
 \begin{align} &\label{f11392016}
 \|f\|_1^2\lesssim\|\sqrt{J}\nabla_{\ml{A}}f\|_{0}^2 +\|\sqrt{J}\mm{div}_{\mathcal{A}}f\|^2_{0}\lesssim\|f\|_1^2 ,\\
& \|f\|^2_0\lesssim\|\sqrt{\bar{\rho}J^{-1}} f\|_0^2\lesssim\|f\|^2_{ 0}. \nonumber
\end{align}
Moreover, $\|\eta\|_4$ is equivalent to the norm $\|\eta\|_{4,0}+\|\partial_1\eta\|_{\underline{3},0}+\|\partial^2_1\eta\|_2$.
 Thus, making use of Cauchy-Schwarz's inequality, \eqref{Hetaprofo}, Lemma \ref{lem:0601}, and \eqref{insoslai}, we conclude
\begin{align}\label{dfifessim1642}
\|\eta\|_4^2 + \|u_t\|_{1}^2 \lesssim  \tilde{\ml{E}}_2^L   \end{align}
 by choosing sufficiently large constants $\kappa_1^L$ and $\kappa^L_2$.
%%%% depending on $\Omega$ and other physical parameters in the perturbation equations.
Combining \eqref{dfifessim1642} with \eqref{dfifessim}, one arrives at
\begin{align*}   \mathcal{E}^L\lesssim \tilde{\ml{E}}_2^L+  {\ml{E}^H}\mathcal{E}^L.   \end{align*}
 In particular, $\mathcal{E}^L \lesssim \tilde{\ml{E}}_2^L $ for sufficiently small $\delta$.
 On the other hand, by Cauchy-Schwarz's inequality we observe that
$\tilde{\mathcal{E}}^{L}_2\lesssim {\mathcal{E}}^{L}$. This gives \eqref{emdsldsnewestima}.

Finally, from \eqref{dfifeddssimlas} we get $\mathcal{D}^L\lesssim\tilde{\mathcal{D}}^L+{\ml{E}^H}\ml{D}^L$,
which implies \eqref{emdsldsne1529w} for sufficiently small $\delta$.
With \eqref{emdsldsnewestima} and \eqref{emdsldsne1529w} in hands, we immediately obtain the desired conclusion from \eqref{emdsldsnew}
by defining $\tilde{\mathcal{E}}^L:=2c_{9}^L\tilde{\mathcal{E}}^L_2/c_{7}^L$ and choosing $\delta\leq c_7^L/2c_{8}^Lc_{9}^L$.
 \hfill$\Box$
\end{pf}

\section{Higher-order energy inequalities}\label{sec:03}
In this section we derive the higher-order energy inequalities for the transformed {Parker problem}.
We shall first establish the higher-order versions of Lemmas \ref{badiseqin}--\ref{lem:dfifessimell} in sequence.
%%%%%%%%%%%%%%%%%%%%%%%%%%%%%%%%%%%%%%%%%%
\begin{lem}\label{Lem:0301}
It holds that
$$   \begin{aligned}  \begin{aligned}
 & \frac{\mm{d}}{\mm{d}t}\left(\|\sqrt{\bar{\rho}} \partial_t^3 u\|_{0}^2-E (u_{tt}) \right)
 +c\| \partial_t^3  u \|^2_{1}\lesssim \sqrt{\mathcal{E}^L} \mathcal{D}^H.    \end{aligned}
 \end{aligned}$$
\end{lem}
\begin{pf}
Multiplying \eqref{01s06p}$_2$ with $j=3$ by $ J\partial_t^3u$, integrating then (by parts) the resulting equations over $\Omega$,
and using \eqref{peturbationequation}$_1$ and \eqref{AklJ=0}, we get
\begin{equation}  \label{estimdba}
\begin{aligned}
 &\frac{1}{2}\frac{\mm{d}}{\mm{d}t}\int \bar{\rho} |\partial_t^3u|^2\mm{d}y
    +\int {J}(\mu_1 |\nabla_{\mathcal{A}}\partial_t^3u|^2+\mu_2 |\mm{div}_{\mathcal{A}}\partial_t^3u|^2)\mm{d}y \\
 & =\int\left(g\mm{div} (\bar{\rho}u_{tt}) \partial_t^3 u_3- P'(\bar{\rho})\mm{div}(\bar{\rho} u_{tt} )  \mm{div}\partial_t^3u\right)
\mm{d}y + \int \partial_t^3  \mathcal{L}_M\cdot \partial_t^3u\mm{d}y \\
  &\quad + \int J\partial_t^3\mathcal{N}\cdot \partial_t^3u\mm{d}y+  \int J N^{t,3}_u  \cdot \partial_t^3u\mm{d}y\\
 &\quad +  \int  (J-1)(g\mm{div}( \bar{\rho}u_{tt}) e_3+ \nabla (P'(\bar{\rho})\mm{div}(\bar{\rho} u_{tt} ) )
 +  \partial_t^3\mathcal{L}_M)\cdot \partial_t^3u\mm{d}y =:\sum_{k=1}^5I^{H}_k,
 \end{aligned}
\end{equation}
where the first two integrals on the right hand side can be written as
\begin{equation*}\label{20160129}
\begin{aligned}
I_1^H = &\frac{1}{2}\frac{\mm{d}}{\mm{d}t}\left(\int g\bar{\rho}'|\partial_t^2 u_3|^2 \mm{d}y
+2\int g\bar{\rho}\partial_t^2 u_3\mm{div}u_{tt}\mm{d}y - \int {P'(\bar{\rho})\bar{\rho}} |\mm{div}u_{tt}|^2\mm{d}y\right)\\
 & - \int \left(  \bar{P}'  + g\bar{\rho}\right)\partial_t^2 u_3\mm{div}\partial_t^3 u\mm{d}y \end{aligned}
\end{equation*}
and
\begin{equation*}\label{20160130}\begin{aligned}I_2^H= -\frac{\lambda}{2}\frac{\mm{d}}{\mm{d}t}\int m^2
(|\partial_1\partial_{t}^2 u_\mm{v}|^2 + |\mm{div}_{\mm{v}}\partial_t^2u_{\mm{v}} |^2
  )\mm{d}y -\lambda\int m m '\partial_t^2 u_3 \mm{div}\partial_t^3u\mm{d}y.
\end{aligned}\end{equation*}
Substituting the above two identities into \eqref{estimdba} and using \eqref{comsteady}, we arrive at
\begin{equation}\label{201603181521MH}
\begin{aligned}
 & \frac{1}{2}\frac{\mm{d}}{\mm{d}t}\left(\|\sqrt{\bar{\rho} } \partial_t^3 u\|_{0}^2-E ( u_{tt}) \right)
 +\|(\sqrt{\mu_1J}\nabla_{\mathcal{A}} \partial_t^3 u,\sqrt{\mu_2J}\mm{div}_{\mathcal{A}}\partial_t^3 u )\|^2_{0}=\sum_{k=3}^5I^{H}_k,
 \end{aligned}    \end{equation}
where the terms on the right hand side can be bounded as follows, using
\eqref{badiseqin3nes1}, \eqref{Nu1743MHtdfsNsg}, \eqref{Jdetemrinat2nes} and  \eqref{Jdetemrinat}.
 \begin{align*}
 &I^{H}_3\lesssim  \|\partial_t^3  \mathcal{N}\|_0\|\partial_t^3u\|_0\lesssim  \sqrt{\mathcal{E}^L} \mathcal{D}^H,\\
 &\begin{aligned}
 I^{H}_4\lesssim\|N^{t,3}_u\|_0  \|\partial_t^3u\|_0  \lesssim   \sqrt{\mathcal{E}^L} \mathcal{D}^H, \end{aligned} \\
& I^{H}_5\lesssim \|J-1\|_2\|u_{tt}\|_{2} \|\partial_t^3u\|_{0} \lesssim \sqrt{\mathcal{E}^L} \mathcal{D}^H.
\end{align*}
Consequently, inserting the above three inequalities into \eqref{201603181521MH} and using \eqref{f11392016}, we
obtain the lemma.  \hfill$\Box$
\end{pf}
%%%%%%%%%%%%%%%%%%%%%%%%%%%%%%%%%%%%
\begin{lem}   \label{Lem:0302}
It holds that$$\begin{aligned}
&\frac{\mm{d}}{\mm{d}t}\left( \int \bar{\rho}J^{-1} \partial_\mm{v}^j \eta \cdot  \partial_\mm{v}^j u\mm{d}y
+ \frac{1}{2} \|\partial_\mm{v}^j(\sqrt{\mu_1}
\nabla \eta,\sqrt{\mu_2}\mm{div}   \eta)\|_{0}^2\right) +  c\|\partial_\mm{v}^j (  \partial_1\eta ,\mm{div} \eta) \|_0^2 \\
&\quad \lesssim  \|  (\partial_1\eta ,\mm{div}\eta)\|_{\underline{j-1},0}^2+\|\partial_\mm{v}^j   u\|^2_{0}+\|u_t\|_{j-1}^2+ \sqrt{\mathcal{E}^L}( \|\eta\|_7^2+\mathcal{D}^H),\;\quad j=4,\ldots,6. \end{aligned}$$
\end{lem}
\begin{pf} We only prove the case $j=6$, while the rest cases can be shown in the same way.
 Following the process of deriving \eqref{estimforhoedsds1st}, we obtain
 \begin{equation}  \label{estimforhoedsds1st1307}
\begin{aligned}
&\frac{\mm{d}}{\mm{d}t}\left(\int \bar{\rho} J^{-1}\partial_\mm{v}^6 \eta \cdot \partial_\mm{v}^6 u \mm{d}y
+\frac{1}{2}\|\partial_\mm{v}^6 (\sqrt{\mu_1}\nabla\eta ,\sqrt{\mu_2}\mm{div}  \eta)\|_{0}^2\right)  -  E (\partial_\mm{v}^6 \eta)\\
&\lesssim c\left(\| \partial_{\mm{v}}^6 u\|^2_{0} + \|\partial_\mm{v}^6 (\partial_1 \eta,\mm{div}\eta)\|_{0} \|
 (\partial_1\eta ,\mm{div}\eta)\|_{\underline{5},0}\right) +\sum_{k=1}^4J^{H}_{k},
\end{aligned}  \end{equation}
where
$$ \begin{aligned}
\sum_{k=1}^4 J^{H}_{k}:= &\int\partial_\mm{v}^6\mathcal{N}\cdot\partial_\mm{v}^6\eta\mm{d}y +\int\partial_\mm{v}^6 N_u\cdot\partial_\mm{v}^6\eta\mm{d}y\\
& -\sum_{ k=0}^5\int (\partial_\mm{v}^{6}(\bar{\rho}J^{-1} u_t )- \bar{\rho}J^{-1}\partial_\mm{v}^6 u_t )\cdot\partial_\mm{v}^6 \eta \mm{d}y
 +\int \bar{\rho} J^{-1}_t\partial_\mm{v}^6\eta \cdot \partial_\mm{v}^6 u \mm{d}y ,\end{aligned}$$
and $J^H_1,\cdots, J^H_4$ can be estimated as follows, using \eqref{Nu1743MH}, \eqref{N21556201630266},
and \eqref{Jtestismtsinver}, \eqref{Jdetemrinat2nes} and the smallness condition.
\begin{align}  \label{201603311524}
J^H_1 &\lesssim \|\partial_\mm{v}^5\mathcal{N} \|_0\|\partial_\mm{v}^6\eta \|_1　\lesssim \sqrt{\mathcal{E}^L}\|\eta\|_7^2 ,\\
J^H_2 &\lesssim \|\partial_\mm{v}^5N_u \|_0\|\partial_\mm{v}^6\eta \|_1 \lesssim\sqrt{\mathcal{E}^L}\|(\eta,u)\|_7^2 ,\\
J^H_3 &\lesssim (\|J^{-1}-1\|_6+1)\|u_t\|_{5}  \|\partial_{\mm{v}}^6\eta\|_0 \lesssim\|u_t\|_{5}\|\partial_{\mm{v}}^6\eta\|_0 ,\\
\label{201603311524new} J^H_4 &\lesssim\|J^{-1}_t\|_2\|\partial_\mm{v}^6 u\|_0\|\partial_\mm{v}^6\eta\|_0\lesssim\sqrt{\mathcal{E}^L}\mathcal{D}^H.
\end{align}
Consequently, if we insert the above four inequalities into \eqref{estimforhoedsds1st1307} and use \eqref{insoslai}, we obtain Lemma \ref{Lem:0302}.  \hfill$\Box$
\end{pf}
%%%%%%%%%%%%%%%%%%%%%%%%%%%%%%%%
\begin{lem}   \label{Lem:0303}
It holds that$$\begin{aligned}
& \frac{\mm{d}}{\mm{d}t}\left(\|\sqrt{ \bar{\rho}J^{-1} }\partial_\mm{v}^j u\|^2_0-E (\partial_{\mm{v}}^j\eta)\right)
+ c\|\partial_\mm{v}^j(\nabla u,\mm{div} u)\|_{0}^2\\
 & \quad \lesssim  \|  ( \partial_1\eta ,\mm{div}\eta)\|_{\underline{j-1},0}^2+\|u_t\|_{j-1}^2+
 \sqrt{\mathcal{E}^L}(\|\eta\|_7^2+\mathcal{D}^H ),\;\quad j=4,\ldots,6. \end{aligned}$$
\end{lem}
\begin{pf}
We only prove the case $j=6$, since the rest cases are similar to deal with. Following the derivation of \eqref{estimforho4e1nwe}, one deduces that
\begin{equation}    \label{estimforho4e1nwe2256}
 \frac{1}{2}\frac{\mm{d}}{\mm{d}t}\left(\|\sqrt{ \bar{\rho}J^{-1} }\partial_\mm{v}^6 u\|^2_0
 - E (\partial_{\mm{v}}^6\eta)\right)+\|\partial_\mm{v}^6 (\sqrt{\mu_1}\nabla u,\sqrt{\mu_2}\mm{div}u)\|_{0}^2
 \leq c \|(\partial_1\eta ,\mm{div}\eta)\|_{\underline{5},0}^2+\sum_{k=1}^4 K^{H}_{k},
\end{equation}
where
\begin{equation*}   \begin{aligned}
\sum_{k=1}^4 K^{H}_{k}:=&\int \partial_\mm{v}^6 \mathcal{N}\cdot\partial_\mm{v}^6 u\mm{d}y +\int\partial_\mm{v}^6 N_u\cdot\partial_\mm{v}^6 u\mm{d}y\\
&\quad +\frac{1}{2}\int \bar{\rho} J^{-1}_t|\partial_\mm{v}^6 u|^2 \mm{d}y -\sum_{ k=0}^5\int
\partial_\mm{v}^{6-k}(\bar{\rho}J^{-1})\partial_\mm{v}^k u_t\cdot \partial_\mm{v}^6 u\mm{d}y .
\end{aligned}
\end{equation*}
On the other hand, similarly to \eqref{201603311524}--\eqref{201603311524new}, we have
\begin{align*}  \sum_{k=1}^4 K^{H}_{k} \lesssim  \|u_t\|_{5}
  \|\partial_{\mm{v}}^6u\|_0 +\sqrt{\mathcal{E}^L} (\|\eta\|_7^2+\mathcal{D}^H). \end{align*}
Inserting the above estimate into \eqref{estimforho4e1nwe2256}, we can use \eqref{insoslai} and Cauchy-Schwarz's
inequality to obtain the desired estimate.   \hfill$\Box$
\end{pf}
%%%%%%%%%%%%%%%%%%%%%%%%%%%%%%%%%%%%%%%%
\begin{lem}\label{lemm:04041447}
It  holds that
$$ \frac{\mm{d}}{\mm{d}t}\mathcal{H}_5(\eta) +\|(\partial_1\eta,\mm{div}\eta)\|_{6}^2+\|u\|_{7}^2
 \lesssim \|(\partial_1 \eta,\mm{div}\eta, \nabla u)\|_{\underline{6},0}^2
 + \| u_t\|_{5}^2+{\mathcal{E}^L}(\|\eta\|_7^2+\mathcal{D}^H), $$
 where the energy functional $\mathcal{H}_5(\eta)$ satisfies
$\|\partial_1^2\eta\|_{5}^2- h_5\|\partial_1\eta_{\mm{v}}\|_{\underline{6},0}^2\lesssim {\mathcal{H}}_{5}(\eta)$
and $h_5$ is a positive constant depending on $\Omega$ and the physical parameters in the perturbation equations.
\end{lem}
\begin{pf} The lemma easily follows from \eqref{Nu1743MH}, \eqref{N21556201630266} and \eqref{omdm12dfs2}.
\hfill$\Box$
\end{pf}
%%%%%%%%%%%%%%%%%%%%%%%%%%%%%%%%%%%%%%%%%%%%%%%%%%%%%%%%
\begin{lem}\label{lem:dfifessim1406} We have
\begin{align} \label{omdm122nsdfsf}
& \sum_{j=0}^2\|\partial_t^j u\|_{6-2j}^2\lesssim\|\eta\|_{6}^2 +\|(u,u_t,\partial_t^3 u)\|_0^2 +{\mathcal{E}^L}\mathcal{E}^H,\\
& \label{highestdids} \sum_{j=1}^2\|\partial_t^j u\|_{7-2j}^2\lesssim\|u\|_5^2 +\|(u_t,\partial_t^3 u)\|_1^2+\mathcal{E}^L\mathcal{D}^H.
 \end{align}
 \end{lem}
\begin{pf} (1) Taking $(i,j)=(4,0)$ in the Stokes estimate \eqref{02181} and using \eqref{Jdetemrinat2nes}, we find that
 \begin{equation*}  \label{omdm122n} \begin{aligned}
 &  \|  u\|_{6}^2 \lesssim \| \eta\|_{6}^2+\|u_t\|_{4}^2  +S_{0,4}.
 \end{aligned}\end{equation*}
On the other hand, it follows from \eqref{Nu1743MHnew} and \eqref{N21556201630266} that
$ S_{0,4}=\|(\mathcal{N},N_u)\|_{4}^2 \lesssim\mathcal{E}^L \mathcal{E}^H$. Therefore,
\begin{equation}  \label{omdjfm122n}\begin{aligned}
 &  \|  u\|_{6}^2 \lesssim \| \eta\|_{6}^2+\|u_t\|_{4}^2    +\mathcal{E}^L \mathcal{E}^H .
 \end{aligned}\end{equation}

If one utilizes \eqref{Jdetemrinat2nes}, \eqref{peturbationequation}$_1$ and the recursion formula \eqref{02181} with $i=4$ for
$j=1$ and $2$, one obtains
\begin{equation*}
\label{}\begin{aligned}
     \sum_{j=1}^2 \| \partial_t^j u\|_{6-2j}^2
  \lesssim  \sum_{j=1}^2  \|\partial_t^{j-1}  u\|_{6-2j}^2+ \|\partial_t^{3}u\|_{0}^2+\sum_{j=1}^2 S_{j,4}.
 \end{aligned}\end{equation*}
On the other hand, from \eqref{Nu1743MHt}, \eqref{N215562016302661Ntt}, \eqref{N215562016302661}, \eqref{Jtestismtsinver} and \eqref{fgestims} we get
$$ \begin{aligned} \sum_{j=1}^2 S_{j,4} =& \sum_{j=1}^2\left( \|\partial_t^j ( \mathcal{N}, N_u)\|_{4-2j}^2
+\left\|\bar{\rho}\sum_{ k =0}^{j-1}\partial_{t}^{j-k} J^{-1}\partial_t^{k+1}u \right\|_{4-2j}^2\right)
  \lesssim  \mathcal{E}^L \mathcal{E}^H.   \end{aligned}$$
Thus,
\begin{equation}   \label{sdfsdfomdm122n}\begin{aligned}
 \sum_{j=1}^2 \| \partial_t^j u\|_{6-2j}^2  \lesssim \sum_{j=1}^2  \|\partial_t^{j-1}  u\|_{6-2j}^2
  +\|\partial_t^{3}u\|_{0}^2+\mathcal{E}^L \mathcal{E}^H, \end{aligned}\end{equation}
which, together with \eqref{omdjfm122n}, yields
\begin{equation}  \label{jfomdm122n} \begin{aligned}
 \sum_{j=0}^2 \| \partial_t^j u\|_{6-2j}^2
  \lesssim \|\eta\|_6^2+\sum_{j=1}^2 \|\partial_t^{j-1}  u\|_{6-2j}^2 + \|\partial_t^3u\|_0^2+\mathcal{E}^L \mathcal{E}^H.
 \end{aligned}\end{equation}
Finally, employing the interpolation inequality, we get \eqref{omdm122nsdfsf} from \eqref{jfomdm122n}.

(2) We now turn to the proof of \eqref{highestdids} for higher-order dissipation estimates.
Making use of \eqref{Jdetemrinat2nes}, \eqref{peturbationequation}$_1$ and the recursion formula \eqref{02181} with $i=5$
for $j=1$ to $2$, we obtain
\begin{equation}    \label{dfeossdsfsd1} \begin{aligned}
 \sum_{j=1}^2  \| \partial_t^j u\|_{7-2j}^2
  \lesssim \sum_{j=1}^2\| \partial_t^{j-1}u\|_{7-2j}^2+\|\partial_t^3 u\|_{1}^2 +\sum_{j=1}^2 S_{j,5}.
 \end{aligned}\end{equation}
On the other hand, it follows from \eqref{Nu1743MHtdfssd}, \eqref{Nu1743MHtdfsNsg}, \eqref{Jtestismtsinver}, \eqref{fgestims} and the
interpolation inequality that
\begin{equation}  \label{201604031311}
\begin{aligned}\sum_{j=1}^2 S_{j,5}=& \sum_{j=1}^2\left( \|\partial_t^j (\mathcal{N},N_u)\|_{5-2j}^2+\left\| \bar{\rho }
\sum_{ k =0}^{j-1}\partial_{t}^{j-k} J^{-1}\partial_t^{k+1}u  \right\|_{5-2j}^2\right)\\
  \lesssim &\mathcal{E}^L\mathcal{D}^H+\|u\|_4\|u_t\|_2 \lesssim \mathcal{E}^L\mathcal{D}^H.
  \end{aligned}   \end{equation}
Finally, exploiting the interpolation inequality again, we get \eqref{highestdids} from \eqref{dfeossdsfsd1}
and \eqref{201604031311}.      \hfill$\Box$
\end{pf}

Now, we are ready to build the higher-order energy inequalities. In what follows the letters $c_i^H$, $i=1,\ldots,8$, will
denote generic constants which may depend on $\Omega$ and the physical parameters in the perturbation equations.
%%%%%%%%%%%%%%%%%%%%%%%%%%%%%%%%%%%%%%%%%%%%%%%%%%%%%
\begin{pro} \label{pro:0301}
Under the assumption \eqref{aprpioses}, if $\delta$ is sufficiently small, then
there is a norm $\tilde{\ml{E}}^H$, which is equivalent to $\mathcal{E}^H$, such that
\begin{equation*}
\frac{\mm{d}}{\mm{d}t} \tilde{\mathcal{E}}^{H}+\mathcal{D}^{H}\lesssim \sqrt{\mathcal{E}^L}\|\eta\|_7^2.
\end{equation*}    \end{pro}
\begin{pf}
Similarly to \eqref{energy1}, we employ Lemmas \ref{Lem:0302}--\ref{Lem:0303}, \eqref{insoslai} and the interpolation inequality
to see that there is a constant $\kappa_{0,1}^H\geq 1$, such that
\begin{equation}   \label{energy1hier}
\begin{aligned}   \frac{\mm{d}}{\mm{d}t}\tilde{\mathcal{E}}_1^H +c_1^H \|(\partial_1\eta , \mm{div}\eta ,\nabla u)\|_{\underline{6},0}^2
  \leq & c_{2}^H\Big( \|(\partial_1\eta ,\mm{div}\eta)\|_0^2 +  \|u_t\|_5^2 +\sqrt{\mathcal{E}^L}(\|\eta\|_7^2+\mathcal{D}^H)\Big)
  \end{aligned}    \end{equation}
for any $\kappa^H_1 \geq \kappa_{0,1}^H$, where
$$\begin{aligned}
&\tilde{\mathcal{E}}_1^H:= \sum_{4\leq \alpha_1+\alpha_2\leq 6}\bigg(\int \bar{\rho}J^{-1}\partial_{2}^{\alpha_1}\partial_{3}^{\alpha_2}
\eta \cdot  \partial_{2}^{\alpha_1}\partial_{3}^{\alpha_2} u\mm{d}y-\kappa^H_1 E ( \partial_{2}^{\alpha_1}\partial_{3}^{\alpha_2}\eta)\\
&\qquad\qquad \qquad \qquad  +\kappa^H_1\|\sqrt{ \bar{\rho}J^{-1} }  \partial_{2}^{\alpha_1}\partial_{3}^{\alpha_2} u\|_0
+\frac{1}{2}\|\partial_{2}^{\alpha_1}\partial_{3}^{\alpha_2}( \sqrt{\mu_1} \nabla \eta,\sqrt{\mu_2}\mm{div}  \eta)\|_{0}^2\bigg) ,
\end{aligned}$$and
$\kappa_{0,1}^H$ depends on $\Omega$ and the physical parameters.

On the other hand, analogously to \eqref{Hetdsa1502}, we use \eqref{highestdids} and the interpolation inequality to deduce
from \eqref{energy1hier}, Lemmas \ref{lemm:04041447} and \ref{Lem:0301} that
 \begin{equation}    \label{Hetdsa150204}
\begin{aligned}
& \frac{\mm{d}}{\mm{d}t}( \mathcal{H}_5(\eta)+c_3^H\tilde{\mathcal{E}}_1^H+c_4^H(\|\sqrt{\bar{\rho}} \partial_t^3 u\|_{0}^2-E (u_{tt}) ))
  + c_5^H(\|(\partial_1\eta, \mm{div}\eta)\|_{6}^2+\|u\|_{7}^2+\|\partial_t^3 u\|_1^2) \\
 &\leq c_{6}^H(\|(\partial_1\eta ,\mm{div}\eta)\|_0^2 +  \|u\|_0^2+\|u_t\|_1^2+ \sqrt{\mathcal{E}^L}(\|\eta\|_7^2+ \mathcal{D}^H),
 \end{aligned}
\end{equation}
which, together with Proposition \ref{pro:0301n}, implies that there is a constant $\kappa_{0,2}^H\geq 1$, such that
 \begin{equation}  \label{emdsldsnewnew}
\frac{\mm{d}}{\mm{d}t}\tilde{\mathcal{E}}^{H}_2+c_{7}^H\tilde{\mathcal{D}}^{H}\leq  c_{6}^H\sqrt{\mathcal{E}^L}(\|\eta\|_7^2+ \mathcal{D}^H)
\quad\mbox{ for any }\kappa^H_2\geq \kappa^H_{0,2},  \end{equation}
where
$$ \begin{aligned}
\tilde{\mathcal{E}}^{H}_2:= \mathcal{H}_5(\eta)+c_3^H\tilde{\mathcal{E}}_1^H
+c_4^H(\|\sqrt{\bar{\rho}} \partial_t^3 u\|_{0}^2-E (u_{tt}) )+\kappa_2^H\tilde{\mathcal{E}}^{L},   \end{aligned}   $$
 $$\begin{aligned}
\tilde{\mathcal{D}}^{H}:= & \|(\partial_1\eta,\mm{div}\eta)\|_{6}^2+\|u\|_{7}^2 + \|( u_t, \partial_t^3u)\|^2_{1}. \end{aligned}$$

Similarly to \eqref{emdsldsnewestima} and \eqref{emdsldsne1529w}, we can utilize Lemma \ref{lem:dfifessim1406}, \eqref{Jdetemrinat22}
and the fact that
$$\|\eta\|_7\mbox{ is equivalent the norm } \|\eta\|_{7,0}+\|\partial_1\eta\|_{\underline{6},0}+\|\partial^2_1\eta\|_5$$
to deduce that $\tilde{\ml{E}}^H_2$ is equivalent to $\mathcal{E}^H$ and ${\ml{D}}^H\leq c_8^H\tilde{\ml{D}}^H$
 by choosing sufficiently large constants $\kappa_1^H$ and $\kappa^H_2$ which depend on $\Omega$ and the physical parameters.
Consequently, one obtains the desired conclusion from \eqref{emdsldsnewnew} by defining
$\tilde{\mathcal{E}}^H:=2c_{8}^H\tilde{\mathcal{E}}^H_2/c_{7}^H$ and choosing $\delta\leq c_7^H/2c_{6}^Hc_8^H$.
 \hfill$\Box$   \end{pf}

Finally, the following lemma will be needed in estimating the initial energy $\mathcal{E}^H(0)$ in the next section.
\begin{lem}\label{esdoval} Under the assumption \eqref{aprpioses}, if $\delta$  is sufficiently small, then
$\mathcal{E}^H\lesssim  \mathcal{E}$, where $\mathcal{E}:= \|\eta\|_{7}^2+ \|u\|_6^2 $.
\end{lem}
\begin{pf}In view of \eqref{omdm122nsdfsf}, we have $\mathcal{E}^H\lesssim\|\eta\|_{7}^2 +\|(u,u_t,\partial_t^3u)\|_0^2+\mathcal{E}^L\mathcal{E}^H$,
which gives
\begin{equation}   \label{EHlesimms}
\mathcal{E}^H\lesssim\|\eta\|_{7}^2+\|(u,u_t,\partial_t^3u)\|_0^2 \qquad\mbox{for sufficiently small }\; \delta.  \end{equation}
Next, we show that the $L^2$-norm of $u_t$ and $\partial_t^3 u$ can be controlled by $\sqrt{\mathcal{E}}$.

Multiplying \eqref{01s06p}$_2$ with $j=2$ by $\partial_t^3 u$ in $L^2$,
we arrive at
$$\begin{aligned}
 \|\sqrt{\bar{\rho}J^{-1}}\partial_t^{3} u\|_0^2 = & \int \Big( \mu_1 \Delta_{\ml{A}}u_{tt} +\mu_2
 \nabla_{\mathcal{A}}\mm{div}_{\mathcal{A}}u_{tt} +\nabla (P'(\bar{\rho}) \mm{div}(\bar{\rho}u_t)) \\
& \quad + g\mm{div}( \bar{\rho}u_t) e_3 + \partial_t^2\mathcal{L}_M+ \mathcal{N}_{tt} + N^{t,2}_u\Big)\cdot\partial_t^3 u\mm{d}y.
\end{aligned}$$
 Employing \eqref{prtislsafdsfsfds}, \eqref{Jdetemrinat22}, \eqref{fgestims} and \eqref{0102}, we see that
$$ \| \partial_t^{3} u\|_0^2 \lesssim \|(u_t,u_{tt})\|_2^2 +\| ( \mathcal{N}_{tt}, N^{t,2}_u )\|^2_0, $$
which, together with \eqref{badiseqin31new} and \eqref{N215562016302661Ntt}, yields
\begin{equation}  \label{uttestimdfd}
\|\partial_t^3 u\|_0^2 \lesssim \|u_{t}\|_2^2+ \| u_{tt}\|_2^2 +\mathcal{E}.
\end{equation}
Here $\|u_{tt}\|_2$ can be bounded as follows, using \eqref{prtislsafdsfsfds}, \eqref{Jdetemrinatnesw}, \eqref{fgestims}
and \eqref{01s06p}$_2$ with $j=1$.
$$ \begin{aligned}  \| u_{tt}\|_2^2=& \|\bar{\rho}^{-1}J( \mu_1 \Delta_{\ml{A}}u_t+\mu_2
 \nabla_{\mathcal{A}}\mm{div}_{\mathcal{A}}u_t +\nabla (P'(\bar{\rho}) \mm{div}(\bar{\rho}u))
 \\
 & +g\mm{div}( \bar{\rho}u) e_3 + \partial_t\mathcal{L}_M+\mathcal{N}_t +  N^{t,1}_u )\|_2^2\\
\lesssim &\|u_t\|_4^2+\mathcal{E}+ \|\mathcal{N}_t + N^{t,1}_u \|_2^2,  \end{aligned}$$
which, together with \eqref{badiseqin31ass} and \eqref{N215562016302661}, results in
\begin{equation}  \label{uteswtesi}
\| u_{tt}\|_2^2\lesssim \| u_t\|_4^2 +\mathcal{E}.
\end{equation}
Similarly, we can show that
\begin{equation}
\label{uteswtesi1549} \| u_{t}\|_4^2= \|\bar{\rho}^{-1}J( \mu_1 \Delta_{\ml{A}}u+\mu_2
 \nabla_{\mathcal{A}}\mm{div}_{\mathcal{A}}u +\nabla (P'(\bar{\rho}) \mm{div}(\bar{\rho}\eta))
  +g\mm{div}( \bar{\rho}\eta) e_3 +  \mathcal{L}_M+\mathcal{N}) \|_4^2 \lesssim \mathcal{E}.  \end{equation}
As a result, we infer from \eqref{uttestimdfd}--\eqref{uteswtesi1549} that
\begin{equation}
\label{udseilaset}
\|\partial_t^3 u\|_0^2
\lesssim   \mathcal{E}.
\end{equation}

Now, substituting \eqref{udseilaset} and \eqref{uteswtesi1549} into \eqref{EHlesimms}, we obtain the lemma.  \hfill $\Box$
\end{pf}

\section{Proof of Theorems \ref{thm3} and \ref{thm3orig}} \label{sec:02}
   This section is devoted to the proof of Theorems \ref{thm3} and \ref{thm3orig}. Theorem \ref{thm3} can be shown by establishing the
       stability estimate \eqref{1.19} and the local well-posedness of the transformed {Parker problem}.
       With Theorem \ref{thm3} in hand, we can easily obtain Theorem \ref{thm3orig} by transforming Lagrangian coordinates to
       Eulerian coordinates.

\subsection{Stability estimate}
In this subsection we show the stability estimate \eqref{1.19} under the assumption
\begin{equation}     \label{201604032029}
   \mathcal{G}_3(T):= {\sup_{0\leq \tau\leq T}\mathcal{E}^H(\tau)}+\mathcal{G}_2(T)\leq \delta^2,
 \end{equation}which is stronger than \eqref{aprpioses}.

 Proposition \ref{pro:0301} implies
$$
\mathcal{G}_1(t)\lesssim \mathcal{E}^H(0)+
\int_0^t\sqrt{\mathcal{E}^L(\tau)}\|\eta(\tau)\|_{7}^2\mm{d}\tau.
$$
Then we can use \eqref{201604032029} to give
$$ \mathcal{G}_1(t) \lesssim\mathcal{E}^H(0)+\int_0^t {\delta}(1+\tau)^{-3/2}\|\eta(\tau)\|_{7}^2\mm{d}\tau
\lesssim \mathcal{E}^H(0)+  {\delta} \mathcal{G}_1(t) , $$
which implies that
\begin{equation}  \label{G3testim}
\mathcal{G}_1(t) \lesssim  \mathcal{E}^H(0).
\end{equation}

We now show the decay estimate of $\mathcal{G}_2(t)$. Note that $\mathcal{E}^L$ can be controlled by $\mathcal{D}^L$
except for $\|\eta\|_{4}$. To control $\|\eta\|_{4}$, we use the interpolation inequality \eqref{inteposlai} to get
$$
\|\eta\|_{4}\lesssim \|\eta\|_{3}^{\frac{3}{4}}\|\eta\|_{7}^{\frac{1}{4}}.$$
On the other hand,  it follows from \eqref{G3testim} and \eqref{emdsldsnewestima} that
$${\mathcal{E}}^{L}+\|\eta\|_{7}^2\lesssim \tilde{\mathcal{E}}^{L}+\|\eta\|_{7}^2\lesssim\mathcal{E}^H(0), $$
whence,
$$\tilde{\mathcal{E}}^L\lesssim{\mathcal{E}}^L\lesssim (\mathcal{D}^L)^{\frac{3}{4}}({\mathcal{E}}^{L}+\|\eta\|_{7}^2)^{\frac{1}{4}}
  \lesssim (\mathcal{D}^L)^{\frac{3}{4}} \mathcal{E}^H(0)^{\frac{1}{4}}. $$
Inserting the above estimate into the lower-order energy inequality in Proposition \ref{pro:0301n}, we obtain
$$\frac{\mm{d}}{\mm{d}t} \tilde{\mathcal{E}}^{L}+
 \frac{ (\tilde{\mathcal{E}}^{L})^{\frac{4}{3}}}{\mathcal{I}_0 ^{1/{3}}}\lesssim 0,$$
which implies
$${\mathcal{E}}^{L}\lesssim \tilde{\mathcal{E}}^{L}\lesssim \frac{\mathcal{I}_0}{\left((\mathcal{I}_0/\mathcal{E}^{L}(0))^{1/3}+  t/3\right)^3}
\lesssim  \frac{   \mathcal{E}^H(0)}{ 1 + t^3},$$
where $\mathcal{I}_0:=c \mathcal{E}^H(0)$ for some positive constant $c$. Hence,
\begin{equation}
\label{etasef12}\mathcal{G}_2(t)\lesssim   \mathcal{E}^H(0). \end{equation}

Now we add \eqref{G3testim} to \eqref{etasef12} to conclude $\mathcal{G}(t):= \mathcal{G}_1(t) +\mathcal{G}_2(t)\lesssim\mathcal{E}^H(0)$.
On the other hand, thanks to Lemma \ref{esdoval}, $\mathcal{E}^H(0)\lesssim \|\eta_0\|_7^2+\|u_0\|_6^2$.
Therefore, $\mathcal{G}(t)\lesssim\|\eta_0\|_7^2+\|u_0\|_6^2$. Consequently, the stability estimate can be summarized as follows.
\begin{pro}  \label{125pro:0401}
Let $(\eta,u)$ be a solution of the transformed {Parker problem} \eqref{peturbationequation}, \eqref{defineditcon}.
Then there is a sufficiently small $\delta_1$, such that $(\eta,u)$ enjoys the following stability estimate:
\begin{equation}\label{decadsyeste}
\mathcal{G}(T_1)\leq C_1( \|\eta_0\|_7^2+\|u_0\|_6^2 ),
\end{equation}
provided that $ \sqrt{\mathcal{G}_3(T_1)}\leq \delta_1 $ for some $T_1>0$, where
$C_1\geq 1$ denotes a constant depending on $\Omega$ and the physical parameters in the perturbation equations.
\end{pro}

\subsection{Local well-posedness}\label{sec:0402}
Now we introduce the local existence of a small classical solution to the transformed {Parker problem}.
\begin{pro} \label{pro:0401n}
There exists a sufficiently small $\delta_2$, such that for any given initial data $(\eta_0,u_0)\in H^7\times H^6$ satisfying
\begin{equation*}
\sqrt{\|\eta_0\|_7^2+\|{u}_0\|_6^2}\leq \delta_2
\end{equation*}
and compatibility conditions (i.e., $\eta(\cdot ,0)|_{\partial\Omega}=0$, $\partial_t^j u(\cdot ,0)|_{\partial\Omega}=0$ for $j=0,1,2$),
there are a $T_2:=T_2(\delta_2)>0$, depending on $\delta_2$, $\Omega$ and the known physical parameters, and a unique classical solution
$(\eta, u)\in C^0([0,T_2]$, $H^7\times H^6 )$ to the transformed {Parker problem} \eqref{peturbationequation}, \eqref{defineditcon}. Moreover,
$\partial_t^iu\in C^0([0,T_2],H^{6-2i})$ for $1\leq i\leq 3$, $\mathcal{E}^H(0)\lesssim \|\eta_0\|_7^2+\|u_0\|_6^2$, and
$$ \sup_{0\leq \tau\leq T_2} \mathcal{E}^H(\tau)+\int_0^{T_2} \mathcal{D}^H(\tau) \mm{d}\tau<\infty, $$
and $\mathcal{G}(t)$ is a continuous function on $[0,T_2]$.
\end{pro}
\begin{rem}
In view of the definition of $\mathcal{G}(t)$, we have
 \begin{equation}  \label{Gestsim}
\mathcal{G}(0)\leq  2\mathcal{E}^H(0) \leq  C_2(\|\eta_0\|_7^2+\|{u}_0\|_6^2),
 \end{equation}
 where $C_2\geq 1$ denotes a constant depending on $\Omega$ and the physical parameters.
\end{rem}
\begin{pf} The transformed {Parker problem} is very similar to the surface wave problem (1.4) in \cite{GYTILW1}.
Moreover, our problem indeed is simpler than the the surface wave problem due to the non-slip boundary condition $u|_{\partial\Omega}=0$.
Using a standard iteration method as in \cite{GYTILW1}, we can easily prove Proposition \ref{pro:0401n}, and hence
we omit the proof here. In addition,
the continuity, such as $(\eta, u)\in C^0([0,T_2], H^7\times H^6)$, $\mathcal{G}(t)$ and so on, can be verified
by using the regularity of $(\eta,u)$, the transformed perturbation equations and a standard regularized method.
\hfill $\Box$
\end{pf}

\subsection{Proof of Theorem \ref{thm3}}
    Now we are in a position to show Theorem \ref{thm3}. Let the initial data $(\eta_0,u_0) \in H^{7}\times H^6$ satisfy
 the assumptions  in Theorem \ref{thm3}, where $\delta$ in Theorem \ref{thm3} further satisfies
\begin{equation} \label{detlaci}
\delta\leq \frac{\min\{\delta_1, \delta_2\}  }{2 C_1C_2}, \end{equation}
and $\delta_1$, $\delta_2$ are the same as in Propositions \ref{125pro:0401} and \ref{pro:0401n} respectively.
Then, by Proposition \ref{pro:0401n}, there exists a solution $(\eta,u)$, defined on $\Omega\times (0,T^{\max})$,
to the transformed {Parker problem} where $T^{\max}$ denotes the maximal existence time,
such that $\mathcal{G}(t)$ is continuous for any $t\in [0,T^{\max})$. Obviously, $T^{\max}\geq T_2(\delta_2)$.
 We denote
$$T^*=\sup\{t\in (0,T^{\max})~|~ \sqrt{\mathcal{G}(t)}\leq  \min\{\delta_1,\delta_2\}\},$$
then $0<T^*\leq T^{\max}$  by \eqref{Gestsim} and the continuity of $\mathcal{G}(t)$. Next we show $T^*=\infty$ by contradiction.

Assume  $T^*<\infty$. By virtue of Proposition \ref{pro:0401n}, $T^{\max}>T^*$. Then one has
$$\mathcal{G}(T^*)=  \min\{\delta_1,\delta_2\}$$ by the continuity of $\mathcal{G}(t)$ on $(0,T^{\max})$.
On the other hand, noting that $\sqrt{\mathcal{G}_3(T^*)}\leq \sqrt{\mathcal{G}(T^*)}\leq\delta_1$, we see that $(\eta,u)$ satisfies \eqref{decadsyeste}
with $T^*$ in place of $T_1$. Thus, by the conditions $\sqrt{\|\eta_0\|_7^2+\|u_0\|_6^2}\leq \delta$ and \eqref{detlaci},
we further obtain a more precise estimate
\begin{equation} \label{Glsadfsa}
\sqrt{\mathcal{G}(T^*)}\leq \min\{\delta_1,\delta_2\}/2,
\end{equation}
 which is a contradiction. Hence $T^*=\infty$, and we obtain a global solution $(\eta,u)$ to the transformed {Parker problem},
 which satisfies \eqref{1.19}. Finally, the uniqueness of the global solution $(\eta,u)$ can be easily verified by the standard energy method.
 This completes the proof of Theorem \ref{thm3}.

\subsection{Proof of Theorem \ref{thm3orig}}
   To show Theorem \ref{thm3orig}, let $(\varrho_0,v_0,N_0)$ satisfy the assumptions of Theorem \ref{thm3orig}. Then, for sufficiently
small $\delta <0$, we can use the standard iteration method as in \cite[Proposition 3.1]{JFJSWWWN} to see that there are a $T>0$ and
a unique classical solution $(\varrho, v,N)\in C^0([0,T], H^6  )$ of the original {Parker problem}. Moreover,
$\partial_t^iv\in C^0([0,T],H^{6-2i})$ for $1\leq i\leq 3$, and
\begin{equation*} \begin{aligned}
\sup_{0\leq  t\leq T}&\left(\|(\varrho,N)\|_6^2+\sum_{k=0}^3 \|\partial_t^k {v}(t)\|_{6-2k}^2 \right)
+\int_0^T \sum_{k=0}^3 \|\partial_t^k {v}(t)\|_{7-2k}^2 \mm{d}t <\infty.
\end{aligned}
\end{equation*}

Let $u_0(y)=v_0(\zeta_0(y))$, where $\zeta_0$ is given by Theorem \ref{thm3orig}.
Noting that $\|u_0\|_{6}^2 \lesssim \|v_0\|_{6}^2$ by the first two initial conditions in Theorem \ref{thm3orig}
for sufficiently small $\delta$,  {one sees that $\eta_0:=\zeta_0-y$ and $u_0$
satisfy the first two conditions in Theorem \ref{thm3}, and $u_0|_{\partial\Omega}=0$
 by the the first two conditions in  Theorem \ref{thm3orig}.}
On the other hand, employing the last two conditions in Theorem \ref{thm3orig}, we can moreover verify that
$(\eta_0,u_0)$ satisfies the third condition in Theorem \ref{thm3} by following the  derivation of \eqref{peturbationequation}$_2$ from \eqref{0103}$_2$.
Hence, $(\eta_0,u_0)$ satisfies the three conditions in Theorem \ref{thm3} for sufficiently small $\delta$.

Now one can use the initial data $(\eta_0,u_0)$ and Theorem \ref{thm3} to construct a classical unique solution
$(\tilde{\eta},\tilde{u})\in C(\mathbb{R}^+_0,H^{7}\times H^6)$
to the transformed {Parker problem}, which satisfies the stability estimate \eqref{1.19} as $(\eta,u)$.
Thus, we can choose a sufficiently small $\delta$, so that there exists a continuously differential invertible
function $\tilde{\zeta}^{-1}$ of $\tilde{\zeta}:=y+\tilde{\eta}$.  Let
\begin{equation}  \label{usedopso}
\tilde{\sigma}:=(\bar{\rho}(\det\nabla \tilde{\zeta})^{-1})|_{y=\tilde{\zeta}^{-1}} -\bar{\rho}(x_3),\;\;
\tilde{v}:=\tilde{u}(\tilde{\zeta}^{-1},t),\;\; \tilde{N}:=(m(\det\nabla \tilde{\zeta})^{-1}
\partial_1\tilde{\zeta})|_{y=\tilde{\zeta}^{-1}}-me_3.   \end{equation}
In view of the relations
\begin{align} & 1\lesssim \det \nabla \tilde{\zeta}\lesssim 1,\\
 &\partial_t( f(y,t)|_{y=\tilde{\zeta}^{-1}})=\partial_t f(y,t)|_{y=\tilde{\zeta}^{-1}}+\partial_if(y,t)|_{y=\tilde{\zeta}^{-1}}\partial_t \tilde{\zeta}^{-1}_i,\\
\label{realdionsf}
&\partial_t \tilde{\zeta}^{-1} =- (\nabla \tilde{\zeta})^{-1}|_{y=\tilde{\zeta}^{-1}}\tilde{v} \mbox{ and }
 \nabla \tilde{\zeta}^{-1} = (\nabla \tilde{\zeta})^{-1}|_{y=\tilde{\zeta}^{-1}},
\end{align}
we see that $(\tilde{\sigma},\tilde{v},\tilde{N})$ is a classical solution to the original {Parker problem}.
Moreover, by a standard uniqueness proof, we have $(\tilde{\sigma},\tilde{v},\tilde{N})=(\varrho,v,N)$,
and therefore, one obtains the existence of a unique global solution in Theorem \ref{thm3orig}.
Finally, recalling that $(\tilde{\eta},\tilde{u})$ satisfies the stability estimate \eqref{1.19} as $(\eta,u)$, we make
use the relations \eqref{usedopso}--\eqref{realdionsf} to deduce the stability estimate \eqref{1.19st}.
This completes the proof of Theorem \ref{thm3orig}.

\section{Proof of Theorem \ref{thm:0101jfjsww101315}}\label{sec:06}
   In this section we prove Theorem \ref{thm:0101jfjsww101315}. To begin with, we introduce the instability of the linearized {Parker problem} \eqref{c0104}--\eqref{lincom} and the local well-posedness of strong solutions to the original {Parker problem} \eqref{0103}--\eqref{0105}.
   %%%%%%%%%%%%%%%%%%%%%%%%%%%%%%%%%%%%%%%%%%%%
\begin{pro}\label{thm:0201201622}
  Under the assumptions of Theorem \ref{thm:0101jfjsww101315},
 the equilibrium state $(\bar{\rho}, {0},\bar{M})$ of the linearized {Parker problem} \eqref{c0104}--\eqref{lincom}
 is unstable, that is, there is an unstable solution in the form
$$({\varrho}, v,N):=e^{\Lambda t}(-\mm{div}(\bar{\rho}\tilde{u})/\Lambda,\tilde{ {u}},(m \partial_1 \tilde{u}
 -\tilde{u}_3 \bar{M} '-\bar{M} \mm{div}\tilde{u})/\Lambda)$$
of \eqref{c0104}--\eqref{lincom}, where $\tilde{u}\in H^4\cap \mathcal{A}$ solves the boundary value problem:
 \begin{equation}
 \label{201604061413}      \left\{  \begin{array}{l}
\Lambda^2\bar{\rho}\tilde{u} =g\bar{\rho}'\tilde{u}_3e_3+\nabla (P'(\bar{\rho})\mm{div}
   ( \bar{\rho}\tilde{u}))+g\bar{\rho}\mm{div}\tilde{u} e_3+
   \lambda m (m \partial_1^2 \tilde{u} -\bar{M} \mm{div}\partial_1\tilde{u}) \\
   \qquad\qquad +\nabla ( \lambda m  (m \partial_2 \tilde{u}_2+m \partial_3 \tilde{u}_3+m'
\tilde{u}_3) )+\Lambda\mu_1\Delta\tilde{u}  +\Lambda\mu_2\nabla\mm{div}{\tilde{u}},\\
 \tilde{u}|_{\partial\Omega}=0   \end{array}\right.
\end{equation}
  with $\Lambda >0$ being a constant satisfying
\begin{equation}
\label{0111nn} \Lambda^2=\sup_{w\in\mathcal{A}}\mathcal{E}_{\mm{c}} (w,\Lambda)= {\mathcal{E}_{\mm{c}}  (\tilde{u},\Lambda)},
 \end{equation}
 where
$$\mathcal{A}:=\left\{w\in H^1_0~\bigg|~\int \bar{\rho}{w}^2\mm{d}x=1\right\}
$$ and
$${\mathcal{E}_{\mm{c}}(w,\Lambda)}:=E(w)-\Lambda  \int (\mu_1|\nabla w|^2+\mu_2|\mm{div}w|^2)\mm{d} x. $$
Moreover, \begin{equation}
 \label{201602081445MH}\tilde{u}_3\neq 0.
 \end{equation}
\end{pro}
\begin{pf} Proposition \ref{thm:0201201622} has been proved in \cite[Theorem 2.1]{JFJSJMFM}, where the considered domain
is a bounded $C^2$-domain and $g$ is a positive constant. Next, we denote the the bounded $C^2$-domain in \cite[Theorem 2.1]{JFJSJMFM}
by $\Omega'$, and show how to modify the proof in \cite{JFJSJMFM} to establish Proposition \ref{thm:0201201622} in a strip domain $\Omega$ here.

(1) Since ${\Omega}'$ is a $C^2$-bounded domain, $H^1({\Omega}')$ is compactly embedded in $L^2({\Omega}')$.
In the derivation of the first conclusion in \cite[Proposition 2.1]{JFJSJMFM}, Jiang et.al. used this fact to derive that
 \begin{equation}\label{201604061120}
 \begin{aligned}&\sup_{w\in \mathcal{A}}\mathcal{E}_{\mm{c}}(w,s)\leq \mathcal{E}_{\mm{c}}(\tilde{w},s)\;\;\mbox{ for the bounded domain }\;
 \Omega'  \end{aligned} \end{equation} and
 \begin{equation}   \label{201604061120new}
 \int_{{\Omega}'} \bar{\rho}\tilde{w}^2\mm{d}x=1,\end{equation}
where $\tilde{w}$ is the limit of a subsequence of a maximizing sequence
$\{w_n\}_{n=1}^\infty $ of $\sup_{w\in \mathcal{A}}\mathcal{E}_{\mm{c}}(w,s)$. More precisely,
$\{w_n\}_{n=1}^\infty \ni w_{n_m}  \rightharpoonup\tilde{w}\;\mbox{ weakly in }\; H^1_0(\Omega')$ and
\begin{equation}    \label{strong201604102123}
w_{n_{m}}\to \tilde{w}\mbox{ strongly in }L^2(\Omega')\;\mbox{ as }\; n_{m}\to \infty.
\end{equation}
Since $g\in C^6_0(\mathbb{R})$ and \eqref{strong201604102123} holds for a bounded subdomain of an unbounded domain,
\eqref{201604061120} still holds for a strip domain $\Omega$. However, by the lower semi-continuity of weak convergence,
we only have $\int_{\Omega}\bar{\rho}\tilde{w}^2\mm{d}x\leq 1$.

To show \eqref{201604061120new} with $\Omega$ in place of ${\Omega}'$, we argue by contradiction.
From the instability condition $\Xi>0$ we get for sufficiently small $s$ that $\sup_{w\in \mathcal{A}}\mathcal{E}_{\mm{c}}(w,s)>0$,
which implies $\int_{{\Omega}} \bar{\rho}\tilde{w}^2\mm{d}x\neq 0$.
Now we assume that $\int_{{\Omega}} \bar{\rho}\tilde{w}^2\mm{d}x<1$, then
 \begin{equation*}
 \begin{aligned}&{\mathcal{E}_{\mm{c}}(\tilde{w},s)}< \frac{\mathcal{E}_{\mm{c}}(\tilde{w},s)}{ \int_{{\Omega}}
 \bar{\rho}\tilde{w}^2\mm{d}x}\leq\sup_{w\in \mathcal{A}}\mathcal{E}_{\mm{c}}(w,s)\leq  \mathcal{E}_{\mm{c}}(\tilde{w},s),
\end{aligned} \end{equation*}
which is a contradiction. Hence, $\int_{\Omega}\bar{\rho}\tilde{w}^2\mm{d}x=1$.

(2) In view of the proof of Proposition 4.2, we observe that $u \in H^4({\Omega}')$
if $u\in H_0^1$ is a weak solution of the boundary value problem \eqref{201604061413}
with a bounded $C^4$-domain ${\Omega}'$ in place of $\Omega$. This conclusion also holds for the strip domain $\Omega$
by employing the same domain expansion technique as in the derivation of (3.7) in \cite{KYLH3295}.

(3) Jiang et.al. in \cite{JFJSJMFM} used Poincar\'e's inequality
$$\|w\|_0\lesssim  \|\nabla w\|_{0}\quad\mbox{ for any }w\in H^1_0({\Omega}')$$
to derive (3.7). Obviously, this inequality still holds for the strip domain $\Omega$, see \eqref{insoslai}.

With the above facts, we easily establish Proposition \ref{thm:0201201622} by following the proof of \cite[Theorem 2.1]{JFJSJMFM}.
\hfill $\Box$
\end{pf}
%%%%%%%%%%%%%%%%%%%%%%%%%%%%%%%%%%%%%%%%%%%%%%%%%%%%
\begin{pro} \label{pro:0401n20160203}
 Under the assumptions of Theorem \ref{thm:0101jfjsww101315}, for any given initial data $(\varrho_0,v_0,N_0)\in H^3$ satisfying
$\inf_{x\in\Omega}\{(\varrho_0 + \bar{\rho})(x)\}>0$, the compatibility conditions $\partial_t^i v(\cdot ,0)|_{\partial\Omega}=0$
for $i=0,1$, and $\mathrm{div}{N}_0=0$, there exist a $T^{\max}>0$ and a unique classical solution
$(\varrho, v, N)\in C^0([0,T^{\max}),H^3 )$ to the original {Parker problem}, where $(\varrho,v,N)$ enjoys the following properties:
\begin{align*}
v_t\in C^0([0,T^{\max}),H^1 )\cap L^2((0,T^{\max}),H^2 )\;\;\mbox{ and }\;\;
0 <\inf_{(x,t)\in {\Omega}\times (0,T^{\max})}\{\varrho({x},t)+\bar{\rho}\},
\end{align*}
where $T^{\max}$ denotes the maximal time of existence of the solution $(\varrho, v, {N})$. In addition, if
\begin{equation}
\label{infer201604061602}
\inf_{x\in \Omega}\{\bar{\rho}(x)\}/2 \leq (\varrho+\bar{\rho})(x,t) \leq 2\sup_{x\in \Omega}\{\bar{\rho}(x)\}\quad
\mbox{on }\;\Omega\times (0,T^{\max}), \end{equation}
then
\begin{equation}  \label{20160204nHK}
\begin{aligned}& \mathcal{S}(\varrho(t),v(t),N(t)):= \sqrt{\|v(t)\|_1^2 + \|(\varrho(t),N(t),\varrho_t,v_t,N_t)\|_0^2}  \\
 &  \leq C_3 \sqrt{\| (\varrho,N)(t) \|_1^2  +\|v(t)\|_2^2 +\|(\varrho,v,N)(t)\|_1^2\|(\varrho,v,N )(t)\|_2^2}
 \end{aligned}
 \end{equation}  for any $t\in (0,T^{\max})$, where the constant $C_3\geq 1$ only depends on $\Omega$ and the known physical parameters
 in the perturbation equations.
\end{pro}
\begin{pf} Proposition \ref{pro:0401n20160203} can be established by a standard iteration scheme as in \cite[Proposition 3.1]{JFJSWWWN},
and hence, we omit its proof here which involves only tedious calculations. Here we only show \eqref{20160204nHK}.

Recalling \eqref{equcomre} and the relation
$$\begin{aligned}P(\varrho+\bar{\rho})
=\bar{P}+P'(\bar{\rho}) \varrho + \int_{0}^{\varrho}(\varrho-z)P''(z+\bar{\rho})\mm{d}z, \end{aligned} $$
we can write \eqref{0103}$_2$ as follows
$$ \begin{aligned}  & (\varrho+\bar{\rho}){v}_t+(\varrho+\bar{\rho}){v}\cdot\nabla
v+\nabla \Big( P'(\bar{\rho})\varrho + \int_{0}^{\varrho}(\varrho -z)P''(z+\bar{\rho})\mm{d}z+\lambda (2mN_1+|N|^2)/2\Big) \\
&  =\mu_1\Delta v+\mu_2\nabla\mm{div} v+ \lambda(\bar{M}' N_3 + m\partial_1 N + N\cdot \nabla N)-\varrho g e_3.
\end{aligned} $$
Multiplying the above identity by $v_t $ in $L^2$ and using \eqref{infer201604061602}, we deduce that, for any $\epsilon>0$,
$$   \begin{aligned}
\int(\varrho+\bar{\rho})| v_t |^2\mm{d}x \leq &C_\epsilon\left(\|(\varrho, N) \|_1^2+\|v\|_2^2\right.
\\
&\left. + \|(\varrho,v,N)\|_1^2 \|(\varrho,v,N)\|_2^2\right)+\epsilon\|v_t\|_0^2,\end{aligned}$$
where the constant $C_\epsilon$ depends on $\epsilon$. So, one can use \eqref{infer201604061602} again to get
$$  \|v_t\|^2_0\leq\|v\|_2^2 + \|(\varrho, N)\|_1^2 + \|(\varrho,v,N)\|_1^2 \|(\varrho,v,N)\|_2^2.   $$
Using \eqref{0103}$_1$ and \eqref{0103}$_3$, we have
$$\|\varrho_t\|_0 \lesssim (1+\|\varrho\|_2)\|v\|_1,\quad\;\; \| N_t\|_0\lesssim (1+\| N\|_2)\|v\|_1. $$
From the above three estimates we get \eqref{20160204nHK} immediately.      \hfill$\Box$
\end{pf}

We next use Propositions \ref{thm:0201201622}--\ref{pro:0401n20160203} to construct unstable solutions to the original {Parker problem}.
To this end, we first use Proposition \ref{thm:0201201622} to give a solution of the linearized perturbation equations \eqref{lincom} in the form:
\begin{equation}\label{0501}
\left(\varrho^\mm{l}, v^\mm{l},N^\mm{l}\right)=e^{{\Lambda t}} (\tilde{\varrho}_0,\tilde{v}_0,\tilde{N}_0)\in H^3,
\end{equation}
where
$$(\tilde{\varrho}_0,\tilde{v}_0,\tilde{N}_0):=(-\mm{div}(\bar{\rho}\tilde{u})/\Lambda,\tilde{ {u}},(m \partial_1 \tilde{u}
 -\tilde{u}_3 \bar{M} '-\bar{M} \mm{div}\tilde{u})/\Lambda).$$
Now, denote
 \begin{equation}
 \label{201602061431MH}
 (\varrho^\mm{a}, v^{\mm{a}}, {N}^{\mm{a}}):=
\delta(\varrho^{\mm{l}},v^{\mm{l}}, {N}^{\mm{l}}),\;\;\quad \delta>0,\mbox{ a constant}.
\end{equation}
Then, $(\varrho^\mm{a}, v^{\mm{a}},N^{\mm{a}})$ solves the linearized perturbation equations \eqref{lincom} with initial data
$\delta (\tilde{\varrho}_0,\tilde{v}_0,\tilde{N}_0) $, and satisfies
 \begin{equation}
 \label{201602071912MH}\|(\varrho^\mm{a}, v^{\mm{a}}, {N}^{\mm{a}})\|_3\lesssim \delta e^{{\Lambda t}}.
 \end{equation}
In the rest of this paper, we call \eqref{201602061431MH} an approximate solution of the original {Parker problem} for fixed $\delta$.

Secondly, we modify  the initial data $\delta (\tilde{\varrho}_0,\tilde{v}_0,\tilde{N}_0) $  to construct a solution of
the original {Parker problem} for sufficiently small $\delta>0$.
%%%%%%%%%%%%%%%%%%%%%%%%%%%%%%%%%%%%%%%%%%%%%%%%%%%%%
\begin{lem}\label{lem:modfied} Let $(\tilde{\varrho}_0,\tilde{ v}_0, \tilde{ N}_0)$ be the same as in \eqref{0501}.
 Then there are an error function $v_\mm{r}$ and a constant ${\delta}_1\in (0,1)$ depending on $(\tilde{\varrho}_0,\tilde{ v}_0, \tilde{N}_0)$,
 such that for any $\delta\in (0, {\delta}_1)$,
\begin{enumerate}
  \item[(1)] The modified initial data
  \begin{equation}\label{mmmode04091215}(  {\varrho}_0^\delta,{v}_0^\delta,{N}_0^\delta ): =\delta
   (\tilde{\varrho}_0,\tilde{v}_0, \tilde{N}_0) + ( 0, \delta^2 v_\mm{r}, 0)
\end{equation}
satisfies the compatibility conditions on boundary of the {Parker problem} (\ref{0103})--(\ref{0105}).
\item[(2)] $ v_\mm{r} $ satisfies the following estimate:
$$ \|v_\mm{r}\|_{3} \leq C_4, $$
where the constant $C_4\geq 1$ depends on $\|(\tilde{\varrho}_0,\tilde{v}_0, \tilde{N}_0)\|_{3}$
and other physical parameters, but is independent of $\delta$.
\end{enumerate}  \end{lem}
\begin{pf}
Note that $(\tilde{\varrho}_0,\tilde{v}_0,\tilde{N}_0)$ satisfies
 $$\tilde{v}_0|_{\partial\Omega}= 0,\;\; (\mu_1\Delta \tilde{v}_0+\mu_2\nabla\mm{div}\tilde{v}_0 +
  \lambda\bar{M}'\tilde{N}^0_3+\lambda m\partial_1\tilde{N}_0-\tilde{\varrho}_0 ge_3 -\nabla (P'(\bar{\rho})\tilde{\varrho}_0+\lambda
 m \tilde{N}^0_1))|_{\partial\Omega}=0.$$
 Hence, if the modified initial data satisfy \eqref{mmmode04091215}, then we expect $v_\mm{r}$ to satisfy the following problem:
\begin{equation} \left\{\begin{array}{l}
 \mu_1\Delta v_\mm{r} +\mu_2\nabla\mm{div} v_\mm{r} -\delta^2(\delta \tilde{\varrho}_0+\bar{\rho}) v_\mm{r}
  \cdot\nabla
v_\mm{r}-\delta(\delta \tilde{\varrho}_0+\bar{\rho})
%(\tilde{ v}_0\cdot\nabla {u}_\mm{r}+
v_\mm{r}\cdot\nabla\tilde{v}_0  \\
%(\delta \tilde{\varrho}_0+\bar{\rho})\tilde{v}_0\cdot\nabla\tilde{v}_0+
\quad=\nabla
 \left(\int_0^{\tilde{\varrho}_0}( \tilde{\varrho}_0-z)P''(\delta z+\bar{\rho})
\mm{d}z+\frac{\lambda |\tilde{N}_0|^2}{2}\right)-\lambda \tilde{ N}_0\cdot \nabla \tilde{N}_0=:F(\tilde{\varrho}_0,\tilde{v}_0,\tilde{N}_0), \\[1mm]
  v_\mm{r}|_{\partial\Omega} = 0.
 \end{array} \right.   \label{js204091215} \end{equation}
Thus the modified initial data naturally satisfy the compatibility conditions on boundary.

Next we shall look for a solution $v_\mm{r}$ to the boundary problem \eqref{js204091215} when $\delta$ is sufficiently small.
%%%
We begin with the linearization of \eqref{js204091215} which reads as
\begin{equation}\label{elliequation04091215}
  \begin{aligned}&\mu_1\Delta v_\mm{r} +\mu_2\nabla\mm{div} v_\mm{r}
  =F(\tilde{\varrho}_0,\tilde{v}_0,\tilde{N}_0)+\delta^2(\delta \tilde{\varrho}_0+\bar{\rho})w
  \cdot\nabla
w+\delta(\delta \tilde{\varrho}_0+\bar{\rho})w\cdot\nabla\tilde{v}_0 \end{aligned}
\end{equation}
with boundary condition
\begin{equation}\label{boundery004091215}
v_\mm{r}|_{\partial\Omega}= 0.
\end{equation}

Let $w\in H^{3}$, then it follows from the elliptic theory that there is a solution $v_\mm{r}$ of
\eqref{elliequation04091215}--\eqref{boundery004091215} satisfying
$$\|v_\mm{r}\|_{3} \leq \|F(\tilde{\varrho}_0,\tilde{v}_0,\tilde{N}_0)+\delta^2(\delta \tilde{\varrho}_0
+\bar{\rho})w\!\cdot\!\nabla w +\delta(\delta \tilde{\varrho}_0+\bar{\rho})w\!\cdot\!\nabla\tilde{v}_0 \|_{1}
\leq c_{1}(1+ \delta^2 \|w\|_{2}^2),\;\;\delta\in (0,1), $$
where $c_1$ depends on $\|(\tilde{\varrho}_0,\tilde{v}_0, \tilde{N}_0)\|_{3}$ and the physical parameters.
Now, taking $\delta\leq \min\{(2c_1)^{-1},1\}$, we have for any $\|w\|_{ {3}}\leq 2c_1$ that
$\|v_\mm{r}\|_3 \leq 2c_1$. Therefore, one can construct an approximate function sequence $v_\mm{r}^n$, such that
\begin{equation*}  \begin{aligned}
& \mu_1\Delta v_\mm{r}^{n+1} +\mu_2\nabla\mm{div} v_\mm{r}^{n+1}
-\delta^2(\delta \tilde{\varrho}_0+\bar{\rho})v_\mm{r}^{n}  \cdot\nabla
v_\mm{r}^{n}-\delta(\delta \tilde{\varrho}_0+\bar{\rho})v_\mm{r}^{n}\cdot\nabla\tilde{v}_0 =F(\tilde{\varrho}_0,\tilde{v}_0,\tilde{N}_0),
\end{aligned} \end{equation*}
and for any $n$, $\|v_\mm{r}^n\|_{{3}}\leq 2c_1$ and $\|v_\mm{r}^{n+1}-v_\mm{r}^{n}\|_{{3}}\leq c_2\delta\|{v}_\mm{r}^{n}-v_\mm{r}^{n-1}\|_{3}$
for some constant $c_2$ independent of $\delta$ and $n$.

Finally, we choose a sufficiently small $\delta$ so that $c_2\delta<1$, and use then a compactness argument to get
a limit function which solves the nonlinear boundary problem \eqref{js204091215}. Moreover $\|v_\mm{r}\|_{3}\leq 2c_1$.
Thus we have proved Lemma \ref{lem:modfied}.    \hfill$\Box$
\end{pf}

Now, in view of the condition $\inf_{x\in {\Omega}}\{\bar{\rho}(x)\}>0$ and the embedding inequality \eqref{esmmdforinfty},
we can choose a sufficiently small $\tau\in (0,\delta_1)$, such that
\begin{equation}
\label{infdensity}
\frac{\inf_{ x\in {\Omega}}\{\bar{\rho}( x)\}}{2}\leq
\inf_{ x\in {\Omega}}\{\varrho_0^\delta( x)+\bar{\rho}( x)\}\leq 2
{\sup_{ x\in {\Omega}}\{\bar{\rho}( x)\}}\quad\mbox{for any }\delta\in (0,\tau).
\end{equation}
Thus, by virtue of Proposition \ref{pro:0401n20160203}, for any given $\delta<\tau$, there exists a unique local-in-time classical solution
$(\varrho^\delta,v^\delta,N^\delta)\in C^0([0,T^{\max}),H^3)$ to the original {Parker problem}, emanating from the initial data $(\varrho_0^\delta,v_0^\delta,N_0^\delta)$.

Thirdly, we estimate the error between the approximate solution $(\varrho^\mm{a},v^{\mm{a}},N^{\mm{a}})$ and $(\varrho^\delta,v^\delta,N^\delta)$. Denote
$(\varrho^{\mathrm{d}}, v^{\mathrm{d}},N^{\mathrm{d}})=(\varrho^{\delta},v^{\delta},N^{\delta})-(\varrho^\mm{a},v^{\mm{a}},N^{\mm{a}} )$.
Then $(\varrho^{\mathrm{d}}, v^{\mathrm{d}},{ {N}}^{\mathrm{d}})$ satisfies the following error equations:
\begin{equation}\label{newfor1232}\left\{\begin{array}{ll}
  \varrho_t^{\mathrm{d}}+\mm{div}(\bar{\rho}v^{\mathrm{d}})=F^\delta, \\[1mm]
  (\varrho^\delta+\bar{\rho})v_t^{\mathrm{d}}+\nabla ( P'(\bar{\rho})\varrho^{\mathrm{d}}+ \lambda
 m  N_1^{\mathrm{d}})-\mu_1\Delta v^{\mathrm{d}}-\mu_2\nabla\mm{div}v^{\mathrm{d}}\\
 \quad  - \lambda \bar{M}' N_3^{\mathrm{d}}-\lambda m \partial_1  N^{\mathrm{d}}+\varrho  g^{\mathrm{d}}e_3 ={G}^\delta,
   \\[1mm]
 {N}_t^{\mathrm{d}} =  m  \partial_1 v^{\mathrm{d}}  -v_3^{\mathrm{d}}  \bar{M} '-\bar{M} \mm{div}v^{\mathrm{d}}   + {H}^\delta ,\\[1mm]
 \mathrm{div} {N}^{\mathrm{d}}={0},
 \end{array}\right.\end{equation}
 where \begin{equation*}\begin{aligned}
 & F^\delta:= -\mm{div}( \varrho^{\delta}v^{\delta}) ,\\
 &
 \begin{aligned}G^\delta:= &\lambda N^{\delta} \cdot\nabla {N}^{\delta} - \nabla
 \left(\int_{0}^{ \varrho^\delta}( \varrho^\delta  -z)P''(z+\bar{\rho})\mm{d}z +\lambda|N^\delta|^2/2\right)
  \\
&   -( \varrho^{\delta}+\bar{\rho}) v^{\delta}\cdot\nabla v^{\delta}-\varrho^\delta v^{\mm{a}}_t,
\end{aligned}\\
 & H^\delta: =N^{\delta}\cdot\nabla v^{\delta} -v^{\delta} \cdot \nabla N^{\delta} - N^{\delta} \mm{div}v^{\delta},
                  \end{aligned}         \end{equation*}
 with initial conditions
 \begin{equation}   \label{errorconditiopn}
  (\varrho^{\mathrm{d}}(0), v^{\mathrm{d}}(0), N^{\mathrm{d}}(0)):= (0,\delta^2 v_{\mathrm{r}},0)     \end{equation}
   and
$$\mm{div}N^{\mm{d}}(0)=0.$$
Moreover,  we can deduce the following error estimate from
the error equations, which plays a key role in the proof of Theorem \ref{thm:0101jfjsww101315}.
\begin{lem}\label{lem:0401}
Let $\epsilon\in (0,1)$ and $\delta\in (0,\tau)\subset(0,1)$. Denote
\begin{eqnarray*} &&  C_5:=(1+\Lambda) (\|(\tilde{\varrho}_0,\tilde{v}_0, \tilde{N}_0) \|_2+ 2C_4)^2, \\[1mm]
&& T^*:=\sup\Big\{t\in (0,T^{\max})|\;\mathcal{S} (\varrho^\delta(\tau),v^\delta(\tau),N^\delta(\tau))
   + \|  v_{s}^\delta \|_{L^2((0,{\tau}),H^1)}  \\
   && \qquad \leq 2 C_3C_5\delta e^{\Lambda \tau}\quad\mbox{for any }\;\tau\in [0,t]\Big\}.
   \end{eqnarray*}
If $T^{\max}=\infty$, and $(\varrho^\delta,v^\delta,N^\delta)$ satisfies \eqref{infer201604061602} and
\begin{equation}
\label{201602071MH}
\|(\varrho^\delta,v^\delta,N^\delta)(t)\|_3\leq \epsilon\mbox{ on }[0,\infty),
 \end{equation}
 then there is a constant $C_6$, such that
\begin{equation}\begin{aligned}\label{0508}
\mathcal{S}(\varrho^\mm{d}(t),v^{\mm{d}}(t),N^{\mm{d}}(t)) + \|v_{\tau}^{\mm{d}} \|_{L^2((0,t),H^1)}
\leq C_6 \sqrt{f(\epsilon, \delta, t)}\quad\mbox{ for any }\; t\in (0,T^*],\end{aligned} \end{equation}
where $\mathcal{S}(\varrho^\mm{d}(t),v^{\mm{d}}(t),N^{\mm{d}}(t))$ is defined by \eqref{20160204nHK} and
$f(\epsilon,\delta, t):=\delta^3 e^{3\Lambda t}+ \epsilon^{\frac{2}{3}}\delta^{ \frac{7}{3} }e^{\frac{7\Lambda t}{3}}.$
\end{lem}
%%%%%%%%%%%%%%%%%%%%%%%%%%%%%%%%%%%%%
\begin{rem}   \label{rem:06011418}
It should be remarked that \eqref{infer201604061602} automatically holds for sufficiently small $\epsilon$. In fact,
by virtue of \eqref{esmmdforinfty}, there is a positive constant $\epsilon_0$, such that
 \eqref{infer201604061602} holds for any $\varrho^\delta$ satisfying $\|\varrho^\delta\|_3\leq\epsilon_0$.
\end{rem}
\begin{pf} From Lemma \ref{lem:modfied} one gets
$$ \begin{aligned}
&\sqrt{ 1 +\|(\varrho_0^\delta,v_0^\delta,N_0^\delta)\|_1^2}\|(\varrho_0^\delta,v_0^\delta,N_0^\delta )\|_2\\
&\leq \delta\sqrt{1+(\|(\tilde{\varrho}_0,\tilde{v}_0,
\tilde{N}_0) \|_2+ C_4)^2 }(\|(\tilde{\varrho}_0,\tilde{v}_0,\tilde{N}_0)\|_2 + C_4) \leq C_5\delta.
\end{aligned}$$
Thus, in view of the regularity of $(\varrho^\delta,v^\delta,N^\delta)$ and the estimate \eqref{20160204nHK}, we see that $T^*>0$.
Moreover, by the definition of $T^*$,
\begin{equation}\label{20160207MH}
\mathcal{S} (\varrho^\delta(t),v^\delta(t),N^\delta(t)) +\|v_{\tau}^\delta \|_{L^2((0,{t}),H^1)}
    \leq 2 C_3C_5\delta e^{\Lambda t}\;\;\mbox{ on } [0,T^*],     \end{equation}
or
\begin{equation}\label{201602092200MH}
\mathcal{S} (\varrho^\delta(T^*),v^\delta(T^*),N^\delta(T^*))
   + \|  v_{\tau}^\delta \|_{L^2((0,T^*),H^1)}
  = 2 C_3C_5\delta e^{\Lambda T^*}\;\mbox{ in the case of } T^*<\infty.
    \end{equation}
%%%%%
We next show \eqref{0508}. By \eqref{newfor1232}$_2$ we find that
\begin{equation}\label{n0310}   \begin{aligned}
&  \frac{\mm{d}}{\mm{d}t}\int (\varrho^\delta+\bar{\rho})| v_t^{\mathrm{d}}|^2
\mathrm{d} x=2\int((\varrho^\delta+\bar{\rho})v_t^{\mm{d}})_t\cdot v_t^{\mm{d}}\mm{d}x-\int \varrho^\delta_t|v_t^{\mm{d}}|^2\mm{d}x\\
&= 2\int ( P'(\bar{\rho})\varrho^{\mathrm{d}}_t+ \lambda m \partial_t N_1^{\mathrm{d}})\mm{div}v^{\mathrm{d}}_t\mm{d}x
+2\int( \lambda \bar{M}' \partial_t N_3^{\mathrm{d}}+  \lambda m\partial_1 N^{\mathrm{d}}_t- g\varrho^{\mm{d}}_t e_3)\cdot v^{\mathrm{d}}_t\mm{d}x   \\
&\qquad +\int (2G_t^\delta-\varrho^\delta_t v_t^{\mm{d}}) \cdot v^{\mathrm{d}}_t\mm{d}x -2\int (\mu_1|\nabla v^{\mathrm{d}}_t|^2+\mu_2| \mm{div}v^{\mathrm{d}}_t|^2)\mm{d}x \\
& =:\tilde{R}_1(t)+\tilde{R}_2(t)+R_3(t)+R_4(t).
\end{aligned}\end{equation}
On the other hand, arguing similarly to that for \eqref{estrhoenn201604111600n}--\eqref{estrhoenn201604111601n}, we use \eqref{newfor1232}$_1$, \eqref{newfor1232}$_3$, $\eqref{comsteady}$
to infer that the first two integrals $\tilde{R}_1(t)$ and $\tilde{R}_2(t)$ can be written as follows.
\begin{equation*}  \begin{aligned}
&\tilde{R}_1(t)+\tilde{R}_2(t) = 2  \int  (-P'(\bar{\rho})\mm{div}( \bar{\rho}v^{\mathrm{d}}) +\lambda m (m\partial_1 v^{\mathrm{d}}_1
 -m'v^{\mathrm{d}}_3-m\mm{div}v^{\mathrm{d}}))\mm{div}v^{\mathrm{d}}_t\mm{d}x \\
&\qquad  +\int (2\lambda m^2 (\partial_1^2 v^{\mathrm{d}}\cdot  v^{\mathrm{d}}_t -\partial_1\mm{div}v^{\mathrm{d}}\partial_tv^{\mathrm{d}}_1)
 +g\mm{div}( \bar{\rho}v^{\mathrm{d}}) \partial_t v^{\mathrm{d}}_3) \mm{d}x +R_2(t) \\
& \quad = \frac{\mm{d}}{\mm{d}t}E (v^{\mathrm{d}}(t))+R_2(t),
\end{aligned}\end{equation*}
where
$$R_2(t):= 2\int F^\delta (P'(\bar{\rho}) \mm{div}v^{\mathrm{d}}_t -g\partial_tv^{\mathrm{d}}_3)\mm{d}x + 2\lambda
\int (m {H}^\delta_1\mm{div}v^{\mathrm{d}}_t +m' {H}^\delta_3\partial_tv^{\mathrm{d}}_1 +m\partial_1 {H}^\delta\cdot v^{\mathrm{d}}_t)\mm{d}x.$$
Thus the equality \eqref{n0310} can be rewritten as
\begin{equation}\label{nnn0314P}
\frac{\mm{d}}{\mm{d}t} \int \left((\varrho^\delta+\bar{\rho})| v_t^{\mathrm{d}}|^2-E(v^{\mathrm{d}}(t))\right)
\mathrm{d} x +2\int (\mu_1|\nabla v^{\mathrm{d}}_t|^2  + \mu_2 | \mm{div}v^{\mathrm{d}}_t|^2)\mm{d}x  = R_2(t)+R_3(t).
\end{equation}
Recalling that $v^\mm{d}(0)=\delta^2 v_{\mm{r}}$, we integrate \eqref{nnn0314P} in time from $0$ to $t$ to get
\begin{equation}\label{0314}
\begin{aligned}     \|\sqrt{\varrho^\delta+\bar{\rho}} v_t^\mm{d}(t)\|^2_{0}+2\int_0^t
(\mu_1\|\nabla v^{\mathrm{d}}_\tau\|^2_0  + \mu_2 \| \mm{div}v^{\mathrm{d}}_\tau\|^2_0)\mm{d}\tau =E(v^{\mathrm{d}}(t))
 + \sum_{i=1}^3\mathfrak{R}_i,
\end{aligned}\end{equation}
where
\begin{equation*}
 \mathfrak{R}_1:= \int(\varrho^\delta+\bar{\rho})| v_t^{\mathrm{d}}|^2\mathrm{d} x \bigg|_{t=0}
-E(\delta^2 v_{\mm{r}}),\quad\mathfrak{R}_2:=\int_0^t R_2 (\tau)\mm{d}\tau ,\;\mbox{ and }\;
 \mathfrak{R}_3:=\int_0^t R_3 (\tau)\mm{d}\tau
\end{equation*}
are bounded from below.

Multiplying \eqref{newfor1232}$_2$ by $v_t^\mm{d}$ in $L^2$, one gets
$$\begin{aligned}
&\int(\varrho^\delta+\bar{\rho})| v_t^{\mathrm{d}}|^2\mm{d}x\\
&  =\int (\mu_1\Delta v^{\mathrm{d}}+\mu_2\nabla\mm{div}v^{\mathrm{d}}
  +\lambda \bar{M}' N_3^{\mathrm{d}}+\lambda m \partial_1  N^{\mathrm{d}}-\varrho  g^{\mathrm{d}}e_3-\nabla ( P'(\bar{\rho})\varrho^{\mathrm{d}}+ \lambda
 m  N_1^{\mathrm{d}}) +G^\delta)\cdot v_t^{\mathrm{d}}\mm{d}x.
 \end{aligned}  $$
Applying \eqref{infer201604061602} and Cauchy-Schwarz's inequality, we obtain
\begin{equation}\label{appesimtsofu12}       \|\sqrt{\varrho^\delta+\bar{\rho}} v_t^{\mathrm{d}}\|^2_0
  \lesssim \|(\varrho^\mm{d},N^{\mm{d}})\|_{1}^2+\| v^{\mm{d}}\|_{2}^2+\|G^\delta\|_0^2. \end{equation}
From the definition of $v_t^{\mm{a}}$ it follows that
\begin{equation}\label{appesimtsofu1857}
\|\partial_{t}^j v^\mm{a}\|_{ k }=\Lambda^j \delta e^{\Lambda t}\|\tilde{v}_0\|_{k }
\lesssim \Lambda^j \delta e^{\Lambda t}\quad\mbox{ for }0\leq k,\, j\leq 3.
\end{equation}
So, using \eqref{appesimtsofu1857}, \eqref{20160207MH}, \eqref{201602071MH},  \eqref{infer201604061602},
the interpolation inequality \eqref{inteposlai}, and the Nirenberg interpolation inequality (see \cite[5.8 Theorem]{ARAJJFF})
\begin{equation}   \label{201602091446MH1}
\|w\|_{L^\infty}\lesssim \|w\|_{0}^{\frac{1}{2}}\| w\|_3^{\frac{1}{2}},
\end{equation}
we arrive at
\begin{equation} \label{appesofu12}
\begin{aligned}    \|G^\delta \|_{0}^2  & \lesssim
  \|v^\delta\|_{L^\infty }^2 \|v^\delta\|_{ 1 }^2+\|\varrho^\delta\|_{0 }^2\|v_t^{\mm{a}}\|_{ 2 }^2
  + \|(\varrho^\delta,N^\delta)\|_{ L^\infty}^2  \|(\varrho^\delta,N^\delta)\|_1^2 \\
& \lesssim \|v^\delta\|_0 \|v^\delta\|_3\|v^\delta\|_1^2 +\|(\varrho^\delta,N^\delta)\|_0^{\frac{7}{3}}\|(\varrho^\delta,N^\delta)\|_3^{\frac{5}{3}}
+\delta^4 e^{4\Lambda t} \\
&\lesssim\delta^4 e^{4\Lambda t} +\epsilon \delta^{3 } e^{{3\Lambda t}} +\epsilon^{\frac{5}{3}}\delta^{\frac{7}{3}}e^{{ \frac{7\Lambda t}{3} }}.
 \end{aligned}
 \end{equation}

Chaining the estimates \eqref{appesimtsofu12} and \eqref{appesofu12} together and taking then to the limit as $t\rightarrow 0$,
we apply Lemma \ref{lem:modfied} and \eqref{errorconditiopn} to obtain the following estimate on $\mathfrak{R}_1$:
\begin{equation}\label{estimeR1}\begin{aligned}
\mathfrak{R}_1=&\lim_{t\rightarrow 0}\|\sqrt{\varrho^\delta+\bar{\rho}} v_t^{\mathrm{d}}(t)\|^2_0 -E(\delta^2 v_{\mm{r}}) \\
\lesssim &\lim_{t\rightarrow 0}\left( \|(\varrho^\mm{d},N^{\mm{d}})\|_{1}^2 +\| v^{\mm{d}}\|_{2}^2+\delta^4 e^{4\Lambda t}
+\epsilon \delta^{3 } e^{{3\Lambda t}} +\epsilon^{\frac{5}{3}}\delta^{\frac{7}{3}} e^{{\frac{7\Lambda t}{3}}}\right) +\delta^4 \\
\lesssim &\delta^4 +\epsilon \delta^{3} +\epsilon^{\frac{5}{3}}\delta^{\frac{7}{3}}\lesssim\delta^3
  +\epsilon^{\frac{2}{3}}\delta^{ \frac{7}{3} } .
\end{aligned}\end{equation}

To bound $\mathfrak{R}_2(t)$, recalling the definition of $\mathfrak{R}_2(t)$, we have
\begin{equation} \label{201604112016n}
\mathfrak{R}_2(t) \lesssim    \int_0^t(\|F^\delta\|_0+\|H^\delta\|_0 )\|\partial_tv^{\mathrm{d}}\|_1 \mm{d}\tau.
 \end{equation}
By the interpolation inequality \cite[Theorem 1.49]{NASII04} in $L^p$,
\eqref{inteposlai} in $H^k$ and the embedding inequality \eqref{esmmdforinftfdsdy}, we see that
$$\|\nabla w\|_{L^3}\lesssim \|\nabla w\|_{0}^{\frac{1}{2}} \|\nabla  w\|_{L^6}^{\frac{1}{2}}\lesssim \| w\|_{1}^{\frac{1}{2}}
\|w\|_{2}^{\frac{1}{2}}\lesssim \|  w\|_{0}^{\frac{1}{2}} \|w\|_{3}^{\frac{1}{2}}, $$
from which, \eqref{201602091446MH1} and the embedding inequality it follows that
\begin{equation}  \label{201604112016n19212000n}
\begin{aligned}
\|F^\delta\|_0+\|H^\delta\|_0\lesssim &  \|v^\delta\|_1\|\nabla (\varrho^\delta,N^\delta)\|_{L^3}
+\|(\varrho^\delta,N^\delta)\|_{L^\infty}\|\nabla v^\delta\|_{0}\\
\lesssim& \|v^\delta\|_1\|(\varrho^\delta,N^\delta)\|_{0}^{\frac{1}{2}} \|(\varrho^\delta,N^\delta)\|_3^{\frac{1}{2}}
+\|(\varrho^\delta,N^\delta)\|_{0}^{\frac{1}{2}} \|  (\varrho^\delta,N^\delta)\|_{3}^{\frac{1}{2}}\|v^\delta\|_{1}  \\
 \lesssim & \epsilon^{\frac{1}{2}} \delta^\frac{3}{2} e^{\frac{3\Lambda t}{2}}.  \end{aligned}  \end{equation}
%%%%%
Plugging (\ref{201604112016n19212000n}) into \eqref{201604112016n}, we make use of \eqref{appesimtsofu1857}, \eqref{20160207MH} and H\"older's inequality to conclude
\begin{equation}  \label{201604112016n1921}
\begin{aligned}  \mathfrak{R}_2(t) \lesssim & \left(\int_0^t \epsilon  \delta^3 e^{{3\Lambda t}} \mm{d}\tau
 \right)^{\frac{1}{2}} \left(\int_0^t \| v_\tau^{\mm{a}}\|_1^2 \mm{d}\tau+\int_0^t \| v_\tau^{\delta }\|_1^2 \mm{d}\tau \right)^{\frac{1}{2}}
 \lesssim   \epsilon^{\frac{1}{2}} \delta^\frac{5}{2} e^{\frac{5\Lambda t}{2}}.
 \end{aligned}  \end{equation}

Now, to control $\mathfrak{R}_3(t)$, we integrate by parts and use H\"older's inequality to find that
$$\begin{aligned}
\mathfrak{R}_{3}(t)= &-2\lambda\int_0^t\int ( N^{\delta}_\tau\cdot\nabla v^{\mathrm{d}}_\tau \cdot {N}^{\delta}
 +{N}^{\delta} \cdot\nabla   v^{\mathrm{d}}_\tau \cdot {N}^{\delta}_\tau -  N^\delta \cdot N^\delta_{\tau} \mm{div} v^{\mathrm{d}}_\tau/2)\mm{d}x\mm{d}\tau \\
&-2\int_0^t \int \left((\varrho^{\delta}+\bar{\rho}) v^{\delta}_t\cdot\nabla v^{\delta} +(\varrho^{\delta}+\bar{\rho}) v^{\delta}\cdot\nabla
   v^{\delta}_t+\varrho^\delta v^{\mm{a}}_{tt} \right)\cdot v^{\mathrm{d}}_\tau\mm{d}x\mm{d}\tau \\
     &-\int_0^t \int \varrho^\delta_\tau\left(\left( v^{\mm{a}}_\tau +v_\tau^{\delta}+2 v^{\delta}\cdot\nabla
   v^{\delta}\right)\cdot v^{\mathrm{d}}_\tau -2\int_{0}^{ \varrho^\delta} P''(z + \bar{\rho})\mm{d}z\mm{div}v_\tau^{\mm{d}}\right)\mm{d}x\mm{d}\tau \\
= &\sum_{i=1}^3\mathfrak{R}_{3,i}(t),
\end{aligned} $$
where the terms on the right hand side can be estimated as follows, employing \eqref{201602091446MH1}, \eqref{appesimtsofu1857},
\eqref{20160207MH} and \eqref{201602071MH}.
\begin{equation}\label{esimtaes11}
\begin{aligned}
\mathfrak{R}_{3,1}(t) \lesssim & \int_0^t  \|N_\tau^\delta\|_0 \|N^\delta\|_\infty\| v_\tau^{\mm{d}}\|_1\mm{d}\tau \\
\lesssim &\left(\int_0^t  \|N_\tau^\delta\|_0^2  \|N^\delta\|_0 \|N^\delta\|_3  \mm{d}\tau \right)^{\frac{1}{2}} \left(\int_0^t \| v_\tau^{\mm{a}}\|_1^2
 \mm{d}\tau+\int_0^t \| v_\tau^{\delta }\|_1^2  \mm{d}\tau \right)^{\frac{1}{2}} \\
\lesssim & \epsilon^{\frac{1}{2}} \delta^\frac{5}{2} e^{\frac{5\Lambda t}{2} },
\end{aligned}  \end{equation}
while one can utilize \eqref{appesimtsofu1857}, \eqref{20160207MH}, \eqref{infer201604061602}, and \eqref{esmmdforinftfdsdy} to deduce that
\begin{equation}\label{esimtaes11newe}
\begin{aligned}    \mathfrak{R}_{3,2}(t)\lesssim & \int_0^t (\|\varrho^\delta \|_0 \|v^\mm{a}_{\tau\tau}\|_1 +\|v^\delta\|_1 \|v_\tau^\delta\|_1)
\|v^\mm{d}_\tau\|_1 \mm{d}\tau \\
 \lesssim & \int_0^t[\delta^2 e^{2\Lambda\tau} (\delta e^{\Lambda\tau} +\|v^\delta_\tau\|_1)
+ \delta e^{ \Lambda \tau}\|  v_\tau^\delta\|_{1}^2]\mm{d}\tau   \\
 \lesssim & \int_0^t\delta e^{\Lambda\tau} (\delta^2 e^{2\Lambda \tau} +\|v_\tau^\delta\|_{1}^2)\mm{d}\tau\lesssim\delta^3 e^{3\Lambda t}.
\end{aligned}    \end{equation}
Finally, similar to the derivation of \eqref{esimtaes11} and \eqref{esimtaes11newe}, one has
\begin{equation}\label{201602071414MH}
\begin{aligned}
\mathfrak{R}_{3,3}(t)\lesssim  &\int_0^t \|\varrho^\delta_\tau\|_0\left( \|v^{\mm{a}}_\tau\|_{L^3}+
 \|v^\delta_\tau\|_{L^3}　+ \|\varrho^{\delta}\|_{L^\infty} +\|v^\delta\|_1 \|v^\delta\|_2 \right)\| v^{\mm{d}}_\tau\|_1 \mm{d}\tau  \\
 \lesssim &\int_0^t \delta e^{ \Lambda \tau}\left( \|\varrho^{\delta}\|_0^{\frac{1}{2}}\|\varrho^{\delta}\|_3^{\frac{1}{2}} \|v_\tau^{\mm{d}}\|_{1}
  +\delta^2 e^{2\Lambda \tau}  +\| v_\tau^\delta\|_{1}^2\right)\mm{d}\tau \\
\lesssim &\delta^\frac{5}{2} e^{\frac{5\Lambda t}{2}} \left(\epsilon^\frac{1}{2}+\delta^\frac{1}{2}e^{\frac{\Lambda t}{2}}
\right). \end{aligned}\end{equation}

Thus, summing up the estimates \eqref{esimtaes11}--\eqref{201602071414MH}, \eqref{201604112016n1921}
and \eqref{estimeR1}, we use Young's inequality to infer
\begin{equation}\label{estimateforhigher}
\begin{aligned}
\sum_{i=1}^3\mathfrak{R}_i\lesssim \delta^3 e^{3\Lambda t}
 +
\epsilon^{\frac{1}{2}} \delta^\frac{5}{2} e^{\frac{5\Lambda t}{2} }
 +\epsilon^{\frac{2}{3}}\delta^{ \frac{7}{3} }
\lesssim f(\epsilon, \delta, t).
\end{aligned}\end{equation}
Combining \eqref{0314} with \eqref{estimateforhigher}, one obtains
$$  \|\sqrt{\varrho^\delta+\bar{\rho}} v_t^\mm{d}(t)\|^2_{0}+2 \int_0^t(\mu_1\|\nabla v^{\mathrm{d}}_\tau\|^2_0
 + \mu_2 \| \mm{div}v^{\mathrm{d}}_\tau\|^2_0)\mm{d}\tau  \leq E(v^{\mathrm{d}}(t))+ cf(\epsilon, \delta, t). $$
Thanks to \eqref{0111nn}, we have
\begin{equation*}\begin{aligned}
 E(v^{\mathrm{d}})\leq &\Lambda^2{\int\bar{\rho}|v^{\mathrm{d}} |^2\mathrm{d}x} + \Lambda(\mu_1\|\nabla v^{\mathrm{d}}\|^2_0
 + \mu_2 \| \mm{div}v^{\mathrm{d}}  \|^2_0) \\
= &\Lambda^2{\int(\varrho^\delta  +\bar{\rho})|v^{\mathrm{d}}|^2\mathrm{d} x} + \Lambda(\mu_1\|\nabla v^{\mathrm{d}}\|^2_0
 + \mu_2 \| \mm{div}v^{\mathrm{d}}\|^2_0) -\Lambda^2{\int\varrho^\delta|v^{\mathrm{d}}|^2\mathrm{d} x} .
\end{aligned}\end{equation*}
On the other hand,
$$\begin{aligned}
{\int\varrho^\delta|v^{\mathrm{d}}|^2\mathrm{d} x}\lesssim & \|\varrho^\delta
\|_{L^\infty}\|v^{\mm{d}}\|_{0}^2\lesssim \|\varrho^\delta\|_0^{\frac{1}{2}}\|\varrho^\delta\|_3^{\frac{1}{2}}
\delta^{2} e^{{2\Lambda t} }\lesssim \epsilon^{\frac{1}{2}} \delta^\frac{5}{2} e^{\frac{5\Lambda t}{2} }.
\end{aligned}$$
Combining the above three inequalities together, we arrive at
\begin{equation}\label{new0311}
\begin{aligned} &\|\sqrt{\varrho^\delta+\bar{\rho}} v_t^\mm{d}(t)\|^2_{0}+2\int_0^t(\mu_1\|\nabla v^{\mathrm{d}}_\tau\|^2_0
 + \mu_2 \| \mm{div}v^{\mathrm{d}}_\tau\|^2_0)\mm{d}\tau \\
&\leq {\Lambda^2} \|\sqrt{\varrho^\delta+\bar{\rho}} v^\mm{d}(t)\|_{0}^2 + {\Lambda} (\mu_1\|\nabla v^{\mathrm{d}}\|^2_0
 + \mu_2 \| \mm{div}v^{\mathrm{d}}\|^2_0)+ cf(\epsilon, \delta, t).
\end{aligned}\end{equation}

Recalling that $v^\mm{d}\in C^0([0,T^{\max}),H^3)$ and $v^{\mm{d}}(0)=\delta^2 v_{\mm{r}}$, we apply Newton-Leibniz's formula and Cauchy-Schwarz's inequality to find that
 \begin{equation}\begin{aligned}  \label{0316}
& \Lambda( \mu_1\|\nabla v^\mm{d}(t)\|_{0}^2  + \mu_2 \| \mm{div}v^{\mathrm{d}}(t)\|^2_0) \\
& =  2\Lambda\mu_1\int_0^t\int \nabla  v^\mm{d} :\nabla  v_\tau^\mm{d} \mm{d} x\mathrm{d}\tau
+ 2\Lambda\mu_2\int_0^t\int \mm{div} v^\mm{d}\mm{div} v_\tau^\mm{d} \mm{d} x\mathrm{d}\tau  \\
&\quad +\delta^4\Lambda( \mu_1\|\nabla v_\mm{r}\|_{0}^2  + \mu_2 \| \mm{div}v_\mm{r}\|^2_0) \\
& \leq \Lambda^2\int_0^t(\mu_1\|\nabla v^\mm{d} \|_{0}^2 +\mu_2\|\mm{div} v^\mm{d}  \|_{0}^2)\mathrm{d}\tau +\int_0^t(
\mu_1\|\nabla v_\tau^\mm{d} \|_{0}^2+\mu_2\|\mm{div} v_\tau^\mm{d} \|_{0}^2) \mathrm{d}\tau.
 \end{aligned}\end{equation}
 Combining \eqref{new0311} with \eqref{0316}, one gets
 \begin{equation}\label{inequalemee}\begin{aligned}
&\frac{1}{\Lambda}\|\sqrt{\varrho^\delta+\bar{\rho}}  v_t^{\mathrm{d}}(t)\|^2_{0 }+
({\mu_1}\|\nabla v^{\mathrm{d}}(t)\|_{0 }^2+\mu_2\|\mm{div} v^{\mathrm{d}}(t)\|_{0 }^2)\\
& \leq   {\Lambda}\|\sqrt{\varrho^\delta+\bar{\rho}} v^{\mathrm{d}}(t)\|^2_{0 }+2 {\Lambda}\int_0^t
(\mu_1\|
\nabla v^{\mathrm{d}}\|_{0 }^2+
\mu_2\|\mm{div} v^{\mathrm{d}}\|_{0 }^2\mm{d}\tau)\mm{d}\tau + cf(\epsilon, \delta, t).
\end{aligned}\end{equation}
 On the other hand,
\begin{equation*}\begin{aligned}\label{}
\frac{\mm{d}}{\mm{d}t}\|\sqrt{\varrho^\delta+\bar{\rho}} v^\mm{d} \|^2_{0}=&2\int
(\varrho^\delta+\bar{\rho}) v^\mm{d} \cdot  v^\mm{d}_t \mm{d} x+\int
\varrho^\delta_t |v^\mm{d}|^2 \mm{d} x\\
\leq&\frac{1}{\Lambda}\|\sqrt{(\varrho^\delta+\bar{\rho})}  v_t^\mm{d} \|^2_{0}
+\Lambda\|\sqrt{\varrho^\delta+ \bar{\rho}} v^\mm{d} \|^2_{0}+\int
\varrho^\delta_t |v^\mm{d}|^2 \mm{d} x
\end{aligned}\end{equation*}
and
$$ \int\varrho^\delta_t |v^\mm{d}|^2 \mm{d} x = -\int\mm{div}(\bar{\rho}v^{\delta}+\varrho^\delta v^\delta) |v^\mm{d}|^2 \mm{d} x
\lesssim (\|v^\delta\|_1+\|\varrho^\delta\|_2\|v^\delta\|_1 )\|v^{\mm{d}}\|_1^2 \lesssim \delta^3 e^{3\Lambda t}. $$

If we put the previous three estimates together, we get the differential inequality
\begin{equation}\label{growallsinequa}
\begin{aligned}
& \frac{\mm{d}}{\mm{d}t} \|\sqrt{\varrho^\delta+\bar{\rho}} v^\mm{d}(t)\|^2_{0}+
(\mu_1\| \nabla v^\mm{d}(t)\|_{0}^2+\mu_2\| \mm{div} v^\mm{d}(t)\|_{0}^2) \\
&\leq 2\Lambda\left( \|\sqrt{\varrho^\delta+\bar{\rho}} v^\mm{d}(t)\|^2_{0}
 +\int_0^t (\mu_1\|\nabla v^\mm{d}\|_{0}^2 +\mu_2\|\mm{div} v^\mm{d}\|_{0}^2)\mathrm{d}\tau\right) +cf(\epsilon,\delta,t).
\end{aligned}
\end{equation}
Recalling $v^{\mm{d}}(0)=\delta^2v_{\mm{r}}$, one can apply Gronwall's inequality to \eqref{growallsinequa} to conclude
 \begin{equation}\label{estimerrvelcoity}
\begin{aligned}
& \|\sqrt{\varrho^\delta+\bar{\rho}} v^{\mathrm{d}}(t)\|^2_{0}+ {\mu}_1\int_0^t\|\nabla v^{\mathrm{d}}\|^2_{0}
\mm{d}\tau + {\mu}_2\int_0^t\|\mm{div}v^{\mathrm{d}}\|^2_{0} \mm{d}\tau \\
& \quad \leq   e^{2\Lambda t}\left(\int_0^t cf(\epsilon, \delta, \tau)e^{-2\Lambda\tau}\mm{d}\tau
+\delta^4\|\sqrt{\varrho^\delta_0+\bar{\rho}} v^{\mathrm{r}}\|^2_{0}\right) \lesssim f(\epsilon,\delta,t), \qquad\forall\, t\in (0,T^*).
 \end{aligned}  \end{equation}
Moreover, using \eqref{insoslai} and \eqref{infer201604061602}, we can further deduce from \eqref{estimerrvelcoity}, \eqref{inequalemee} and \eqref{new0311} that
\begin{eqnarray}\label{uestimate1n}
\|v^{\mathrm{d}}(t)\|_{1 }^2+\| v_t^{\mathrm{d}}(t)\|^2_0 +\|v^{\mathrm{d}}_\tau\|^2_{L^2((0,t),H^1)}\lesssim f(\epsilon,\delta,t).
\end{eqnarray}

Next we turn to the derivation of the error estimates for the perturbation density and magnetic field.
It follows from the equations \eqref{newfor1232}$_1$ that
\begin{equation*}\begin{aligned}
\frac{1}{2}\frac{\mm{d}}{\mm{d}t}\|\varrho^{\mathrm{d}}(t)\|_{0}^2= & - \int
 \mm{div}(\bar{\rho}v^{\mathrm{d}} +\varrho^{\delta}v^{\delta}) )\varrho^{\mathrm{d}}\mm{d}x\\
= & -\int (\mm{div}(\bar{\rho}v^{\mathrm{d}}  +\varrho^{\mm{a}}v^{\delta} )+ \varrho^{\mm{d}}\mm{div}v^{\delta}/2 )\varrho^{\mathrm{d}}\mm{d}x\\
\lesssim &(\| v^{\mathrm{d}}\|_1 + \|\varrho^{\mm{a}}\|_2\| v^{\delta}\|_1 + \|\varrho^{\mm{d}}\|_{L^\infty}\|v^{\delta}\|_1 )\|\varrho^{\mathrm{d}}\|_{0}.
\end{aligned}\end{equation*}
Therefore, from \eqref{uestimate1n}, \eqref{201602091446MH1}, \eqref{20160207MH} and \eqref{201602071MH}, it follows that
\begin{equation}\begin{aligned}\label{erroresimts}
 \|\varrho^{\mathrm{d}}(t)\|_{0 }\lesssim  &  \int_0^t\left(\| v^{\mathrm{d}}\|_{1}+ \| \varrho^{\mm{a}}\|_{2}\| v^{\delta}\|_{1}+
 \|  \varrho^{\mm{d}} \|_{0}^{\frac{1}{2}}\|  \varrho^{\mm{d}} \|_3^{\frac{1}{2}}\|v^{\delta} \|_{1}\right)\mm{d}\tau
\lesssim  \sqrt{  f(\epsilon, \delta, t)}.
\end{aligned}\end{equation}
Using \eqref{newfor1232}$_1$ again, we can argue analogously to \eqref{201604112016n19212000n} to deduce
\begin{equation}  \label{201602071919MH}
\begin{aligned}\| \varrho_t^{\mathrm{d}}\|^2_0 = \| \mm{div}(\bar{\rho}v^{\mathrm{d}}+ \varrho^{\delta}v^{\delta})\|^2_0
 \lesssim  f(\epsilon, \delta, t).
 \end{aligned}   \end{equation}
Similarly to \eqref{erroresimts} and \eqref{201602071919MH}, we get from \eqref{newfor1232}$_3$ that
\begin{equation}
\label{201602071933MH}
\begin{aligned}\|N^{\mathrm{d}}(t)\|_{0 }+\| {N}_t^{\mathrm{d}}   \|_0
\lesssim
   \sqrt{ f(\epsilon, \delta, t)}.
\end{aligned}
\end{equation}

Now, \eqref{0508} follows from the estimates \eqref{uestimate1n}--\eqref{201602071933MH} immediately.
This completes the proof of Lemma \ref{lem:0401}.  \hfill$\Box$
\end{pf}

Finally, let
\begin{equation}\label{201602092112MH}
T^\delta:=\frac{1}{\Lambda}\ln \frac{\epsilon}{\delta},\quad\mbox{i.e.,}\quad \epsilon=\delta e^{\Lambda T^\delta},
\end{equation}
where
 \begin{equation}\label{defined}
\epsilon:=\min\left\{\epsilon_0,\frac{(C_3C_5)^2}{4 C_6^2},\frac{\|\tilde{v}_{0}\|_{0}^2}{16C_6^2},1 \right\}>0,
 \end{equation}
 and $\epsilon_0$ is the same as in Remark \ref{rem:06011418}. Then we have the following conclusion concerning the relation
 between $T^\delta$ and $T^*$.
 %%%%%%%%%%%%%%%%%%%%%%%%%%%%%%%%%%%%%%%%
\begin{lem}   \label{lem:201602091151MH}
Under the assumptions of Lemma \ref{lem:0401}, $T^\delta<T^*$.
\end{lem}
\begin{pf}
It suffices to show that $T^*<\infty$, and we prove it by contradiction.
Suppose that $T^\delta\geq T^*$. By virtue of the definitions of $(\varrho^\mm{a},v^\mm{a},N^\mm{a})$,
$C_5$, $f(\epsilon, \delta,t)$ and $\epsilon$ in \eqref{201602092112MH}, it is easy to verify that
 $$\begin{aligned}
 &\mathcal{S}(\varrho^\mm{a}(t),v^\mm{a}(t),N^\mm{a}(t))+\| v_\tau^\mm{a} \|_{L^2((0,t),H^1)} \\
 &\leq  \left(\sqrt{\| \tilde{v}_0 \|_1^2+\|( \tilde{\varrho}_0 ,  \tilde{N}_0,\Lambda  \tilde{\varrho}_0,\Lambda
\tilde{v}_0 , \Lambda  \tilde{N}_0 ) \|_0^2}+\sqrt{\frac{\Lambda}{2}}\| \tilde{v}_0\|_1\right) \delta e^{\Lambda t}  \\
 &\leq 2\sqrt{(1+\Lambda)}\|( \tilde{\varrho}_0 , \tilde{v}_0, \tilde{N}_0)\|_1\leq   C_5
\delta e^{\Lambda t}\leq C_3C_5  \delta e^{\Lambda t}\end{aligned} $$
 and
\begin{equation}   \label{201602081013MH}
f(\epsilon, \delta,t) \leq  2{\epsilon}  \delta^{2} e^{{2\Lambda t} }\mbox{ for any }t\leq T^*.  \end{equation}
Consequently,
 \begin{equation*}\begin{aligned}
& \mathcal{S} (\varrho^\delta(T^*),v^\delta(T^*), N^\delta(T^*)) +\| v_\tau^\delta \|_{L^2((0,T^*),H^1)} \\
&\leq \mathcal{S}(\varrho^\mm{a}(T^*),v^\mm{a}(T^*), N^\mm{a}(T^*)) +\|v_\tau^\mm{a}\|_{L^2((0,T^*),H^1)}
    +  \mathcal{S} (\varrho^\mm{d}(T^*),v^\mm{d}(T^*),N^\mm{d}(T^*))
   +\|  v_\tau^\mm{d} \|_{L^2((0,T^*),H^1)} \\
&\leq  C_3C_5 \delta e^{\Lambda T^*} + C_6\sqrt{f(\epsilon,\delta,T^*)} <\delta e^{\Lambda T^*}
(C_3C_5 + 2C_6 \sqrt{\epsilon} ) <2 C_3C_5 \delta e^{\Lambda T^*},
 \end{aligned} \end{equation*}
which contradicts \eqref{201602092200MH}. Hence, $T^\delta<T^*$. The proof is complete. \hfill $\Box$
\end{pf}

Now we are in a position to show Theorem \ref{thm:0101jfjsww101315}. Let $\epsilon$ be given by \eqref{defined}, and
$$\varepsilon :=  \epsilon\min\{1,\|\tilde{v}_{0}\|_{0}/2\},$$
then $\varepsilon>0$ by \eqref{201602081445MH} and the fact $\tilde{v}_{0}=\tilde{u}$.
For any given $\delta\in (0,\iota)$,

(1) if the solution of the original perturbation equations
$(\varrho^{\delta},{v}^{\delta},{ {N}}^{\delta})$ satisfies
\begin{equation}\label{201602081002MH}
\|(\varrho^{\delta}, {v}^{\delta},N^{\delta})(T)\|_3> \epsilon\quad \mbox{ for some }T\in (0,T^{\max}),
 \end{equation}
 then Theorem \ref{thm:0101jfjsww101315} automatically holds by virtue of $\epsilon\geq\varepsilon$.

(2) if \eqref{201602081002MH} fails, then one has $T^{\max}=\infty$ by Proposition \ref{pro:0401n20160203} and Remark \ref{rem:06011418}, and
 $$\|(\varrho^{\delta},{v}^{\delta},{ {N}}^{\delta})(t)\|_3\leq \epsilon\quad\mbox{ for any }t>0. $$
Thus, we can use \eqref{201602092112MH}--\eqref{defined} and Lemmas \ref{lem:0401}--\ref{lem:201602091151MH} to deduce that
 \begin{equation*}\begin{aligned}
 \|v^{\delta}(T^\delta)\|_{0}\geq & \|v^{\mathrm{a}}(T^{\delta})\|_{0}-\|v^{\mm{d}}(T^{\delta})\|_{0 }
= \delta\|v^{\mathrm{l}}(T^{\delta})\|_{0 }-\|v^{\mm{d}}(T^{\delta})\|_{0} \\
  >&\|\tilde{v}_{0}\|_{0} \delta e^{\Lambda T^\delta}- 2C_6\sqrt{\epsilon}\delta e^{\Lambda T^\delta}  \\
 = & (\|\tilde{v}_{0}\|_{0} -2C_6\sqrt{\epsilon})\epsilon >\epsilon\|\bar{v}_{0}\|_{0} /2\geq \varepsilon,
 \end{aligned}      \end{equation*}
 which proves Theorem \ref{thm:0101jfjsww101315}.  \hfill $\Box$

\section{Positivity of the energy functional}\label{201605171506}

This section is devoted to verification of $\Xi>0$ under the assumptions of Theorem \ref{thm:201605151841}--\ref{201605161704n}.
For this purpose, we first (Fourier) transform the energy functional $E$ to a new energy functional with frequency.
It should be pointed out that in the 2D case, $E$ is defined by
$$E (w) := \int g\bar{\rho}'w_2^2\mm{d} x +\int 2 g\bar{\rho}w_2\mm{div}w\dx
 -\int\Big( \gamma \bar{P}|\mm{div}w|^2 +\lambda m^2(|\partial_1 w_2|^2+|\partial_2 w_2|^2) \Big) \dx,$$
 where
 $\bar{\rho}:=\bar{\rho}(x_2)$, $\mm{div}w=\partial_1w_1 +\partial_2 w_2$, $w=(w_1,w_2)$ and $x=(x_1,x_2)$.

 We fix a spatial frequency $\xi=(\xi_1,\xi_2)\in \mathbb{R}^2$, and define the new unknowns
\begin{equation*}
\tilde{u}_1( {x})=-\mathrm{i}\varphi(x_3)e^{\mm{i} {x}_h\cdot\xi}, \;\;
\tilde{u}_2( {x})=-\mathrm{i}\theta(x_3)e^{\mm{i} {x}_h\cdot\xi},\;\;\tilde{u}_3( {x})=\psi(x_3)e^{\mm{i} {x}_h\cdot\xi},
\end{equation*}
where $\mm{i}^2=-1$.  {Substituting the above ansatz into} \eqref{201604061413}, we see that $\varphi$, $\theta$, $\psi$ and
$\Lambda$ satisfy the following system of ODEs:
\begin{eqnarray}\label{2.2.1}
\left\{ \begin{array}{ll}  \Lambda^2\bar{\rho}\varphi+\xi_1 P'(\bar{\rho})\bar{\rho}(\xi_1\varphi+
\xi_2\theta+\psi')+\Lambda\mu_1 (|\xi|^2\varphi-\varphi'')=  \xi_1g\bar{\rho}\psi - \Lambda\mu_2 \xi_1(\xi_1\varphi+\xi_2\theta+\psi') ,\\[2mm]
\Lambda^2\bar{\rho}\theta+\xi_2 P'(\bar{\rho}) \bar{\rho}(\xi_1\varphi+\xi_2\theta+\psi')+\Lambda\mu_1 (|\xi|^2\theta-\theta'') \\
\quad =\xi_2g\bar{\rho}\psi -\Lambda\mu_2 \xi_2(\xi_1 \varphi+\xi_2 \theta+ \psi')-\lambda m^2|\xi|^2\theta
-\lambda m^2\xi_2\psi' ,  \\[2mm]
\Lambda^2\bar{\rho}\psi-(\lambda m^2(\xi_2\theta+\psi')+P'(\bar{\rho})\bar{\rho}(\xi_1\varphi+\xi_2\theta+\psi'))'
+\Lambda\mu_1 (|\xi|^2\psi-\psi'') \\
\quad = g\bar{\rho}(\xi_1\varphi+\xi_2\theta) -\Lambda\mu_2((\xi_1\varphi +\xi_2\theta)'+\psi'')-\lambda m^2|\xi_1|^2\psi
\end{array}   \right.
\end{eqnarray}
with
$$ \varphi(-l)=\theta(-l) = \psi(-l) = \varphi(l)=\theta(l)=\psi(l)=0. $$

We multiply \eqref{2.2.1} with $(\varphi,\theta,\psi)$ in $L^2(-l,l)$ respectively to deduce that
$$ \begin{aligned}
 \Lambda^2\int_{-l}^l\bar{\rho}(\varphi^2 +\theta^2 +\psi^2) \mm{d}x_3=&\tilde{E}(\varphi,\theta,\psi)-\Lambda\mu_2
\int_{-l}^l ( \xi_1\varphi+\xi_2\theta - \psi'  )^2\mm{d}x_3\\
& -\Lambda\mu_1 \int_{-l}^l (|\xi|^2(\varphi^2+\theta^2+\psi^2)+ |\varphi'|^2+ |\theta'|^2+ |\psi'|^2)\mm{d}x_3,
\end{aligned} $$
where
$$ \begin{aligned} \tilde{E}(\varphi,\theta,\psi):=& \int_{-l}^l
 \left(g\bar{\rho}' \psi^2 + 2 g\bar{\rho}\psi (\xi_1\varphi+\xi_2\theta+\psi')\right.  \\
 & \qquad\quad  \left.-\gamma \bar{P}\left(\xi_1\varphi+\xi_2\theta+\psi' \right)^2- \lambda m^2\left(
\xi_1^2(\theta^2+\psi^2)+(\xi_2\theta+\psi')^2\right)\right) \mm{d}x_3,
\end{aligned}$$
which we call the energy functional with frequency. The next lemma gives the positivity of $E$ based on the $\tilde{E}$.
%%%%%%%%%%%%%%%%%%%%%%%%%%%%%%%%%%%%%%
\begin{lem}\label{201605171718}
Under the assumptions of Theorem \ref{thm:201605151841}, there is a function $w_0\in H^1_0$, such that
\begin{equation}  \label{201605161659}  E(w_0)>0.  \end{equation}
\end{lem}
\begin{pf}
(1) We show \eqref{201605161659} under Schwarzschild's condition. To this end, we take $\xi_1=L^{-1}_1$ and $\xi_2=L_2^{-1}$,
then we can rewrite $\tilde{E}$ as follows.
\begin{equation}\label{201605141107}
\begin{aligned} \tilde{E}(\varphi,\theta,\psi):=& \int_{-l}^l
\left(\left(\frac{g^2\bar{\rho}^2}{\gamma  \bar{P}}+g\bar{\rho}' -\lambda \xi_1^2 m^2\right)\psi^2-\frac{\lambda \xi_1^2 m^2}{|\xi|^2}|\psi'|^2\right.\\
&\left.\qquad\quad -\gamma \bar{P}\left(\xi_1\varphi+\xi_2\theta + \psi' - \frac{g\bar{\rho}\psi}{\gamma \bar{P}}\right)^2-\lambda|\xi|^2 m^2
\left( \theta+\frac{\xi_2\psi'}{|\xi|^2}\right)^2\right)\mm{d}x_3.
\end{aligned} \end{equation}

Recalling the definition of $\xi_{3D}$ and the fact $\xi_1<\xi_{3D}$, we see that there is a function $\psi_0\in H_0^1(-l,l)$, such that
\begin{equation*}
 \int_{-l}^l \left(\left( \frac{g^2\bar{\rho}^2}{\gamma \bar{P}}+g\bar{\rho}'-\lambda\xi_1^2m^2\right)\psi^2_0
 -\frac{\lambda \xi_1^2 m^2}{|\xi|^2}|\psi'_0|^2\right)\mm{d}x_3>0.
\end{equation*}
Denoting
$$\theta_0:=-\frac{\xi_2 \psi'_0}{|\xi|^2}\mbox{ and }\varphi_0:=\xi_1^{-1}\left(
\frac{g\bar{\rho}\psi_0}{\gamma
 \bar{P} }-\psi'_0-\xi_2\theta_0\right),$$
\eqref{201605141107} thus reduces to  \begin{equation}
\label{05112209}
\tilde{E}(\varphi_0,\theta_0,\psi_0)=  \int_{-l}^l
\left(\left( \frac{g^2\bar{\rho}^2}{\gamma \bar{P}}+g\bar{\rho}'-\lambda\xi_1^2m^2\right)\psi^2_0-\frac{\lambda \xi_1^2 m^2}{|\xi|^2}|\psi'_0|^2\right)\mm{d}x_3>0.
\end{equation}
We mention that Newcomb regarded \eqref{Newcombinstablity} as an instability condition based the fact \eqref{05112209} \cite{NWACTOPF}.

Now, define
$$ w^0_1:= \varphi_0\sin(\xi\cdot x_h),\quad w^0_2:=\theta_0\sin(\xi\cdot x_h),\quad w^0_3:=\psi_0\cos(\xi\cdot x_h).$$
By \eqref{05112209} and the fact
$$\int_{\mathcal{T}}\sin^2(\xi\cdot x_h)\mm{d}x_h=\int_{\mathcal{T}}\cos^2(\xi\cdot x_h)\mm{d}x_h=2\pi^2 L_1 L_2>0,$$
we have
$${E}(w_0)=\tilde{E}(\varphi_0,\theta_0,\psi_0)\int_{\mathcal{T}}\cos^2(\xi\cdot x_h)\mm{d}x_h>0\qquad (w_0:=(w^0_1,w^0_2,w^0_3)),$$
from which \eqref{201605161659} follows.

(2) We proceed to show \eqref{201605161659} under Tserkovnikov's condition  \eqref{Tserkovnikovs}. Taking $\xi_1=0$ and $\xi_2=1/L_2$, one has
$$\tilde{E}(\varphi,\theta,\psi)= \int_{-l}^l \Big(\big( g\bar{\rho}'+\frac{g^2\bar{\rho}^2}{\gamma \bar{P}+\lambda m^2}\big)\psi^2
-(\gamma  \bar{P}+ \lambda m^2)\big(\xi_2\theta+\psi' - \frac{g\bar{\rho}\psi}{\gamma  \bar{P}+\lambda m^2}\big)^2\Big)\mm{d}x_3.$$
Recalling Tserkovnikov's condition, we find that there is a function $\psi_0\in H_0^1(-l,l)$, such that
$$  \int_{-l}^l  \left(g\bar{\rho}'+\frac{g^2\bar{\rho}^2}{\gamma  \bar{P}+\lambda m^2}\right)\psi^2_0 \mm{d}x_3>0.$$
Denoting
$$\theta_0:=\xi_2^{-1}\left(\frac{g\bar{\rho}\psi_0}{\gamma \bar{P}+\lambda m^2}-\psi'_0\right),$$
we easily see that for any $\varphi$,
\begin{equation} \label{05112209n}  \tilde{E}(\varphi,\theta_0,\psi_0)>0. \end{equation}

Denoting $w_0=(w^0_1,w^0_2,w^0_3):=(1,\theta_0\sin(\xi_2 x_2),\psi_0\cos(\xi_2 x_2))$, we use \eqref{05112209n} and the fact
$$ \int_{\mathcal{T}}\cos^2(\xi_2x_2)\mm{d}x_{h}=2\pi^2 L_1 L_2>0 $$
to find that
$${E}(w_0)=\tilde{E}(\varphi,\theta_0,\psi_0)\int_{\mathcal{T}}\cos^2(\xi_2x_2)\mm{d}x_{h}>0,$$
which proves the lemma.   \hfill$\Box$
 \end{pf}

\begin{lem}\label{201605171716}
Under the assumptions of Theorem \ref{201605161704n}, there is a function $w_0\in H^1$, such that $E(w_0)>0$.
\end{lem}
\begin{pf}
In the 2D case, $\tilde{E}$ is defined by
$$ \begin{aligned} \tilde{E}(\varphi, \psi):=& \int_{-l}^l
 \left(g\bar{\rho}' \psi^2 + 2 g\bar{\rho}\psi (\xi_1\varphi+ \psi') -\gamma \bar{P}\left(\xi_1\varphi +\psi' \right)^2- \lambda m^2\left(
\xi_1^2 \psi^2 + |\psi'|^2\right)\right) \mm{d}x_2 , \end{aligned}$$
which can be rewritten as
\begin{equation*}  \begin{aligned}
\tilde{E}(\varphi,\psi):=& \int_{-l}^l\left(\left(\frac{g^2\bar{\rho}^2}{\gamma \bar{P}} +g\bar{\rho}'-\lambda \xi^2_1 m^2\right)\psi^2
-{\lambda m^2}|\psi'|^2 -\gamma \bar{P}\left(\xi_1\varphi+\psi' - \frac{g\bar{\rho}\psi}{\gamma \bar{P}}\right)^2\right)\mm{d}x_2.
\end{aligned} \end{equation*}

Taking $\xi_1=L^{-1}_1$, and taking into account the definition of $\xi_{2D}$ and the assumption $L^{-1}_1<\xi_{2D}$,
one sees that there is a function $\psi_0$ satisfying
\begin{equation*}
 \int_{-l}^l
\left(\left( \frac{g^2\bar{\rho}^2}{\gamma \bar{P}}+g\bar{\rho}'-\lambda\xi_1^2m^2\right)
\psi^2_0- {\lambda m^2} |\psi'_0|^2\right)\mm{d}x_2>0.
\end{equation*}
Denoting
$$\varphi_0:=\xi^{-1}_1\left(
\frac{g\bar{\rho}\psi_0}{\gamma
 \bar{P} }-\psi'_0\right),$$
one has
\begin{equation*} \tilde{E}(\varphi_0,\psi_0)=  \int_{-l}^l \left(\left( \frac{g^2\bar{\rho}^2}{\gamma  \bar{P}}+g\bar{\rho}'
-\lambda\xi^2_1m^2\right)\psi^2_0-{\lambda m^2}|\psi'_0|^2\right)\mm{d}x_2>0.
\end{equation*}

Now, if we denote $w^0_1:=\varphi_0\sin(\xi_1 x_1)$ and $w^0_2:=\psi_0\cos (\xi_1 x_1)$,
then we have by \eqref{05112209} and the fact
$$\int_{2\pi L_1\mathbb{T}}\sin^2(\xi_1 x_1)\mm{d}x_1=\int_{2\pi L_1\mathbb{T}}\cos^2(\xi_1 x_1)\mm{d}x_1=\pi L_1>0$$
that
$$ {E}(w_0)=\tilde{E}(\varphi_0,\psi_0)\int_{2\pi L_1\mathbb{T}}\cos^2(\xi_1x_1)\mm{d}x_1>0\qquad (w_0:=(w^0_1,w^0_2)), $$
which proves the lemma. \hfill$\Box$
\end{pf}

In view of Lemmas \ref{201605171718} and \ref{201605171716}, we conclude $\Xi>0$ under the assumptions
of Theorem \ref{thm:201605151841}--\ref{201605161704n}.
Therefore, we obtain Theorems \ref{thm:201605151841} and \ref{201605161704n} by following the proof of Theorem \ref{thm:0101jfjsww101315}.

\vspace{4mm} \noindent\textbf{Acknowledgements.}
The research of Fei Jiang   was supported by NSFC (Grant
Nos 11471134 and 11671086)  and the NSF of Fujian Province of China (Grant
No. 2016J06001),
 and the research of Song Jiang by the National Basic Research Program (Grant
No 2014CB745002) %, 2011CB309705)
  and NSFC (Grant
Nos 11371065 and 11631008).
%This work was finished when Fei Jiang was visiting the City University of Hong Kong.
%Fei Jiang would like to thank Prof.   for his kind hospitality.
 %Finally, the authors would like to thank the anonymous referee for invaluable suggestions,
 % which improve the presentation of this paper.

\renewcommand\refname{References}
\renewenvironment{thebibliography}[1]{%
\section*{\refname}
\list{{\arabic{enumi}}}{\def\makelabel##1{\hss{##1}}\topsep=0mm
\parsep=0mm
\partopsep=0mm\itemsep=0mm
\labelsep=1ex\itemindent=0mm
\settowidth\labelwidth{\small[#1]}%
\leftmargin\labelwidth \advance\leftmargin\labelsep
\advance\leftmargin -\itemindent
\usecounter{enumi}}\small
\def\newblock{\ }
\sloppy\clubpenalty4000\widowpenalty4000
\sfcode`\.=1000\relax}{\endlist}
\bibliographystyle{model1b-num-names}
%\bibliography{refs}

\end{document}